\documentclass[leqno,11pt]{article}  %putting formula numbers on the left side
\usepackage{amsmath}
\usepackage{amssymb}

\usepackage[german,english]{babel}
\usepackage{amsthm}
\title{Frobenius and Monodromy operators in rigid analysis, and Drinfel'd's symmetric space}

\author{\textsc{Elmar Grosse-Kl\"onne}}
\date{}

\theoremstyle{plain} 
\newtheorem{satz}{Theorem}[section]  %@@
\newtheorem{kor}[satz]{Corollary}  %@@
\newtheorem{lem}[satz]{Lemma}  %@@
\newtheorem{pro}[satz]{Proposition}  %@@
\newcommand{\spec}{\mbox{\rm Spec}}  %@@
\newcommand{\proj}{\mbox{\rm Proj}}
\newcommand{\quot}{\mbox{\rm Quot}}  %@@
\newcommand{\spm}{\mbox{\rm Sp}}  %@@
\newcommand{\spf}{\mbox{\rm Spf}}  %@@
\newcommand{\spwf}{\mbox{\rm Spwf}}
\newcommand{\bi}{\mbox{\rm im}}  %@@
\newcommand{\ke}{\mbox{\rm Ker}}  %@@
\newcommand{\sym}{\mbox{\rm Sym}}
\newcommand{\id}{\mbox{\rm id}}

\newcommand{\dlog}{\mbox{\rm dlog}}
\newcommand{\Gr}{\mbox{\rm Gr}}

\theoremstyle{remark}

\theoremstyle{definition}

\begin{document}
\maketitle
\footnote[0]
    {2000 \textit{Mathematics Subject Classification}.
    14F30}                               %2000MSC 
\footnote[0]{\textit{Key words and phrases}. Hyodo-Kato isomorphism, Drinfeld space, dagger space, log rigid cohomology}

\footnote[0]{Part of this work was done during my visit at the University of California, Berkeley, in the academic year 2000/2001. I wish to thank the members of this institution, in particular Robert Coleman, for their warm hospitality. I am grateful to Arthur Ogus for discussions on the Hyodo-Kato isomorphism and on the application of crystalline methods to the Hodge-type decomposition conjecture. Thanks also go to Matthias Strauch for a discussion on Drinfel'd's $p$-adic symmetric spaces. I thank the referee for pointing on  several points to be clarified. I am grateful to the Deutsche Forschungs Gemeinschaft for supporting my stay at Berkeley.}

\begin{abstract} 

We define Frobenius and monodromy operators on the de Rham cohomology of $K$-dagger spaces (rigid spaces with overconvergent structure sheaves) with strictly semistable reduction $Y$, over a complete discrete valuation ring $K$ of mixed characteristic. For this we introduce log rigid cohomology and generalize the so called Hyodo-Kato isomorphism to versions for non-proper $Y$, for non-perfect residue fields, for non-integrally defined coefficients, and for the various strata of $Y$. We apply this to define and investigate crystalline structure elements on the de Rham cohomology of Drinfel'd's symmetric space $X$ and its quotients. Our results are used in a critical way in the recent proof of the monodromy-weight conjecture for quotients of $X$ given by de Shalit \cite{desh}.

\end{abstract}

% \addtolength{\textwidth}{1.0in}
% \setlength{\hoffset}{-.5in}
%
%\addtolength{\topmargin}{-34pt}
%\addtolength{\textheight}{68pt}
%\baselineskip13,6pt
%\setlength{\parindent}{0pt}
%\setlength{\parskip}{0.8ex plus 0.2ex minus 0.2ex}
%\setlength{\oddsidemargin}{3cm}
%\frenchspacing
%\sloppy
%\pagestyle{myheadings}
%\markboth{}{}

%\textwidth14cm
%\textheight9in
%\hoffset+1in
%\voffset-0,8cm

\begin{center} {\bf Introduction} \end{center}

Let $A$ be a complete discrete valuation ring of mixed characteristic $(0,p)$, with perfect residue field $k$ and quotient field $K$, let $A_0=W(k)$ and $K_0=\quot(A_0)$. Let $X$ be a proper strictly semistable $A$-scheme. The Hyodo-Kato isomorphism is an isomorphism $\rho$ (depending on the choice of a uniformizer $\pi\in A$) between the de Rham cohomology $H_{dR}^*(X_K)$ of the generic fibre $X_K$ of $X$ and the (scalar extended) Hyodo-Kato cohomology $H^*_{HK}(Y)\otimes_{K_0}K$ of the special fibre $X_k=Y$ of $X$ endowed with its canonical log structure. It plays an important role in the Fontaine-Jannsen conjecture $C_{st}$, now proven (independently) by Tsuji and Faltings: $\rho$ provides $H_{dR}^*(X_K)$ with the structure of a filtered $(\phi,N)$-module in the sense of Fontaine, taking $H^*_{HK}(Y)$ as $K_0$-lattice with $(\phi,N)$-structure. $C_{st}$ says (in particular) that the $p$-adic \'{e}tale cohomology group $H_{et}^*(X_{\overline{K}},\mathbb{Q}_p)$ together with its Gal$(\overline{K}/K)$-action can be reconstructed from this filtered $(\phi,N)$-module.

In this paper we ask for $K_0$-lattices with $(\phi,N)$-structure in the de Rham cohomology of $K$-rigid (or dagger) spaces {\it not necessarily proper}. That $p$-adic Hodge theory should encompass more general $K$-rigid spaces than just smooth proper $K$-schemes was already suggested in Tate's article \cite{tate} and is strongly evidenced by the book \cite{rz} of Rapoport and Zink on $p$-adic period domains. From the paper \cite{berko} of Berkovich it became clear that the study of general rigid spaces can often be reduced to those having strictly semistable reduction. Since de Rham cohomology of rigid spaces should be defined using overconvergent functions we work with weak formal schemes and dagger spaces (\cite{crelle}) rather than formal schemes and rigid spaces. Thus, we start with a strictly semistable weak formal $A$-scheme ${\mathfrak X}$ with associated $K$-dagger space ${\mathfrak X}_{\mathbb{Q}}$ and reduction $Y$. We allow coefficients: local systems $F$ of $K$-vector spaces on ${\mathfrak X}_{\mathbb{Q}}$ arising from representations of $\Pi_1^{top}({\mathfrak X}_{\mathbb{Q}})$, the topological fundamental group of (the Berkovich analytic space associated with) ${\mathfrak X}_{\mathbb{Q}}$ (these $F$ need not be integrally defined, i.e. need not come from crystals on $Y$). We do not require that $k$ be perfect. Yet another new aspect is that besides for $(K_0,\phi,N)$-structures on $H^*_{dR}({\mathfrak X}_{\mathbb{Q}},F)$ alone we ask for such structures on the entire {\it canonical Cech spectral sequence}\begin{gather}E_{1}^{rs}=F(]Y^{r+1}[_{{\mathfrak X}})\otimes_KH^s_{dR}(]Y^{r+1}[_{{\mathfrak X}})\Longrightarrow H^{s+r}_{dR}({\mathfrak X}_{\mathbb{Q}},F).\tag{$*$}\end{gather}
Here $Y^t$ denotes the $t$-fold intersections of irreducible $Y$-components, and $]Z[_{\mathfrak X}$ for a subscheme $Z$ of $Y$ is the preimage of $Z$ under the specialization map $sp:{\mathfrak X}_{\mathbb{Q}}\to Y$. To motivate this we mention that Coleman and Iovita use $(*)$ to describe rigid analytically the Hyodo-Kato monodromy operator $N$ on $H^{1}_{dR}({\mathfrak X}_{\mathbb{Q}})$ for proper ${\mathfrak X}$ of relative dimension $d=1$; in \cite{mono} we ask for the interaction of $(*)$ with Frobenius and monodromy if $d>1$. Also in work of de Shalit, $(*)$ is of central interest, see below.\\We define log rigid cohomology of $E=sp_*F$ as the appropriate substitute for log crystalline cohomology adapted to our purposes. Subschemes $Z$ of $Y$ are endowed with their induced structure of log scheme over the log point $S^0=(\spec(k),1\mapsto 0)$. The interesting thickenings of $S^0$ are the bases ${\mathfrak S}^{0} =(\spf(A_0),1\mapsto 0)$ and ${\mathfrak S}^{\pi} =(\spf(A),1\mapsto\pi)$: while ${\mathbb{R}}\Gamma_{rig}(Z/{\mathfrak S}^{0},E)$ is a $(K_0,N)$-structure, and also a $\phi$-structure if $E$ carries a Frobenius structure, only ${\mathbb{R}}\Gamma_{rig}(Z/{\mathfrak S}^{\pi},E)$ can {\it a priori} be canonically identified with $F(]Z[_{{\mathfrak X}})\otimes_KH^s_{dR}(]Z[_{{\mathfrak X}})$. We have a spectral sequence\begin{gather}E_{1}^{rs}=E(Y^{r+1})\otimes_{K_0}H^s_{rig}(Y^{r+1}/{\mathfrak S}^{0})\Longrightarrow H^{s+r}_{rig}(Y/{\mathfrak S}^{0},E).\tag{$**$}\end{gather}

\begin{satz}\label{lochyka}\label{strahk} (a) (Theorem \ref{hkrigis}) Let $M$ be the intersection of some irreducible components of $Y$. There is an isomorphism (depending on $\pi$)$$\rho_M:{\mathbb{R}}\Gamma_{dR}(]M[_{\mathfrak X})={\mathbb{R}}\Gamma_{rig}(M/{\mathfrak S}^{\pi})\cong {\mathbb{R}}\Gamma_{rig}(M/{\mathfrak S}^{0} )\otimes_{K_0}K.$$(b) (Theorem \ref{ssiso}, Corollary \ref{rigs}) There is an isomorphism (depending on $\pi$)$${\mathbb{R}}\Gamma_{dR}({\mathfrak X}_{\mathbb{Q}},F)={\mathbb{R}}\Gamma_{rig}(Y/{\mathfrak S}^{\pi},E)\cong {\mathbb{R}}\Gamma_{rig}(Y/{\mathfrak S}^{0},E ).$$It comes along with an isomorphism of spectral sequences $(*)\cong(**)$.
\end{satz}

Hence a $(K_0,N)$-structure (resp. $(K_0,\phi,N)$-structure) on $(*)$ since there is such a structure on $(**)$. Except for abelian varieties the existence of the resulting monodromy operator on $H^*_{dR}({\mathfrak X}_{\mathbb{Q}})$ seemed to be unknown before (for non-perfect $k$) even for proper ${\mathfrak X}$. If ${\mathfrak X}$ is proper we have $H^*_{rig}(Y/{\mathfrak S}^{0})=H^*_{HK}(Y)$ and ${\mathbb{R}}\Gamma_{dR}({\mathfrak X}_{\mathbb{Q}})$ coincides with the de Rham cohomology of the rigid space associated to the dagger space ${\mathfrak X}_{\mathbb{Q}}$, and with the de Rham cohomology of an underlying $K$-scheme if such an algebraization exists; we recover the Hyodo-Kato isomorphism (for perfect $k$).

According to the literature, the Hyodo-Kato isomorphism was considered to be a delicate convergence theorem on Frobenius in the log crystalline cohomology of certain log smooth and {\it proper} varieties over perfect fields, see Hyodo-Kato \cite{hyoka} and Ogus \cite{oglog}. On the other hand, the corresponding complex analytic comparison isomorphism relies on a comparison with singular cohomology, see Steenbrink \cite{steen}. Both methods break down in our setting (even if $F=K$, $Y$ is proper and $k$ is perfect since for example $M$ is not log smooth over $S^0$; it is ideally smooth in the sense of Ogus and $H_{crys}^*(M/(\spf(W(k)),1\mapsto0))$ is well behaved, but $H_{crys}^*(M/(\spf(W(k)),1\mapsto p))$ is at present not well understood). We introduce an entirely new method, in fact a {\it geometric} approach. The most conceptual way to describe it is in terms of log schemes with boundary (see \cite{colo}): An $S$-log scheme with boundary $(P,V)$ is an $S$-log scheme $V$ together with a log schematically dense open immersion $V\to P$. Let $S=(\spec(k[t],1\mapsto t)$ and ${\mathfrak S}=(\spwf(A_0[t]^{\dagger}),1\mapsto t)$ (with $(.)^{\dagger}$ denoting weak completion). To the $S^0$-log scheme $M$ we assign finitely many $S$-log schemes with boundary $(P^{J'}_M,V^{J'}_M)$ (here $J'$ is an index), and to the $S^0$-log scheme $Y$ we assign {\it canonically} a simplicial $S$-log scheme with boundary $(P_{\bullet},V_{\bullet})$. For these constructions Falting's interpretation of log structures through line bundles is essential. We define their log rigid cohomology ${\mathbb{R}}\Gamma_{rig}((P^{J'}_M,V^{J'}_M)/{\mathfrak S})$ and ${\mathbb{R}}\Gamma_{rig}((P_{\bullet},V_{\bullet})/{\mathfrak S},E)$ relative to ${\mathfrak S}$. We have natural restriction maps$${\mathbb{R}}\Gamma_{rig}(M/{\mathfrak S}^{0})\otimes_{K_0}K\leftarrow{\mathbb{R}}\Gamma_{rig}((P^{J'}_M,V^{J'}_M)/{\mathfrak S})\otimes_{K_0}K\rightarrow{\mathbb{R}}\Gamma_{rig}(M/{\mathfrak S}^{\pi})$$
$${\mathbb{R}}\Gamma_{rig}(Y/{\mathfrak S}^{0},E)\leftarrow{\mathbb{R}}\Gamma_{rig}((P_{\bullet},V_{\bullet})/{\mathfrak S},E)\rightarrow{\mathbb{R}}\Gamma_{rig}(Y/{\mathfrak S}^{\pi},E)$$and we prove that they are all isomorphisms; once $(P^{J'}_M,V^{J'}_M)$ and $(P_{\bullet},V_{\bullet})$ are found this is by more or less standard local arguments.

As an illustration of how the full strength of the above generalizations of the Hyodo-Kato isomorphism can be applied we consider, for $K$ a finite extension of $\mathbb{Q}_p$, Drinfel'd's $p$-adic symmetric space $X=\Omega^{(d+1)}_K$ of dimension $d$  over $K$ ---  the complement in ${\bf P}^d_K$ of the union of all $K$-rational hyperplanes --- and its strictly semistable weak formal model ${\mathfrak Q}$, which is not proper. The cohomology of $X$ is of great representation theoretical importance. In \cite{ss}, P. Schneider and U. Stuhler computed it as a ${\rm {\rm PGL}}_{d+1}(K)$-representation for an arbitrary cohomology theory satisfying certain minimal axioms. Examples are de Rham cohomology and $\ell$-adic ($\ell\ne p$) cohomology. Moreover they showed that $Gal(\overline{K}/K)$ acts on $H^s_{et}(\overline{X},{\mathbb Q}_{\ell})$ through the $s$-th power of the cyclotomic character. Using Theorem \ref{lochyka} we obtain a $q$-th power Frobenius endomorphism $\phi$ on $H_{dR}^*(X)$, where $|k|=q$. We show $\phi=q^s$ on $H_{dR}^s(X)$  (Corollary \ref{sympur}). The proof relies on a recent acyclicity theorem of E. de Shalit \cite{ds}. We also investigate the $(\phi,N)$-structure on $H_{dR}^{d}(X_{\Gamma},F)$ for quotients $X_{\Gamma}=\Gamma\backslash X$ of $X$ by (sufficiently small) discrete cocompact subgroups $\Gamma$ of ${\rm PGL}_{d+1}(K)$ and coefficients $F$ defined by finite dimensional $K[\Gamma]$-modules ${\bf F}$ (note $\Gamma=\Pi^{top}_1(X_{\Gamma})$). There is a covering spectral sequence\begin{gather} E_2^{rs}=H^r(\Gamma,{\bf F}\otimes_K H_{dR}^s(X))\Longrightarrow H_{dR}^{r+s}(X_{\Gamma},F).\tag*{$(G)_{\pi}$}\end{gather}Only the cohomology $H_{dR}^{d}(X_{\Gamma},F)$ in middle degree $d$ is interesting (\cite{schn}). Let $(F_{\Gamma}^r)_{r\ge 0}$ be the filtration which $(G)_{\pi}$ defines on $H_{dR}^{d}(X_{\Gamma},F)$. By the above it must be the slope and the weight filtration for the Frobenius endomorphism $\phi$ on $H_{dR}^{d}(X_{\Gamma},F)$. Moreover we show that it also coincides with the filtration defined by $(*)$ (with ${\mathfrak X}_{\mathbb{Q}}=X_{\Gamma}$ there).

Schneider and Stuhler conjectured that $(F_{\Gamma}^r)_{r\ge 0}$ is opposite to the Hodge filtration $(F_{Hdg}^j)_{j\ge 0}$. Our results allow us to reconsider this conjecture in terms of $p$-adic Hodge theory. For example we see that it makes a prediction on the $p$-adic \'{e}tale cohomology, as a $Gal(\overline{K}/K)$-representation, of the relative Tate module of the universal $p$-divisible group over $X\otimes_K\widehat{K}^{ur}$, cf. page \pageref{pdiv}. For $d=1$, where the conjecture is known for ${\bf F}$ which are "algebraic", this was exploited in recent work of Iovita and Spiess. Conversely we give a condition on the filtered $(\phi,N)$-module $H_{dR}^{d}(X_{\Gamma},F)$ which implies the Hodge-type decomposition. In particular we reprove it (Theorem \ref{conjhod}) in the case ${\bf F}=K$ (for an earlier and completely different proof see Iovita and Spiess \cite{iovspi}.

In case ${\bf F}=K$, inserting the computations of the vector space dimensions of the graded pieces for $(F_{\Gamma}^r)_{r\ge 0}$ given in \cite{schn}, we also describe the matrix of the monodromy operator $N$ on $H_{dR}^{d}(X_{\Gamma})$:  If $d$ is odd we let $\overline{H}_{dR}^d(X_{\Gamma})=H_{dR}^d(X_{\Gamma})$, if $d$ is even we let $\overline{H}_{dR}^d(X_{\Gamma})$ be the quotient of $H_{dR}^d(X_{\Gamma})$ by a certain one dimensional subspace. All structure elements pass from $H_{dR}^d(X_{\Gamma})$ to $\overline{H}_{dR}^d(X_{\Gamma})$. Then (see \ref{netcoh}):

\begin{satz} (a) The filtrations $({\overline F}_{\Gamma}^r)_{r\ge0}$ and $({\overline F}_{Hdg}^j)_{j\ge 0}$ are opposite.\\(b) The filtration $({\overline F}_{\Gamma}^r)_{r\ge0}$ is stable for $\overline{\phi}$; we have $\overline{\phi}=q^{d-r}$ on ${\overline F}_{\Gamma}^r/\overline {F}^{r+1}_{\Gamma}$.\\(c) The filtration $({\overline F}_{\Gamma}^r)_{r\ge0}$ coincides with both the kernel and the image filtration for $\overline{N}$: for all $r$ we have $$\overline {F}^r_{\Gamma}=\ke(\overline{N}^{d+1-r})=\bi(\overline{N}^{r}).$$  
\end{satz}

For statement (c) we make use of the monodromy-weight conjecture for $X_{\Gamma}$ which has recently been proven by T. Ito \cite{ito} and independently by E. de Shalit \cite{desh}. Ito reduces the problem to standard cohomological conjectures which he then proves for the particular varieties in question. The completely different approach of de Shalit to $X$ and $X_{\Gamma}$, developed in a series of articles, is by concepts of harmonic analysis (also the spectral sequence $(*)$ (for ${\bf F}=K$) appears). It has the advantage of being rather explicit and seems to have some potential for more general semistable reduction situations. de Shalit's proof of the monodromy-weight conjecture relies crucially on our Theorem \ref{lochyka}. Note that our work also responds to the problems 6.14 and 6.15 posed in Ito's paper \cite{ito}.

Coming back to Theorem \ref{strahk} we now sketch how the isomorphisms $\rho_M$ are constructed in two particular cases, namely when $M$ is a single irreducible component of $Y$ (case 1), or when $\dim(Y)=1$ and $M$ is the intersection of two distinct irreducible components of $Y$ (case 2). These two cases require two different ideas (which for more general $M$ must be merged).\\Case 1: Let $M^{\heartsuit}$ be the open subscheme of $M$ complementary to all other irreducible components. Then it turns out that the restriction maps$${\mathbb{R}}\Gamma_{rig}(M/{\mathfrak S}^{0} )\longrightarrow {\mathbb{R}}\Gamma_{rig}(M^{\heartsuit}/{\mathfrak S}^{0} )$$$${\mathbb{R}}\Gamma_{rig}(M/{\mathfrak S}^{\pi} )\longrightarrow {\mathbb{R}}\Gamma_{rig}(M^{\heartsuit}/{\mathfrak S}^{\pi} )$$are isomorphisms: It is at this point where one must work in the overconvergent setting (weak formal/dagger), as opposed to the convergent setting (formal/rigid). Now ${\mathbb{R}}\Gamma_{rig}(M^{\heartsuit}/{\mathfrak S}^{0} )$ resp. ${\mathbb{R}}\Gamma_{rig}(M^{\heartsuit}/{\mathfrak S}^{\pi} )$ are in fact canonically isomorphic to the classical (non logarithmic) rigid cohomology ${\mathbb{R}}\Gamma_{rig}(M^{\heartsuit}/K_0)$ resp. ${\mathbb{R}}\Gamma_{rig}(M^{\heartsuit}/K)$. We take $\rho_M$ to be the composition of isomorphisms$${\mathbb{R}}\Gamma_{rig}(M/{\mathfrak S}^{0} )\otimes_{K_0}K\cong {\mathbb{R}}\Gamma_{rig}(M^{\heartsuit}/K_0)\otimes_{K_0}K\cong {\mathbb{R}}\Gamma_{rig}(M^{\heartsuit}/K)\cong {\mathbb{R}}\Gamma_{rig}(M/{\mathfrak S}^{\pi} ).$$Case 2: For simplicity we assume $M$ is $k$-rational and consists of a single point (its underlying scheme is $\spec(k)$). Let $U$ be the two dimensional open unit polydisk over $K$ with coordinates $x_1, x_2$, viewed as a dagger- (or rigid) analytic space. Let $U^0$ be its closed subspace where $x_1x_2=0$, let $U^{\pi}$ be its closed subspace where $x_1x_2=\pi$. Let ${{\Omega}'}^{\bullet}_U$ be the de Rham complex on $U$ with logarithmic poles along the divisor $U^0$, and let ${\Omega}^{\bullet}_U$ be the quotient of ${{\Omega}'}^{\bullet}_U$ by its sub-${\cal O}_{U}$-algebra generated by $\dlog(x_1x_2)$. It restricts to complexes $\Omega^{\bullet}_{U^0}$ on $U^0$ and $\Omega^{\bullet}_{U^{\pi}}$ on $U^{\pi}$. Of course, $\Omega^{\bullet}_{U^{\pi}}$ is nothing but the classical de Rham complex on the classically smooth space $U^{\pi}$. We may view the $S^0$-log scheme $M$ as exact closed log subscheme of the log smooth $\spf(A_0[x_1,x_2]/(x_1x_2))$ over ${\mathfrak S}^{0} $, resp. the log smooth $\spf(A[x_1,x_2]/(x_1x_2-\pi))$ over ${\mathfrak S}^{\pi} $, with tube $U^0$ resp. $U^{\pi}$. Thus $${\mathbb{R}}\Gamma_{rig}(M/{\mathfrak S}^{0} )\otimes_{K_0}K={\mathbb{R}}\Gamma(U^0,\Omega^{\bullet}_{U^0})$$$${\mathbb{R}}\Gamma_{rig}(M/{\mathfrak S}^{\pi} )={\mathbb{R}}\Gamma(U^{\pi},\Omega^{\bullet}_{U^{\pi}}).$$Now $H^*(U^0,\Omega^{\bullet}_{U^0})\cong H^*(U^{\pi},\Omega^{\bullet}_{U^{\pi}})$ for all $*$ by explicit computations; for example $H^1(U^0,\Omega^{\bullet}_{U^0})$ and $H^1(U^{\pi},\Omega^{\bullet}_{U^{\pi}})$ are one dimensional $K$-vector spaces generated by the class of $\dlog(x_1)$. But we need a construction on the level of complexes. Let $$P=({\bf P}^1_K\times{\bf P}^1_K)^{an}=((\spec(K[x_1])\cup\{\infty\})\times(\spec(K[x_2])\cup\{\infty\}))^{an}$$and let ${{\Omega}'}_P^{\bullet}$ be the de Rham complex on $P$ with logarithmic poles along the divisor $$(\{0\}\times {\bf P}_K^1)\cup ({\bf P}_K^1\times\{0\})\cup (\{\infty\}\times {\bf P}_K^1)\cup ({\bf P}_K^1\times\{\infty\}).$$The section $\dlog(x_1x_2)\in{{\Omega}'}_U^{\bullet}(U)={{\Omega}'}_P^{\bullet}(U)$ extends canonically to a section $\dlog(x_1x_2)\in{{\Omega}'}_P^{\bullet}(P)$ and we let $\Omega_P^{\bullet}$ be the quotient of ${{\Omega}'}^{\bullet}_P$ by its sub-${\cal O}_{P}$-algebra generated by $\dlog(x_1x_2)$. It turns out that the canonical restriction maps $${\mathbb{R}}\Gamma(U^0,\Omega^{\bullet}_{U^0})\leftarrow {\mathbb{R}}\Gamma(P,\Omega^{\bullet}_{P})\longrightarrow {\mathbb{R}}\Gamma(U^{\pi},\Omega^{\bullet}_{U^{\pi}})$$
are both isomorphisms. This yields the wanted $\rho_M$ also in this case.\\  
The paper is organized as follows. In section 1 we recall log schemes with boundary and define log rigid cohomology. In section 2 we introduce $Y$ and certain canonical projective bundles on its component intersections $M$, giving rise to the various $(P_M^{J'},V_M^{J'})$. Liftings of these bundles to characteristic zero are meant to play the role of $P$ in case 2 of the above example. In section 3 the $(P_M^{J'},V_M^{J'})$ are glued into a simplicial $S$-log scheme with boundary. Theorem \ref{strahk} is stated, and \ref{strahk}(a)$\Longrightarrow$\ref{strahk}(b) is explained. In section 4 we prove \ref{strahk}(a) in two steps: first we show that we may restrict everything to the maximal open subscheme $M^{\heartsuit}$ of $M$ which avoids all irreducible components of $Y$ not fully containing $M$ --- this is analogous to case 1 above. Then it follows from explicit computations, the crucial point of the whole method is Lemma \ref{nochiso}. In section 5 we develop the usual weight-monodromy formalism from \cite{steen}, \cite{mokr} for $H^*_{rig}(Y/{\mathfrak S}^{0},E)$, in particular write it as the abutment of a spectral sequence which begins with classical, non-logarithmic rigid cohomology. In sections 6-8 we explain the applications to $\Omega_K^{(d+1)}$.\\

{\it Notations:} For basics on log algebraic geometry we refer to K. Kato \cite{kalo}. For a log scheme $(X,{\cal N}_X\to{\cal O}_X)$ we just write $X$ if it is clear from the context to which log structure on $X$ we refer. In this text, all log schemes and morphisms of log schemes have charts for the {\it Zariski} topology (rather than only for the \'{e}tale topology; since we are interested in {\it strictly} semistable log schemes, see section \ref{bassec}, this class of log schemes will be enough). By an {\it exactification} of a closed immersion of fine log schemes $i':X\to Y'$, we mean a factorization $X\stackrel{i}{\to}Y\stackrel{g}{\to}Y'$ of $i'$ with $i$ an exact closed immersion and $g$ log \'{e}tale. If $i$ has a chart, exactifications exist by the proof of \cite{kalo} 4.10. All this similarly for (weak) formal log schemes.\\
We let $A_0$ be a complete discrete valuation ring with residue field $k$ of characteristic $p>0$, maximal ideal $m_{A_0}$ and fraction field $K_0$ of characteristic zero. If $k$ is perfect, one may take $A_0=W(k)$, the ring of Witt vectors of $k$. Write $\widetilde{W}$ for both the abstract scheme $\spec(A_0)$ and the log scheme $(\spec(A_0),\mbox{trivial log structure})$. Let $A$ be a complete discrete valuation ring which is a totally ramified finite extension of $A_0$. Throughout this paper, {\it we fix a uniformizer $\pi$ of $A$}. Define the log scheme$$\widetilde{S}=(\spec(A_0[t]),(\mathbb{N}\to A_0[t], 1\mapsto t)).$$We will often view $t$ as an element of its log structure. Denote by ${\mathfrak S}$ (resp. ${\mathfrak W}$) the weak completion (see below) of $\widetilde{S}$ (resp. of $\widetilde{W}$). Let $${\mathfrak S}^{0}=(\spf(A_0),(\mathbb{N}\to A_0, 1\mapsto 0)),\quad\quad{\mathfrak S}^{\pi}=(\spf(A),(\mathbb{N}\to A, 1\mapsto \pi)),$$ exact closed (weak) formal log subschemes of ${\mathfrak S}$. Denote by $S$ the exact closed log subscheme of $\widetilde{S}$ (or ${\mathfrak S}$) defined by $(m_{A_0})$, and by $S^0$ the log point over $k$, i.e. the exact closed log subscheme of $S$ defined by $(t)$. Identify $S^0$ with the exact closed log subscheme of ${\mathfrak S}^{a}$ defined by the maximal ideal of ${\cal O}_{{\mathfrak S}^{a}}$ (for $a\in\{0,\pi\}$). Write $K_0(0)=K_0$, $K_0(\pi)=K$. Often (not always) we tried to follow this notational pattern: a roman capital ($Y$, $M$, $V,\ldots$) denotes a $k$-scheme, the same roman capital with superscript \, $\widetilde{}$ \, \,  ($\widetilde{Y}$, $\widetilde{M}$, $\widetilde{V},\ldots$) denotes a scheme flat over $A_0$ together with a fixed embbeding of that $k$-scheme, and the corresponding Fraktur capital ($\mathfrak{Y}$, $\mathfrak{M}$, $\mathfrak{V},\ldots$) denotes the (weak) completion of the \, $\widetilde{}$-scheme. For a section $s$ of a coherent ${\cal O}$-module on a scheme, ${\bf V}(s)$ denotes the zero set of $s$.

\section{Log Rigid Cohomology}
\label{secrico}

\addtocounter{satz}{1}{\bf \arabic{section}.\arabic{satz}} For a weak formal $A_0$-scheme (see Meredith \cite{mer}) we require that it is locally of the form $\spwf(B)$ where $B$ is a quotient of some algebra $A_0[X_1,\ldots,X_n]^{\dagger}$, the weak completion with respect to $m_{A_0}$ of $A_0[X_1,\ldots,X_n]$. Weak completions of noetherian $A_0$-algebras or -schemes are always weak completions with respect to $m_{A_0}$. In the same way a $p$-adic formal $A_0$-scheme of topologically finite type gives rise to a $K_0$-rigid space ("generic fibre"), a weak formal $A_0$-scheme ${\mathfrak X}$ gives rise to a $K_0$-dagger space ${\mathfrak X}_{\mathbb{Q}}$ as defined in \cite{crelle}. In particular, for a weakly complete algebra $B$ as above, the algebra $B\otimes{\mathbb{Q}}$ is a $K_0$-dagger algebra and gives rise to the affinoid $K_0$-dagger space $\spm(B\otimes{\mathbb{Q}})$ (notation and terminology from \cite{crelle}). As in the formal/rigid context, there is a specialization map $sp:{\mathfrak X}_{\mathbb{Q}}\to{\mathfrak X}$. For a $k$-subscheme $Z$ of ${\mathfrak X}$ we write $]Z[_{\mathfrak X}=sp^{-1}(Z)$, the tube of $Z$ in ${\mathfrak X}$, an admissible open subspace of ${\mathfrak X}_{\mathbb{Q}}$. All this analogously for $A$ and $K$ instead of $A_0$ and $K_0$.

\begin{lem}\label{tubwohl} Let $a\in\{0,\pi\}$ and ${\mathfrak X}$ a fine weak formal ${{\mathfrak S}^{a} }$-log scheme and let $f:Z\to{\mathfrak X}$ be a closed immersion of a ${S}^{0} $-log scheme $Z$. Let $Z\to{\mathfrak X}'\to{\mathfrak X}$ and $Z\to{\mathfrak X}''\to{\mathfrak X}$ be two exactifications of $f$. Then there is a canonical isomorphism $]Z[_{{\mathfrak X}'}\cong]Z[_{{\mathfrak X}''}$.
\end{lem}

{\sc Proof:} Comparing with the tube of $Z$ in an exactification of the diagonal embedding $Z\to{\mathfrak X'}\times_{{\mathfrak X}}{\mathfrak X}''$, we see that we may assume that there is a log \'{e}tale map $q:{\mathfrak X}'\to{\mathfrak X''}$ compatible with the $Z$-embeddings and with charts for the Zariski topology. It follows from \cite{kalo} 3.8 that in a neighbourhood of $Z$, $q$ is \'{e}tale in the classical sense. Let $q_Z:]Z[_{{\mathfrak X}'}\to]Z[_{{\mathfrak X}''}$ be the morphism of dagger spaces induced by $q$. From \cite{berco} we get that the morphism of rigid spaces associated with $q_Z$ --- this is the morphism between the tubes of $Z$ in the respective $p$-adic completions ---  is an isomorphism. Therefore $q_Z$ is an isomorphism by \cite{crelle}. 
 
\addtocounter{satz}{1}{\bf \arabic{section}.\arabic{satz}} Denote by $\Delta$ the category whose objects are the ordered sets $[m]=\{0<\ldots<m\}$ for $m\in\mathbb{Z}_{\ge0}$ and whose morphisms are the {\it injective} order preserving maps of sets. For our purposes, a simplicial scheme $X_{\bullet}$ is a functor from $\Delta^{op}$ to the category of schemes. On objects we write it as $[m]\mapsto X_m$. Similarly we define simplicial (log)-schemes over a fixed base, simplicial dagger spaces and so on. 

Let $a\in\{0,\pi\}$. Let $X$ be a fine ${S}^{0} $-log scheme. Suppose we have an open covering $X=\cup_{i\in I}V_{i}$ and for every $i\in I$ an exact closed immersion $V_{i}\to{\mathfrak V}_{i}$ into a log smooth weak formal ${{\mathfrak S}^{a} }$-log scheme ${\mathfrak V}_{i}$. Choose (perhaps after refining the covering) for each non empty finite subset $H\subset {I}$ an exactification $$V_H=\cap_{i\in H}V_{i}\stackrel{\iota}{\longrightarrow}{\mathfrak V}_H\stackrel{f}{\longrightarrow}\times_{{{\mathfrak S}^{a} }}({\mathfrak V}_i)_{i\in H}$$of the diagonal embedding $V_H{\to}\times_{{{\mathfrak S}^{a} }}({\mathfrak V}_{i})_{i\in H}$. Let $\Omega^{\bullet}_{{\mathfrak V}_H/{\mathfrak S}^{a} }$ be the relative logarithmic de Rham complex of the morphism of weak formal log schemes ${\mathfrak V}_H\to{{\mathfrak S}^{a} }$. This is a sheaf complex on ${\mathfrak V}_H$, and tensoring with $\mathbb{Q}$ induces a sheaf complex $\Omega^{\bullet}_{{\mathfrak V}^{}_{H,\mathbb{Q}}}$ on the $K_0(a)$-dagger space ${\mathfrak V}^{}_{H,\mathbb{Q}}$, the generic fibre of ${\mathfrak V}_H$, as follows: If $\spwf(B)\subset{\mathfrak V}_H$ is open, then the affinoid dagger space $U=\spm(B\otimes\mathbb{Q})$ is admissible open in ${\mathfrak V}^{}_{H,\mathbb{Q}}$, and $\Omega^{q}_{{\mathfrak V}^{}_{H,\mathbb{Q}} }|_{U}$ is associated with the coherent $B\otimes\mathbb{Q}$-module $\Omega^{q}_{{\mathfrak V}_H/{\mathfrak S}^{a} }(\spwf(B))\otimes\mathbb{Q}$. For varying open affines $\spwf(B)$ this construction glues. By \ref{tubwohl} the tube $]V_H[_{{\mathfrak V}_H}$ and the restriction $\Omega^{\bullet}_{{\mathfrak V}^{}_{H,\mathbb{Q}}}|_{]V_H[_{{\mathfrak V}_H}}=\Omega^{\bullet}_{]V_H[_{{\mathfrak V}_H}}$ of $\Omega^{\bullet}_{{\mathfrak V}^{}_{H,\mathbb{Q}}}$ to $]V_H[_{{\mathfrak V}_H}$ depend only on the system $\{V_i\to{\mathfrak V}_i\}_i$, not on the chosen exactification $(\iota,f)$. For $H_1\subset H_2$ one has a canonical projection map $\delta_{H_1H_2}:]V_{H_2}[_{{\mathfrak V}_{H_2}}\to ]V_{H_1}[_{{\mathfrak V}_{H_1}}$ and a natural map $\delta_{H_1H_2}^{-1}\Omega^{\bullet}_{]V_{H_1}[_{{\mathfrak V}_{H_1}}}\to \Omega^{\bullet}_{]V_{H_2}[_{{\mathfrak V}_{H_2}}}$. Choosing a well ordering of $I$ we get as usual a simplicial dagger space $(]V_H[_{{\mathfrak V}_H})_{H\subset {J}}=]V_{\bullet}[_{{\mathfrak V}_{\bullet}}$; furthermore, multiplying the transition maps with alternating signs as usual, we get a sheaf complex $\Omega^{\bullet}_{]V_{\bullet}[_{{\mathfrak V}_{\bullet}}}$ on $]V_{\bullet}[_{{\mathfrak V}_{\bullet}}$. Given a locally constant sheaf $E$ of $K_0(a)$-vector spaces on $X_{Zar}$ we form the complex $$E\otimes_{K_0(a)}\Omega^{\bullet}_{]V_{\bullet}[_{{\mathfrak V}_{\bullet}}}=(sp^{-1}(E|_{V_{\bullet}})\otimes_{K_0(a)}\Omega^{\bullet}_{]V_{\bullet}[_{{\mathfrak V}_{\bullet}}},1\otimes d)$$ on $]V_{\bullet}[_{{\mathfrak V}_{\bullet}}$ (where $sp:]V_{{\bullet}}[_{{\mathfrak V}_{\bullet}}\to V_{\bullet}$ is the specialization map). Now suppose we are given another set of data $X'=\cup_{i\in{I}'}{V}'_i$ with embeddings $\{{V}'_i\to{\mathfrak V}'_i\}$ as above and a ${S}^{0}$-morphism $f:X'\to X$. Let $E'=f^{-1}E$. For $(i,i')\in I\times{I}'$ let $V^{\times}_{(i,i')}={V}'_{i'}\cap f^{-1}{V}_i$ and choose an exactification $V^{\times}_{(i,i')}\to{\mathfrak V}^{\times}_{(i,i')}$ of the embedding $(\id\times f):V^{\times}_{(i,i')}\to{\mathfrak V}'_{i'}\times_{{{\mathfrak S}^{a} }}{\mathfrak V}_i$. Again we get a simplicial dagger space $]V_{\bullet}^{\times}[_{{\mathfrak V}^{\times}_{\bullet}}$. There are projection maps ${\rm pr}_1:]V_{\bullet}^{\times}[_{{\mathfrak V}^{\times}_{\bullet}}\to]V_{\bullet}[_{{\mathfrak V}_{\bullet}}$ and ${\rm pr}_2:]V_{\bullet}^{\times}[_{{\mathfrak V}^{\times}_{\bullet}}\to]V'_{\bullet}[_{{\mathfrak V}'_{\bullet}}$. These give rise to \begin{gather}{\rm pr}_1^{-1}(E\otimes_{K_0(a)}\Omega^{\bullet}_{]V_{\bullet}[_{{\mathfrak V}_{\bullet}}})\longrightarrow E'\otimes_{K_0(a)}\Omega^{\bullet}_{]V^{\times}_{\bullet}[_{{\mathfrak V}^{\times}_{\bullet}}}\tag{$i$}\\{\rm pr}_2^{-1}(E'\otimes_{K_0(a)}\Omega^{\bullet}_{]V'_{\bullet}[_{{\mathfrak V}'_{\bullet}}})\longrightarrow E'\otimes_{K_0(a)}\Omega^{\bullet}_{]V^{\times}_{\bullet}[_{{\mathfrak V}^{\times}_{\bullet}}}\tag{$ii$}.\end{gather}

\begin{lem} ${\mathbb{R}}\Gamma(]V^{\times}_{\bullet}[_{{\mathfrak V}^{\times}_{\bullet}},(ii))$ is an isomorphism. In particular the definitions$${\mathbb{R}}\Gamma_{rig}(X/{\mathfrak S}^{a},E)={\mathbb{R}}\Gamma(]V_{\bullet}[_{{\mathfrak V}_{\bullet}},E\otimes_{K_0(a)}\Omega^{\bullet}_{]V_{\bullet}[})$$$${\mathbb{R}}\Gamma_{rig}(X/{\mathfrak S}^{a} )={\mathbb{R}}\Gamma_{rig}(X/{\mathfrak S}^{a},K_0(a))$$are independent on the covering $X=\cup_{i\in I}V_{i}$ and the embeddings $V_{i}\to{\mathfrak V}_{i}$. There is a natural map ${\mathbb{R}}\Gamma_{rig}(X/{\mathfrak S}^{a},E)\to {\mathbb{R}}\Gamma_{rig}(X'/{\mathfrak S}^{a},E')$.
\end{lem} 

{\sc Proof:} That ${\mathbb{R}}\Gamma(]V^{\times}_{\bullet}[_{{\mathfrak V}^{\times}_{\bullet}},(ii))$ is an isomorphism is a local claim, hence one may assume $E=K_0(a)$. But then the reasoning of \cite{berco} carries over: the key is that each projection ${\mathfrak V}^{\times}_{(i,i')}\to{\mathfrak V}_{i'}$ is strict and {\it classically} smooth near $V^{\times}_{(i,i')}$, hence $]V^{\times}_{(i,i')}[_{{\mathfrak V}^{\times}_{(i,i')}}\to]V^{\times}_{(i,i')}[_{{\mathfrak V}_{i'}}$ is a relative open polydisk so that the Poincar\'{e} lemma applies. Applied to the identity $X'=X\to X$ we get the well definedness of ${\mathbb{R}}\Gamma_{rig}$. The map ${\mathbb{R}}\Gamma_{rig}(X/{\mathfrak S}^{a},E)\to {\mathbb{R}}\Gamma_{rig}(X'/{\mathfrak S}^{a},E')$ is the one induced by $(i)$. 

\addtocounter{satz}{1}{\bf \arabic{section}.\arabic{satz}}\newcounter{rigvebe1}\newcounter{rigvebe2}\setcounter{rigvebe1}{\value{section}}\setcounter{rigvebe2}{\value{satz}} For a simplicial scheme $X_{\bullet}$ and a ring $R$ denote by ${\cal LS}(X_{\bullet},R)$ the category of sheaves on $X_{\bullet}$ with values in the category of finitely generated $R$-modules, locally constant on $(X_m)_{Zar}$ for each $m\ge0$. Let $X_{\bullet}$ be a simplicial fine ${S}^{0} $-log scheme and $E\in{\cal LS}(X_{\bullet},K_0(a))$. Suppose that there exists an open covering $X_0=\cup_{i\in I}V_i$ satisfying the following two conditions:\\(1) for any $m\ge0$, any $i\in I$ and for any $\lambda:[0]\to[m]$, if $\sigma(\lambda):X_m\to X_0$ denotes the corresponding structure morphism, then $\sigma(\lambda)^{-1}V_i=V_{m,i}$ is independent of the choice of $\lambda$.\\(2) For any $m$ and $i$ there exists an exact closed embedding $\iota_{m,i}:V_{m,i}\to{\mathfrak V}_{m,i}$ into a log smooth weak formal ${{\mathfrak S}^{a} }$-log scheme ${\mathfrak V}_{m,i}$.\\We wish to define the rigid cohomology ${\mathbb{R}}\Gamma_{rig}(X_{\bullet}/{\mathfrak S}^{a},E)$ of $X_{\bullet}$ with values in $E$.\\For the moment let us fix a non empty finite subset $H\subset {I}$. Define inductively$$\widetilde{\mathfrak Z}_{0,H}=\times_{i\in H}{\mathfrak V}_{0,i},\quad\quad\widetilde{\mathfrak Z}_{m+1,H}=(\times_{i\in H}{\mathfrak V}_{m+1,i})\times(\widetilde{\mathfrak Z}_{m,H})^{m+2}.$$With the $m+2$ natural projections $\widetilde{\mathfrak Z}_{m+1,H}\to\widetilde{\mathfrak Z}_{m,H}$ we get a simplicial weak formal ${{\mathfrak S}^{a} }$-log scheme $\widetilde{\mathfrak Z}_{\bullet,H}$. Let $V_{m,H}=\cap_{i\in H}V_{m,i}$ and let $t_{m,H}:V_{m,H}\to \times_{i\in H}{\mathfrak V}_{m,i}$ be the diagonal embedding. Define inductively embeddings $j_{m,H}:V_{m,H}\to\widetilde{\mathfrak Z}_{m,H}$ as follows: $$j_{0,H}=t_{m,H},\quad\quad j_{m+1,H}=(t_{m+1,H}\times(t_m\circ\sigma_0)\times\ldots\times(t_m\circ\sigma_{m+1}))$$where $\sigma_s:V_{m+1,H}\to V_{m,H}$ for $s=0,\ldots, m+1$ are the structure projections (obtained by restricting those of $X_{\bullet}$). We have defined a closed embedding of simplicial objects $j_{\bullet,H}:V_{\bullet,H}\to\widetilde{\mathfrak Z}_{\bullet,H}$. For any $m$ choose an exactification $$V_{m,H}\longrightarrow{\mathfrak Z}_{m,H}\longrightarrow\widetilde{\mathfrak Z}_{m,H}.$$The tubes $]V_{m,H}[_{{\mathfrak Z}_{m,H}}$ form a simplicial $K_0(a)$-dagger space $]V_{\bullet,H}[_{{\mathfrak Z}_{\bullet,H}}$. Now we vary $H$: for $H_1\subset H_2$ one has a canonical projection map of simplicial $K_0(a)$-dagger spaces $\delta_{H_1H_2}:]V_{\bullet,H_2}[_{{\mathfrak V}_{\bullet,H_2}}\to ]V_{\bullet,H_1}[_{{\mathfrak V}_{\bullet,H_1}}$ and if we fix a well ordering of $I$ we get as usual a simplicial simplicial $K_0(a)$-dagger space $]V_{\bullet,\bullet}[_{{\mathfrak Z}_{\bullet,\bullet}}=(]V_{m,H}[_{{\mathfrak Z}_{m,H}})_{m,H}$ (a functor from $\Delta^{op}$ to the category of simplicial $K_0(a)$-dagger spaces). As before there is natural logarithmic de Rham complex $E\otimes_{K_0(a)}\Omega^{\bullet}_{]V_{\bullet,\bullet}[_{{\mathfrak V}_{\bullet,\bullet}}}$ on $]V_{\bullet,\bullet}[_{{\mathfrak V}_{\bullet,\bullet}}$ and we set$${\mathbb{R}}\Gamma_{rig}(X_{\bullet}/{\mathfrak S}^{a},E)={\mathbb{R}}\Gamma(]V_{\bullet,\bullet}[_{{\mathfrak V}_{\bullet,\bullet}},E\otimes_{K_0(a)}\Omega^{\bullet}_{]V_{\bullet,\bullet}[_{{\mathfrak V}_{\bullet,\bullet}}}).$$That this is well defined is shown as before.\\
We define ${\mathbb{R}}\Gamma_{conv}(X_{\bullet}/{\mathfrak S}^{a},E)$ for $E\in {\cal LS}(X_{\bullet},K_0(a))$ by the same procedure, using $p$-adic formal schemes and rigid spaces instead of weak formal schemes and dagger spaces. Thus, for a fine ${S}^{0} $-log scheme $X$ our ${\mathbb{R}}\Gamma_{conv}(X/{\mathfrak S}^{a},K)$ for the constant coefficient $E=K$ is what Shiho calls the analytic cohomology of $X/{\mathfrak S}^{a}$ in \cite{shiho}. In particular, by \cite{shiho} it is isomorphic with the log convergent cohomology of $X/{\mathfrak S}^{a}$ in the sense of Ogus \cite{oglog}, and we have comparison isomorphisms with log crystalline cohomology if $a=0$. For a $k$-scheme $X$ with trivial log structure, if in all our constructions we replace the base ${\mathfrak S}^{a}$ by the base $\widetilde{W}$ and work with the constant sheaf $E=K_0$ on $X$, then we obtain cohomology objects which we denote by ${\mathbb{R}}\Gamma_{rig}(X/K_0)$ resp. by ${\mathbb{R}}\Gamma_{conv}(X/K_0)$. There is a canonical map from our ${\mathbb{R}}\Gamma_{rig}(X/K_0)$ to the object ${\mathbb{R}}\Gamma_{rig}(X/K_0)$ defined by Berthelot \cite{berco}, and it follows from \cite{crelle} 5.1 that this is an isomorphism. 

\addtocounter{satz}{1}{\bf \arabic{section}.\arabic{satz}}\newcounter{defbs1}\newcounter{defbs2}\setcounter{defbs1}{\value{section}}\setcounter{defbs2}{\value{satz}} See \cite{colo} for details on the following concept. A {\it log scheme with boundary} is a morphism of quasi-compact log schemes $i:(X,{\cal N}_X)\to(\overline{X},{\cal N}_{\overline{X}})$ such that\\(a) the underlying morphism of schemes is a schematically dense open immersion,\\(b) ${\cal N}_{\overline{X}}\to i_*{\cal N}_X$ is injective, $i^*{\cal N}_{\overline{X}}={\cal N}_X$ and $(i_*{\cal N}_X)^{gp}={\cal N}_{\overline{X}}^{gp}$.\\Let $T=(T,{\cal N}_T)=\widetilde{S}$ or $T=(T,{\cal N}_T)={S}$. A $T${\it -log scheme with boundary} is a log scheme with boundary together with a morphism of log schemes $X\to T$. For short we denote it by $(\overline{X},X)$. It is called fine if $(\overline{X},{\cal N}_{\overline{X}})$ is fine. Let $\Omega^{\bullet}_{\overline{X}/\widetilde{W}}$ be the logarithmic de Rham complex of the morphism of log schemes $\overline{X}\to \widetilde{W}$. The section $\dlog(t)\in\Omega^{1}_{\overline{X}/\widetilde{W}}(X)$ uniquely extends to a section $\dlog(t)\in\Omega^{1}_{\overline{X}/\widetilde{W}}(\overline{X})$. We define $\Omega^{\bullet}_{(\overline{X},X)/T}$ as the quotient of $\Omega^{\bullet}_{\overline{X}/\widetilde{W}}$ by its  ${\cal O}_{\overline{X}}$-subalgebra generated by this $\dlog(t)$.\\
A morphism of $T$-log schemes with boundary $f:(\overline{X},X)\to(\overline{X}',X')$ is a morphism of log schemes $f:\overline{X}\to \overline{X}'$ restricting to a morphism of $T$-log schemes $X\to {X}'$ (in particular $X\subset f^{-1}(X')$). It is called an exact closed immersion if $f$ is one. Taking $X$ to $(X,X)$ is a fully faithful functor from the category of $T$-log schemes to the category of $T$-log schemes with boundary. For fine $T$-log schemes with boundary $(\overline{X}_1,X_1)$, $(\overline{X}_2,X_2)$ there exists a product $(\overline{X}_1\overline{\times}_T\overline{X}_2,X_1\times_TX_2)$ in the category of fine $T$-log schemes with boundary: $\overline{X}_1\overline{\times}_T\overline{X}_2$ is the schematic closure of $X_1\times_TX_2$ in $\overline{X}_1{\times}_{\widetilde{W}}\overline{X}_2$, endowed with the image of the log structure of $\overline{X}_1{\times}_{\widetilde{W}}\overline{X}_2$ in the push forward log structure of $X_1\times_TX_2$.

\addtocounter{satz}{1}{\bf \arabic{section}.\arabic{satz}}\newcounter{defsmo1}\newcounter{defsmo2}\setcounter{defsmo1}{\value{section}}\setcounter{defsmo2}{\value{satz}} A morphism of $T$-log schemes with boundary $(\overline{Y},Y)\to(\overline{X},X)$ is called a  boundary exact closed immersion if $\overline{Y}\to\overline{X}$ is an exact closed immersion and if for every open neighbourhood $U$ of $Y$ in $X$, there exists an open neighbourhood $\overline{U}$ of $\overline{Y}$ in $\overline{X}$ with $U$ schematically dense in $\overline{U}$. A fine $T$-log scheme with boundary $(\overline{X},X)$ is called {\it smooth} if it satisfies the following conditions (1)-(3): (1) $\overline{X}$ is locally of finite presentation over $\widetilde{W}$. (2) For every morphism $\eta:(\overline{L}',L')\to(\overline{L},L)$ such that $\overline{L'}\to\overline{L}$ is an exact closed immersion defined by a square zero ideal in ${\cal O}_{\overline{L}}$ and for every morphism $\mu:(\overline{L}',L')\to(\overline{X},X)$ there is \'{e}tale locally on $\overline{L}$ a morphism $\epsilon:(\overline{L},L)\to(\overline{X},X)$ such that $\mu=\epsilon\circ\eta$. (3) For all morphisms $(\overline{Y},Y)\to(\overline{X},X)$ and all boundary exact closed immersions $(\overline{Y},Y)\to(\overline{V},V)$ of fine $T$-log schemes with boundary, there exists \'{e}tale locally on $(\overline{X}\overline{\times}_T\overline{V})$ an exactification $$\overline{Y}\to Z\to(\overline{X}\overline{\times}_T\overline{V})$$of the diagonal embedding $\overline{Y}\to(\overline{X}\overline{\times}_T\overline{V})$ (a morphism of log schemes in the usual sense) such that the projection $Z\to(\overline{X}\overline{\times}_T\overline{V})\to \overline{V}$ is strict and log smooth, or equivalently: strict and smooth on underlying schemes.

\begin{pro}\label{keykoh}(\cite{colo} Proposition 2.6) Let $(\overline{Y},Y)\to(\overline{X}_i,X_i)$ be boundary exact closed immersions into smooth $T$-log schemes with boundary ($i$ in a finite index set $I$). There exist locally on $\overline{\times}_T(\overline{X}_i)_{i\in I}$ factorizations$$(\overline{Y},Y)\stackrel{\iota}{\longrightarrow}(\overline{Z},Z)\stackrel{}{\longrightarrow}(\overline{\times}_T(\overline{X}_i)_{i\in I},{\times}_T({X_i})_{i\in I})$$of the diagonal embedding such that $\iota$ is a boundary exact closed immersion, the map $\overline{Z}\to\overline{\times}_T(\overline{X}_i)_{i\in I}$ is log \'{e}tale, and the projections $p_i:\overline{Z}\to\overline{X}_i$ are strict and log smooth, hence smooth on underlying schemes.
\end{pro}

\addtocounter{satz}{1}{\bf \arabic{section}.\arabic{satz}}\newcounter{smocrit1}\newcounter{smocrit2}\setcounter{smocrit1}{\value{section}}\setcounter{smocrit2}{\value{satz}} A chart $(t\in P^{gp}\supset P)$ for $(\overline{X},X)$ over $T$ is a chart $\lambda:P\to\Gamma(\overline{X},{\cal N}_{\overline{X}})$ for $(\overline{X},{\cal N}_{\overline{X}})$ together with an element $t_{P}\in P^{gp}$ such that $\lambda^{gp}(t_P)=t$ in $\Gamma(\overline{X},{\cal N}_{\overline{X}})^{gp}$ (with $t$ on the right hand side the image of $t\in\Gamma(T,{\mathcal N}_T)$). In \cite{colo} Theorem 2.5 it is shown that $(\overline{X},X)/T$ is smooth if locally on $\overline{X}$ there are charts $(t_P\in P^{gp}\supset P)$ satisfying the following two conditions: (i) The subgroup $(t_{P})\subset P^{gp}$ generated by $t_{P}$, as well as the torsion part of $P^{gp}/(t_{P})$ are finite groups of orders invertible on $\widetilde{W}$. (ii) The morphism of schemes $\overline{X}\to\spec(A_0[P])$ if $T=\widetilde{S}$, resp. $\overline{X}\to\spec(k[P])$ if $T={S}$, is smooth. 

\addtocounter{satz}{1}{\bf \arabic{section}.\arabic{satz}}\newcounter{dflorico1}\newcounter{dflorico2}\setcounter{dflorico1}{\value{section}}\setcounter{dflorico2}{\value{satz}}\label{ricode} Let $(\overline{X},X)$ be a ${S}$-log scheme with boundary. Suppose that there exists an open covering $\overline{X}=\cup_{i\in I}\overline{V}_{i}$ and for every $i\in I$, if we set ${V}_{i}=X\cap\overline{V}_{i}$, a boundary exact closed immersion $(\overline{V}_{i},{V}_{i})\to(\overline{\widetilde{V}}_{i},\widetilde{V}_{i})$ into a smooth ${\widetilde{S}}$-log scheme with boundary $(\overline{\widetilde{V}}_{i},\widetilde{V}_{i})$. (We do not know if such local embeddings exist for any ${S}$-log scheme with boundary; in our applications we will have explicit such local embeddings at hand.) Choose (perhaps after localizing) for each finite subset $H\subset {I}$ a factorization$$(\overline{V}_H,V_H)=(\cap_{i\in H}\overline{V}_{i},\cap_{i\in H}{V}_{i})\stackrel{\iota}{\longrightarrow}(\overline{\widetilde{V}}_H,{\widetilde{V}}_H)\stackrel{f}{\longrightarrow}(\overline\times_{\widetilde{S}}(\overline{\widetilde{V}}_i)_{i\in H},\times_{\widetilde{S}},{\widetilde{V}}_i)_{i\in H})$$ of the diagonal embedding as in \ref{keykoh}. Let $(\overline{\mathfrak{V}}_H,{\mathfrak{V}}_H)$ be the weak completion of $(\overline{\widetilde{V}}_H,{\widetilde{V}}_H)$: this might be called a weak formal $\mathfrak{S}$-log scheme with boundary. Weak completion of the de Rham complex $\Omega^{\bullet}_{(\overline{\widetilde{V}}_H,{\widetilde{V}}_H)/{\widetilde{S}}}$ gives a de Rham complex $\Omega^{\bullet}_{(\overline{\mathfrak{V}}_H,{\mathfrak{V}}_H)/{\mathfrak{S}}}$ on $\overline{\mathfrak{V}}_H$. As before, tensoring with $\mathbb{Q}$ induces a sheaf complex $\Omega^{\bullet}_{(\overline{\mathfrak{V}}_H,{\mathfrak{V}}_H)/{\mathfrak{S}}}\otimes\mathbb{Q}$ on the generic fibre (as $K_0$-dagger space) of $\overline{\mathfrak{V}}_H$. The tube $]\overline{V}_H[_{\overline{\mathfrak{V}}_H}$ with sheaf complex $\Omega^{\bullet}_{(\overline{\mathfrak{V}}_H,{\mathfrak{V}}_H)/{\mathfrak{S}}}\otimes\mathbb{Q}|_{]\overline{V}_H[_{\overline{\mathfrak{V}}_H}}$ is independent on the chosen exactification $(\iota,f)$. For varying $H$ one has natural transition maps, hence a simplicial dagger space $]\overline{V}_{\bullet}[_{\overline{\mathfrak{V}}_{\bullet}}$ and, given $E\in {\cal LS}(\overline{X},K_0)$, a sheaf complex $$E\otimes_{K_0}\Omega^{\bullet}_{(\overline{\mathfrak{V}}_{\bullet},{\mathfrak{V}}_{\bullet})/{\mathfrak{S}}}\otimes\mathbb{Q}|_{]\overline{V}_{\bullet}[_{\overline{\mathfrak{V}}_{\bullet}}}=sp^{-1}(E|_{\overline{X}_{\bullet}})\otimes_{K_0}\Omega^{\bullet}_{(\overline{\mathfrak{V}}_{\bullet},{\mathfrak{V}}_{\bullet})/{\mathfrak{S}}}\otimes\mathbb{Q}|_{]\overline{V}_{\bullet}[_{\overline{\mathfrak{V}}_{\bullet}}}$$ on $]\overline{V}_{\bullet}[_{\overline{\mathfrak{V}}_{\bullet}}$. As before one shows that$${\mathbb{R}}\Gamma_{rig}((\overline{X},X)/{\mathfrak{S}},E)={\mathbb{R}}\Gamma(]\overline{V}_{\bullet}[_{\overline{\mathfrak{V}}_{\bullet}},E\otimes_{K_0}\Omega^{\bullet}_{(\overline{\mathfrak{V}}_{\bullet},{\mathfrak{V}}_{\bullet})/{\mathfrak{S}}}\otimes\mathbb{Q}|_{]\overline{V}_{\bullet}[_{\overline{\mathfrak{V}}_{\bullet}}})$$is independent of choices. (Caution: Despite what our notation might suggest there is no ${\cal O}_{\mathfrak S}$-action on ${\mathbb{R}}\Gamma_{rig}((\overline{X},X)/{\mathfrak{S}},E)$).\\More generally, let $(\overline{X}_{\bullet},X_{\bullet})$ be a simplicial ${S}$-log scheme with boundary. Suppose that there is an open covering $\overline{X}=\cup_{i\in I}\overline{V}_i$ satisfying the following two conditions:\\(1) for any $m\ge0$, any $i\in I$ and for any $\lambda:[0]\to[m]$, if $\sigma(\lambda):\overline{X}_m\to \overline{X}_0$ denotes the corresponding structure morphism, then $\sigma(\lambda)^{-1}\overline{V}_i=\overline{V}_{m,i}$ is independent of the choice of $\lambda$.\\(2) For any $m$ and $i$, if we set ${V}_{m,i}=X_m\cap\overline{V}_{m,i}$, there exists a boundary exact closed immersion $(\overline{V}_{m,i},{V}_{m,i})\to(\overline{\widetilde{V}}_{m,i},\widetilde{V}_{m,i})$ into a smooth ${\widetilde{S}}$-log scheme with boundary.\\From \arabic{rigvebe1}.\arabic{rigvebe2} it is now clear how we define $${\mathbb{R}}\Gamma_{rig}((\overline{X}_{\bullet},X_{\bullet})/{\mathfrak{S}},E)$$ for $E\in {\cal LS}(\overline{X}_{\bullet},K_0)$.

\addtocounter{satz}{1}{\bf \arabic{section}.\arabic{satz}}\newcounter{fro1}\newcounter{fro2}\setcounter{fro1}{\value{section}}\setcounter{fro2}{\value{satz}} Fix $q\in p^{\mathbb{N}}$ and suppose that $\sigma:A_0\to A_0$ lifts the $q$-th power Frobenius endomorphism of $k$. We also denote by $\sigma$ the unique endomorphism of the formal log scheme ${\mathfrak S}^{0} $ which equals $\sigma$ on the structure sheaf and for which multiplication with $q$ on the standard chart $\mathbb{N}$ of ${\mathfrak S}^{0} $ is a chart. Let $X$ be a fine ${S}^{0} $-log scheme. Denote by $F\mbox{-}{\cal LS}(X,K_0)$ the category of pairs $(E,\phi)$, where $E\in{\cal LS}(X,K_0)$ and $\phi:E\to E$ is a $\sigma$-linear automorphism. The datum of $\phi$ is equivalent with a $K_0$-linear isomorphism $\phi:E\otimes_{K_0,\sigma}K_0\to E$. Given $(E,\phi)\in F\mbox{-}{\cal LS}(X,K_0)$ we define the $\sigma$-linear endomorphism $\phi$ on $H^*_{rig}(X/{\mathfrak S}^{0},E)$ as the composite$$H^*_{rig}(X/{\mathfrak S}^{0},E)\stackrel{H^*(\sigma)}{\longrightarrow} H^*_{rig}(X'/{\mathfrak S}^{0},E)\stackrel{H^*(\phi')}{\longrightarrow} H^*_{rig}(X/{\mathfrak S}^{0},E)$$with $X'=X\times_{{\mathfrak S}^{0},\sigma}{\mathfrak S}^{0}$, where $H^*(\sigma)$ is induced by the base change $\sigma:{\mathfrak S}^{0} \to {\mathfrak S}^{0} $, and where $H^*(\phi')$ is induced by the structure map $\phi:E\otimes_{K_0}K_0\to E$ and the first arrow in the natural factorization $X\to X'\to X$ of the $q$-th power Frobenius endomorphism.

\section{Natural projective bundles on a semistable scheme}
\label{bassec}

\addtocounter{satz}{1}{\bf \arabic{section}.\arabic{satz}}\newcounter{basob1}\newcounter{basob2}\setcounter{basob1}{\value{section}}\setcounter{basob2}{\value{satz}} Our basic object of study in this paper is $Y$: a strictly semistable log scheme over $k$. That is, $Y$ is a fine ${S}^{0} $-log scheme $(Y,{\cal N}_Y)$ which allows a {\it Zariski} open covering by open subschemes $Y'\subset Y$ with the following property: there exist integers $m\ge 1$ and charts $\mathbb{N}^m\to{\cal N}_{Y}(Y')$ for ${\cal N}_Y|_{Y'}$ such that\\(i) if on the log scheme ${S}^{0} $ we use the chart $\mathbb{N}\to k, 1\mapsto 0$, the diagonal morphism $\mathbb{N}\stackrel{\delta}{\to}\mathbb{N}^m$ is a chart for the structure morphism of log schemes $Y'\to {S}^{0} $, and\\(ii) the induced morphism of schemes $$f:Y'\longrightarrow\spec(k) \times_{\spec(k[t])}\spec(k[t_1,\ldots,t_m])$$ is smooth in the classical sense. 
If not said otherwise, {\it we endow subschemes of} $Y$ {\it with the pull back structure of} ${S}^{0}$-{\it log scheme structure induced by that of} $Y$. For simplicity we suppose that $Y$ is connected. By $\{Y_j\}_{j\in \Upsilon}$ we denote the set of irreducible components of $Y$: all of them are classically smooth over $k$; this follows from the existence of charts as above for the Zariski topology.

Conversely, suppose for the moment that we require the existence of charts as above only for the {\it \'{e}tale} topology of $Y$. We claim that then the classical smoothness of all $Y_j$ automatically implies the existence of charts as above for the Zariski topology of $Y$. Indeed, let $y\in Y$ be a closed point. Suppose it lies on the irreducible components $Y_1,\ldots, Y_r$, but not on any other irreducible component. For any $1\le j\le r$ let $t_j\in{\mathcal O}_{Y,y}$ be an element with ${\mathbb V}(t_j)=Y_j\cap\spec({\mathcal O}_{Y,y})$. The smoothness of all $Y_j$ implies that also $Z=Y_1\cap\ldots\cap Y_r$ is smooth. Therefore we find an open affine neighbourhood $\spec(B)$ of $y$ in $Z$ together with an \'{e}tale morphism of $k$-algebras $b:k[t_{r+1},\ldots,t_d]\to B$. Shrinking $\spec(B)$ if necessary we may write $\spec(B)=\spec(C)\cap Z$ for an affine open subscheme $\spec(C)$ of $Y$ with $t_1,\ldots,t_r\in C$. For any $r+1\le j\le d$ choose a preimage $t_j\in C$ of $b(t_j)$ under the surjection $C\to B$ and consider the map$$k[t_1,\ldots,t_d]/(t_1\cdot\ldots\cdot t_r)\longrightarrow C,\quad t_j\mapsto t_j.$$By construction, it induces an \'{e}tale map from an open neigbourhood of $y$ in $\spec(C)$ to $\spec(k[t_1,\ldots,t_d]/(t_1\cdot\ldots\cdot t_r))$ and our claim is proved.

For $i\ge 1$ we let $$Y^i=\coprod_{|{J}|=i}\times_Y(Y_{j})_{j\in {J}}$$ (disjoint sum) where ${J}$ runs through the set of subsets of $\Upsilon$ with precisely $i$ elements. In \cite{kalo} p.222/223 it is explained how the log structure on $Y$ gives rise, for every $j\in \Upsilon$, to an invertible ${\cal O}_Y$-module ${\cal L}_j$ and a global section $s_j$ of ${\cal L}_j$ with ${\bf V}(s_j)=Y_j$. These data form a log structure in the sense of Faltings. This assignment goes as follows: For $j\in \Upsilon$ let ${\cal N}_{Y,j}$ be the subsheaf of ${\cal N}_Y$ which is the preimage of $\ke({\cal O}_Y\to{\cal O}_{Y_j})$. This ${\cal N}_{Y,j}$ is a principal homogeneous space over ${\cal O}_Y^{\times}$, and its associated line bundle is the dual ${\cal L}_j^{-1}$ of ${\cal L}_j$. The section $s_j$ is the one corresponding to ${\cal N}_{Y,j}\to{\cal O}_Y$. There is a structure isomorphism ${\cal O}_Y\cong\otimes({\cal L}_j)_{j\in \Upsilon}$ which sends $1_{{\cal O}_Y}$ to $\underset{j\in \Upsilon}\otimes s_j$. It constitutes the ${S}^{0} $-log scheme structure of $Y$.

\addtocounter{satz}{1}{\bf \arabic{section}.\arabic{satz}}\newcounter{lam1}\newcounter{lam2}\setcounter{lam1}{\value{section}}\setcounter{lam2}{\value{satz}} Now fix a non-empty subset ${J}\subset \Upsilon$ of $\Upsilon$ with $i$ elements for some $i\ge 1$. Suppose $$M=M_{J}=\cap_{j\in {J}}Y_j=\times_Y(Y_{j})_{j\in {J}}$$ is non-empty and let $\widehat{J}=\Upsilon-{J}$. For $j\in \Upsilon$ denote the restriction of ${\cal L}_j$ to $M$ by abuse of notation again by ${\cal L}_j$ (where we should write ${\cal L}_j\otimes_{{\mathcal O}_Y}{\mathcal O}_M$), an invertible ${\cal O}_M$-module. For $j\in {J}$ we regard the affine line bundle $$V_M^j=\spec(\sym_{{\cal O}_M}({\cal L}^{-1}_j))$$ over $M$ as an open subscheme of the projective line bundle $$P_M^j=\proj(\sym_{{\cal O}_M}({\cal O}_M\bigoplus{\cal L}_j^{-1}))$$over $M$ as follows. Let $\sym_{{\cal O}_M}({\cal O}_M\oplus{\cal L}_j^{-1})[1_{{\cal O}_M}^{-1}]$ be the $\mathbb{Z}$-graded algebra obtained by inverting the degree one element $1_{{\cal O}_M}$ of $\sym_{{\cal O}_M}({\cal O}_M\oplus{\cal L}_j^{-1})$, and let $\sym_{{\cal O}_M}({\cal O}_M\oplus{\cal L}_j^{-1})[1_{{\cal O}_M}^{-1}]_0$ be its degree zero part. Then $\spec(\sym_{{\cal O}_M}({\cal O}_M\oplus{\cal L}_j^{-1})[1_{{\cal O}_M}^{-1}]_0)$ is an open subscheme of $P_M^j$ which we identify with $V_M^j$ by means of the isomorphism $$\sym_{{\cal O}_M}({\cal L}^{-1}_j)\longrightarrow \sym_{{\cal O}_M}({\cal O}_M\bigoplus{\cal L}_j^{-1})[1_{{\cal O}_M}^{-1}]_0,\quad s\mapsto 1_{{\cal O}_M}^{-1}\otimes s\quad\mbox{for } s\in {\cal L}_j^{-1}.$$For a subset ${J'}\subset {J}$ let $$\label{pmchecki}P^{{J'}}_M=\times_M(P_M^j)_{j\in {J'}}.$$With the above embeddings $V_M^j\to P_M^j$ we may regard$$V_M^{J'}=\spec(\sym_{{\cal O}_M}(\bigoplus_{j\in {J'}}{\cal L}^{-1}_j))\cong\times_M(V_M^j)_{j\in {J'}}$$ as an open subscheme of $P_M^{J'}$. (Attention: for $j\in J$ do not confuse $P^{\{j\}}$ resp. $V^{\{j\}}$ with $P^j$ resp. $V^j$). We let $P_M=P_M^J$, $V_M=V_M^J$. For $j\in {J}$ let $N_{j,\infty}$ be the divisor on $P_M$ which is the pullback to $P_M$ of the divisor $P_M^j-V_M^j$ on $P_M^j$. Let $N_{j,0}$ be the divisor on $P_M$ which is the pullback to $P_M$ of the zero section divisor $M\to V_M^j\to P_M^j$ on $P_M^j$. Set $$N_{\infty}=\cup_{j\in {J}}N_{j,\infty}\quad\quad\quad N_0=\cup_{j\in {J}}N_{j,0}.$$Viewing $P^{{J'}}_M$ as a closed subscheme of $P_M$ by identifying it with the intersection (in $P_M$) of all $N_{j,0}$ for $j\in {J}-{J'}$, we get natural closed embeddings$$M=P_M^{\emptyset}\to P^{{J_1'}}_M\to P_M^{J_2'}$$for $J'_{1}\subset J'_2\subset J$. The divisor $M\cap\cup_{j\in\widehat{J}}Y_j$ on $M$ defines a divisor $D_M$ on $P_M$ by pull back via the structure map $P_M\to M$. Endow $P_M$ with the log structure ${\cal N}_{P_M}$ defined by the normal crossings divisor $N_0\cup N_{\infty}\cup D_M$ on $P_M$.\\ 
The sections $s_j$ of ${\cal L}_j$ for $j\in \widehat{J}$ define a map $\otimes({\cal L}_j^{-1})_{j\in \widehat{J}}\stackrel{\eta}{\to}{\cal O}_M$. Hence a map $${\cal O}_M\cong(\otimes({\cal L}_j^{-1})_{j\in {J}})\otimes(\otimes({\cal L}_j^{-1})_{j\in \widehat{J}})\stackrel{1\otimes\eta}{\longrightarrow}\otimes({\cal L}_j^{-1})_{j\in {J}}\longrightarrow\sym_{{\cal O}_M}(\bigoplus_{j\in {J}}{\cal L}_j^{-1})$$of ${\cal O}_M$-modules, where the first isomorphism is obtained by dualizing the inverse of the structure isomorphism ${\cal O}_M\cong\otimes({\cal L}_j)_{j\in \Upsilon}$ (the restriction to $M$ of the global structure isomorphism). By sending $t$ to the image of $1_{{\cal O}_M(M)}$ under this map we define a morphism $\lambda:V_M\to S.$ We denote $\lambda^*(t)\in\sym_{{\cal O}_M}(\oplus_{j\in {J}}{\cal L}_j^{-1})(M)={\cal O}_{V_M}(V_M)$ again by $t$. Note $M\cap\cup_{j\in\widehat{J}}Y_j={\bf V}(\prod_{j\in \widehat{J}}s_j)$ and thus ${\bf V}(t)=(N_0\cup N_{\infty}\cup D_M)\cap V_M$, therefore $\lambda$ is a morphism of log schemes (for the log structure ${\cal N}_{P_M}|_{V_M}$ on $V_M$). The zero section $M\to V_M$ is an exact closed embedding and its composite with $\lambda$ factors as $M\stackrel{\kappa_M}{\to}{S}^{0}\to S$ where ${\kappa_M}$ is the morphism of log schemes $M\to Y\to {S}^0$. Viewing $t$ as an element of ${\cal N}_{P_M}({V_M})$ it extends uniquely to an element ${\cal N}^{gp}_{P_M}({P_M})$. In fact $(P_M,V_M)$ is an $S$-log scheme with boundary. Taking pull back log structures, each $(P_M^{J'},V_M^{J'})$ for $J'\subset J$ becomes an $S$-log scheme with boundary.

\addtocounter{satz}{1}{\bf \arabic{section}.\arabic{satz}} A local description. Let $Y'\subset Y$ be an open subscheme with a trivialization $f$ as above. The log structure on $Y'$ can be identified with the pullback, via $f$ and the canonical inclusion, of the log structure ${\cal N}_{{\bf A}_k^m}$ on the affine space ${\bf A}_k^m=\spec(k[t_1,\ldots,t_m])$ given by its divisor ${\bf V}(t_1\ldots t_m)\subset{\bf A}_k^m$. Each $t_u$ can be regarded as an element of ${\cal N}_{{\bf A}_k^m}({\bf A}_k^m)$ and we denote by $\tilde{t}_u\in{\cal N}_Y(Y')$ its image. It generates a principal homogeneous space ${\cal O}_{Y'}^{\times}.\tilde{t}_u$ over ${\cal O}_{Y'}^{\times}$, and its associated line bundle is ${\cal L}^{-1}_{r_u}|_{Y'}$ for a uniquely determined $r_u\in \Upsilon$. After a change of indexation we have an inclusion $\{1,\ldots,m\}\subset\Upsilon$ such that $r_u=u$. Now suppose $M\cap Y'\ne\emptyset$. Then $J\subset\{1,\ldots,m\}$ inside $\Upsilon$. We can view the element $\otimes_{j\in {J}}\tilde{t}_{j}\in\otimes({\cal L}^{-1}_j)_{j\in {J}}(Y')$ as a homogeneous element of degree $|J|$ in $\sym_{{\cal O}_{Y}}(\oplus_{j\in {J}}{\cal L}_j^{-1})(Y')$ and by restriction we get a homogeneous element $a_1\in\sym_{{\cal O}_M}(\oplus_{j\in {J}}{\cal L}_j^{-1})(M\cap Y')$ of degree $|{J}|$. On the other hand let $a_2'=\prod_uf^*(t_u)\in{\cal O}_Y(Y')$ where the product runs through all $1\le u\le m$ with $u\notin {J}$. Let $a_2\in {\cal O}_M(M\cap Y')$ be its restriction, which we view as an element of degree zero in $\sym_{{\cal O}_M}(\oplus_{j\in {J}}{\cal L}_j^{-1})(M\cap Y')$. Multiplying we get the element $t=a_1.a_2\in\sym_{{\cal O}_M}(\oplus_{j\in {J}}{\cal L}_j^{-1})(M\cap Y')$.

\addtocounter{satz}{1}{\bf \arabic{section}.\arabic{satz}}\newcounter{locemb1}\newcounter{locemb2}\setcounter{locemb1}{\value{section}}\setcounter{locemb2}{\value{satz}} For $J\subset \Upsilon$ and $M=M_J$, each $(P_M^{J'},V_M^{J'})$ for $J'\subset J$ admits locally a boundary exact closed immersion into a smooth $\widetilde{S}$-log scheme with boundary. Indeed, we may suppose $M=\spec(B)$ is affine and that there is for any $j\in\widehat{J}$ a $\nu_j\in B$ such that $M\cap Y_j={\bf V}({\nu}_{j})$, and for any $j\in J$ a $\mu_{j}\in {\cal L}_j^{-1}(M)$ which generates ${\cal L}_j^{-1}$ (so that this is a trivial line bundle) and such that $(\prod_{j\in \widehat{J}}\nu_{j})(\bigotimes_{j\in {J}}\mu_{j})=t$ in ${\cal O}_{V_{M}}$. Lift $B$ to a smooth affine $A_{0}$-algebra $\widetilde{B}$ (see \cite{elk}) and the elements $\nu_j\in B$ to elements $\widetilde{\nu}_j\in \widetilde{B}$ such that ${\bf V}(\prod_{j\in \widehat{J}}\widetilde{\nu}_{j})$ is a normal crossings divisor (relative $\widetilde{W}$) on $\widetilde{M}=\spec(\widetilde{B})$. Define the product of projective bundles
$$\widetilde{P}_{M}=\times_{\widetilde{M}}(\proj (\widetilde{B}[z_{j,0},z_{j,1}]))_{j\in {J}}$$over $\widetilde{M}$, a smooth $A_0$-scheme. We identify the reduction of $\widetilde{P}_{M}$ modulo $m_{A_0}$ with $P_{M}=\times_{M}(\proj(\sym_{{\mathcal O}_M}({\mathcal O}_M\oplus {\cal L}_j^{-1})))_{j\in {J}}$ by viewing $z_{j,1}$ as a lift of $\mu_{j}$ and viewing $z_{j,0}$ as a lift of the degree one element $1_{{\mathcal O}_M}$. Define the open subscheme $$\widetilde{V}_{M}=\times_{\widetilde{M}}(\spec( \widetilde{B}[x_{j}]))_{j\in {J}}=\spec( \widetilde{B}[x_{j}]_{j\in {J}})\quad \mbox{ with } \quad x_{j}=\frac{z_{j,1}}{z_{j,0}}$$of $\widetilde{P}_{M}$. Give $\widetilde{P}_{M}$ the log structure ${\cal N}_{\widetilde{P}_{M}}$ defined by the divisor $$\widetilde{D}_{M}\cup(\widetilde{P}_{M}-\widetilde{V}_{M})\cup\widetilde{N}_0$$with $\widetilde{D}_{M}={\bf V}(\prod_{j\in \widehat{J}}\widetilde{\nu}_{j})$ (zero set in $\widetilde{P}_{M}$) and $\widetilde{N}_0$ the closure of ${\bf V}(\prod_{j\in {J}}x_{j})$ in $\widetilde{P}_{M}$. As before,$$t\mapsto\prod_{j\in \widehat{J}}\widetilde{\nu}_{j}\prod_{j\in {J}}x_{j}\in\Gamma({\widetilde{V}_{M}},{\cal O}_{\widetilde{V}_{M}})$$defines a morphism of log schemes $\lambda:\widetilde{V}_{M}\to\widetilde{S}$ restricting on $V_{M}$ to $\lambda$ from \arabic{lam1}.\arabic{lam2}. For $J'\subset J$ the composite $(P_M^{J'},V_M^{J'})\to(P_M,V_M)\to (\widetilde{P}_{M},\widetilde{V}_M)$ is a boundary exact closed immersion of $\widetilde{S}$-log schemes with boundary. Moreover, applying the criterion \arabic{smocrit1}.\arabic{smocrit2} we see that $(\widetilde{P}_{M},\widetilde{V}_M)/\widetilde{S}$ is smooth. Indeed, local charts are of the form $P=\mathbb{N}^{s}$ with $t_P=\oplus_{n=1}^s\pm 1_n\in P^{gp}=\mathbb{Z}^{s}$ (some $s\ge0$).

\section{The Hyodo-Kato isomorphisms}

\begin{satz}\label{hkrigis} For $\emptyset\ne J'\subset J\subset \Upsilon$ and $a\in\{0,\pi\}$ the map$${\mathbb{R}}\Gamma_{rig}((P_{M}^{{J'}},V_{M}^{{J'}})/{\mathfrak S})\otimes_{K_0}K_0(a){\longrightarrow}{\mathbb{R}}\Gamma_{rig}(M/{\mathfrak S}^{a})$$defined by restricting to the zero section $M=M_J\to P_M^{J'}$ and specializing $t\mapsto a$ is an isomorphism. In particular, we have a natural isomorphism, independent on the choice of $J'$ but depending on the choice of $\pi$,$${\mathbb{R}}\Gamma_{rig}(M/{\mathfrak S}^{0})\otimes_{K_0}K\cong{\mathbb{R}}\Gamma_{rig}(M/{\mathfrak S}^{\pi}).$$
\end{satz}

The proof will be given in section \ref{secreshea}.

\addtocounter{satz}{1}{\bf \arabic{section}.\arabic{satz}}\newcounter{redsim1}\newcounter{redsim2}\setcounter{redsim1}{\value{section}}\setcounter{redsim2}{\value{satz}} For $m\ge0$ define $$\Lambda_m(\Upsilon)=\{\lambda=(J_0(\lambda),\ldots,J_{m}(\lambda));\quad\emptyset\ne J_0(\lambda)\subsetneq J_1(\lambda)\subsetneq\ldots\subsetneq J_{m}(\lambda)\subset\Upsilon\},$$$$(P_{\lambda},V_{\lambda})=(P_{M_{{J_m}(\lambda)}}^{{J_0}(\lambda)},V_{M_{{J_m}(\lambda)}}^{{J_0}(\lambda)})\quad\mbox{for}\quad\lambda\in\Lambda_m(\Upsilon).$$If $m>0$ "forgetting $J_s(\lambda)$" defines for $0\le s\le m$ a map $d_s:\Lambda_{m}(\Upsilon)\to\Lambda_{m-1}(\Upsilon)$. We have natural exact closed embeddings $(P_{\lambda},V_{\lambda})\to(P_{d_s(\lambda)},V_{d_s(\lambda)})$, hence a simplicial $S$-log scheme with boundary$$(P_{\bullet},V_{\bullet})=((P_m,V_m)=\coprod_{\lambda\in\Lambda_m(\Upsilon)}(P_{\lambda},V_{\lambda}))_{m\ge0}.$$The simplicial $S^0$-log scheme$$M_{\bullet}=(M_m=\coprod_{\lambda\in\Lambda_m(\Upsilon)}M_{{J_m}(\lambda)})_{m\ge0}$$comes with an augmentation $M_{\bullet}\to Y$ and a morphism (namely the exact closed immersions $M_{{J_m}(\lambda)}\to P_{M_{{J_m}(\lambda)}}^{{J_0}(\lambda)}$) of simplicial $S$-log schemes with boundary$$M_{\bullet}\to(P_{\bullet},V_{\bullet}).$$

It is clear that $M_{\bullet}$ and $(P_{\bullet},V_{\bullet})$ are functorial in open subschemes $Y'$ of $Y$: If $M'_{\bullet}$ and $(P'_{\bullet},V'_{\bullet})$ denote the corresponding simplicial $S$-log schemes (with boundary) constructed from $Y'$ instead of $Y$, then there are natural open embeddings $M'_{\bullet}\to M_{\bullet}$ and  $(P'_{\bullet},V'_{\bullet})\to(P'_{\bullet},V'_{\bullet})$, and moreover these open embeddings form cartesian diagrams with the simplicial structure maps: for each of the $m+1$ structure maps $M_{m}\to M_{m-1}$, the pullback via $M'_{m-1}\to M_{m-1}$ is the corresponding structure map $M'_{m}\to M'_{m-1}$; similarly for $(P_{\bullet},V_{\bullet})$. On the other hand, $Y$ is covered by open subschemes $Y'$ all of whose component intersections satisfy the hypotheses of \arabic{locemb1}.\arabic{locemb2}. Hence, the corresponding objects $(P'_{\bullet},V'_{\bullet})$ can be embedded as in \arabic{locemb1}.\arabic{locemb2}. 

It follows that $M_{\bullet}$ and $(P_{\bullet},V_{\bullet})$ satisfy the hypotheses of \arabic{dflorico1}.\arabic{dflorico2} so that we can define their rigid cohomology.\\

\addtocounter{satz}{1}{\bf \arabic{section}.\arabic{satz}} Let $E\in{\cal LS}(Y,K)$. For any $J\subset\Upsilon$ the pullback $E|_{M_J}$ to $M_J$ is constant; therefore, the datum $E$ is equivalent with the collection of the values $E(M_J)$ for each $J\subset\Upsilon$, together with the restriction maps $E(M_{J_1})\to E(M_{J_2})$ for $J_1\subset J_2$. In particular $E$ gives rise to objects $E\in{\cal LS}(M_{\bullet},K)$ and $E\in{\cal LS}(P_{\bullet},K)$. Restricting scalars we view $E$ also as objects in ${\cal LS}(.,K_0)$. For a weak formal scheme ${\mathfrak X}$ and a closed subscheme $Z$ of its special fibre, if $Z=\cup_{i\in I}Z_i$ is a finite covering by closed subschemes of $Z$, then the dagger space covering $]Z[_{\mathfrak X}=\cup_{i\in I}]Z_i[_{\mathfrak X}$ is admissible open. This follows from the corresponding formal/rigid fact (\cite{berco} 1.1.14.) together with the comparison principles from \cite{crelle}. Therefore the closed covering $Y=\cup_{j\in\Upsilon}Y_j$ leads to Cech spectral sequences\begin{gather}E_{1}^{rs}=E(Y^{r+1})\otimes_{K_0}H^s_{rig}(Y^{r+1}/{\mathfrak S}^{0} )\Longrightarrow H^{r+s}_{rig}(Y/{\mathfrak S}^{0},E)\tag{$1$}\\E_{1}^{rs}=E(Y^{r+1})\otimes_KH^s_{rig}(Y^{r+1}/{\mathfrak S}^{\pi} )\Longrightarrow H^{r+s}_{rig}(Y/{\mathfrak S}^{\pi},E).\tag{$2$}\end{gather}

\begin{satz}\label{ssiso} For $a\in\{0,\pi\}$ the natural maps$${\mathbb{R}}\Gamma_{rig}(Y/{\mathfrak S}^{a},E)\stackrel{\alpha^a}{\longrightarrow} {\mathbb{R}}\Gamma_{rig}(M_{\bullet}/{\mathfrak S}^{a},E)\stackrel{\beta^a}{\longleftarrow}{\mathbb{R}}\Gamma_{rig}((P_{\bullet},V_{\bullet})/{\mathfrak S},E)$$are isomorphisms. In particular, we have a natural isomorphism (depending on $\pi$)$${\mathbb{R}}\Gamma_{rig}(Y/{\mathfrak S}^{0},E)\cong{\mathbb{R}}\Gamma_{rig}(Y/{\mathfrak S}^{\pi},E).$$It comes along with a canonical isomorphism of spectral sequences $(1)\cong(2)$.
\end{satz}
   
{\sc Proof:} That $\beta^a$ is an isomorphism follows from \ref{hkrigis}. To see that $\alpha^a$ is an isomorphism we may argue locally, hence assume that there is an exact closed embedding of $Y$ into a smooth weak formal ${\mathfrak S}^a$-log scheme ${\mathfrak Y}$. Then ${\mathbb{R}}\Gamma_{rig}(M_{\bullet}/{\mathfrak S}^{a},E)$ can be computed through sheaf complexes on the tubes $]M_J[_{{\mathfrak Y}}$ (for the various $J\subset\Upsilon$) and the claim follows from \ref{cechel}. In particular we get isomorphisms between the abutment terms of $(1)$ and $(2)$. For the $E_1$-terms we get isomorphisms from \ref{hkrigis}, compatible with the differentials (i.e. restriction maps).

\begin{lem}\label{cechel} Let $X$ be a (Grothendiek) topological space, $X=\cup_{j\in \Upsilon}U_j$ an (admissible) open covering, ${\mathcal K}^{\bullet}$ an abelian sheaf complex on $X$. For $J\subset{\Upsilon}$ set $U_J=\cup_{j\in J}U_j$. For the simplicial space $U_{\bullet}=(U_m=\coprod_{\lambda\in\Lambda_m(\Upsilon)}U_{J_m(\lambda}))_{m\ge0}$ we have an isomorphism$${\mathbb{R}}\Gamma(X,{\mathcal K}^{\bullet})\cong{\mathbb{R}}\Gamma(U_{\bullet},{\mathcal K}^{\bullet}|_{U_{\bullet}}).$$
\end{lem}

{\sc Proof:} This is a variant on usual Cech cohomology. For example one can compare with the cohomology of the full Cech complex associated with the (highly redundant!) covering $X=\cup_{J\subset\Upsilon}U_J$. 

\addtocounter{satz}{1}{\bf \arabic{section}.\arabic{satz}}\newcounter{lifsit1}\newcounter{lifsit2}\setcounter{lifsit1}{\value{section}}\setcounter{lifsit2}{\value{satz}} Let ${\mathfrak X}$ be a strictly semistable weak formal $A$-scheme: locally for the Zariski topology it admits \'{e}tale maps to $\spwf(A[X_1,\ldots,X_n]^{\dagger}/(X_1\ldots X_r-\pi))$ with $1\le r\le n$. Its special fibre $Y$ defines a log structure on ${\mathfrak X}$ rendering ${\mathfrak X}$ log smooth over ${\mathfrak S}^{\pi}$; pulling back to $Y$ we get a strictly semistable log scheme over $k$ as in \arabic{basob1}.\arabic{basob2} for which we keep our previous notations. Denote by ${\mathfrak X}_{\mathbb{Q}}$ the generic fibre of ${\mathfrak X}$, a $K$-dagger space. Denote by ${\cal LS}({\mathfrak X}_{\mathbb{Q}},K)$ the category of sheaves $F$ on ${\mathfrak X}_{K}$ with values in the category of finite dimensional $K$-vector spaces, with the property that for each irreducible component $Y_j$ of $Y$, the restriction of $F$ to the tube $]Y_j[_{\mathfrak X}$ is constant. Since for each $]Y_j[_{\mathfrak X}$ the associated Berkovich analytic space is contractible (see \cite{berko}), ${\cal LS}({\mathfrak X}_{\mathbb{Q}},K)$ is equivalent with the category of locally constant sheaves of finite dimensional $K$-vector spaces on the Berkovich analytic space  associated with ${\mathfrak X}_{\mathbb{Q}}$. Moreover, the functors $sp_*$ and $sp^{-1}$ induce an equivalence of categories ${\cal LS}({\mathfrak X}_{\mathbb{Q}},K)\cong {\cal LS}(Y,K)$. For $F\in{\cal LS}({\mathfrak X}_{\mathbb{Q}},K)$ we let$${\mathbb{R}}\Gamma_{dR}({\mathfrak X}_{\mathbb{Q}},F)={\mathbb{R}}\Gamma({\mathfrak X}_{\mathbb{Q}},(F\otimes_K\Omega^{\bullet}_{{\mathfrak X}_{\mathbb{Q}}},1\otimes d)).$$For $s\ge1 $ write $]Y^{s}[_{\mathfrak X}=\coprod_u]Y^{s}_u[_{\mathfrak X}$ where the coproduct runs through the set $\{Y^{s}_u\}_u$ of connected components of $Y^{s}$; these components are subschemes of $Y$. Setting ${E}=sp_*F$ the spectral sequence $(2)$ can be identified with\begin{gather}E_{1}^{rs}=F(]Y^{r+1}[_{\mathfrak X})\otimes_KH^s_{dR}(]Y^{r+1}[_{\mathfrak X})\Longrightarrow H^{r+s}_{dR}({\mathfrak X}_{\mathbb{Q}},F).\tag{$2_{\mathfrak X}$}\end{gather}

\begin{kor}\label{rigs} There are canonical isomorphisms $(1)\cong (2_{\mathfrak X})$ and$${\mathbb{R}}\Gamma_{dR}({\mathfrak X}_{\mathbb{Q}},F)\cong {\mathbb{R}}\Gamma_{rig}(Y/{\mathfrak S}^{0},E).$$
\end{kor}

\addtocounter{satz}{1}{\bf \arabic{section}.\arabic{satz}}\newcounter{coueco1}\newcounter{coueco2}\setcounter{coueco1}{\value{section}}\setcounter{coueco2}{\value{satz}} Let ${\mathfrak X}'$ be the completion of ${\mathfrak X}$ along its special fibre and let ${\mathfrak X}'_{\mathbb{Q}}$ be its generic fibre, a $K$-rigid space. We have spectral sequences\begin{gather}E_{1}^{rs}={E}(Y^{r+1})\otimes_{K_0}H^s_{conv}(Y^{r+1}/{\mathfrak S}^{0} )\Longrightarrow H^{r+s}_{conv}(Y/{\mathfrak S}^{0},E)\tag{$1'$}\\E_{1}^{rs}=F(]Y^{r+1}[_{\mathfrak X})\otimes_KH^s_{dR}(]Y^{r+1}[_{{\mathfrak X}'})\Longrightarrow H^{r+s}_{dR}({\mathfrak X}'_{\mathbb{Q}},F)\tag{$2'_{\mathfrak X}$}\end{gather}with ${\mathbb{R}}\Gamma_{dR}({\mathfrak X}'_{\mathbb{Q}},F)={\mathbb{R}}\Gamma({\mathfrak X}'_{\mathbb{Q}},(F\otimes_K\Omega^{\bullet}_{{\mathfrak X}'_{\mathbb{Q}}},1\otimes d))$. Now suppose in addition that {\it all irreducible components of $Y$ are proper over $k$}. Then the tubes $]Y^{p+1}[_{\mathfrak X}$ are partially proper, so \cite{crelle} 3.2 (applied to the morphism between the respective Hodge-de Rham spectral sequences) tells us$${\mathbb{R}}\Gamma_{dR}(]Y^{j}[_{\mathfrak X})\cong {\mathbb{R}}\Gamma_{dR}(]Y^{j}[_{{\mathfrak X}'}),\quad\quad {\mathbb{R}}\Gamma_{dR}({\mathfrak X}_{\mathbb{Q}},F)\cong {\mathbb{R}}\Gamma_{dR}({\mathfrak X}'_{\mathbb{Q}},F).$$On the other hand, \ref{corit} and \ref{maint} tell us$${\mathbb{R}}\Gamma_{rig}(Y^j/{\mathfrak S}^{0} )\cong {\mathbb{R}}\Gamma_{conv}(Y^j/{\mathfrak S}^{0} ),\quad\quad {\mathbb{R}}\Gamma_{rig}(Y/{\mathfrak S}^{0},E )\cong {\mathbb{R}}\Gamma_{conv}(Y/{\mathfrak S}^{0},E).$$

\begin{kor}\label{rigss} There are canonical isomorphisms $(1')\cong (1)\cong (2_{\mathfrak X})\cong (2'_{\mathfrak X})$ and$${\mathbb{R}}\Gamma_{dR}({\mathfrak X}'_{\mathbb{Q}},F)\cong {\mathbb{R}}\Gamma_{conv}(Y/{\mathfrak S}^{0},E).$$
\end{kor}

\addtocounter{satz}{1}{\bf \arabic{section}.\arabic{satz}}\newcounter{phien1}\newcounter{phien2}\setcounter{phien1}{\value{section}}\setcounter{phien2}{\value{satz}} The interest in \ref{rigs}, \ref{rigss} is that the spectral sequence $(1)$ (or $(1)'$) carries additional structure {\it a priori} not present in $(2_{\mathfrak X})$ (or $(2'_{\mathfrak X})$). Namely (see section \ref{fromo}):\\(i) Each term $H^*_{rig}(Y^{m}/{\mathfrak S}^{0},E)$ is finite dimensional if $Y$ is quasi-compact; if $Y$ is proper also $H^{*}_{rig}(Y/{\mathfrak S}^{0},E)$ is finite dimensional.\\(ii) There is a natural $K_0$-linear monodromy operator $N$ on $H^{*}_{rig}(Y/{\mathfrak S}^{0},E)$.\\(iii) If $q, \sigma$ are as in \arabic{fro1}.\arabic{fro2} and if $(E,\phi)\in F\mbox{-}{\cal LS}(Y,K_0)$ then we have the $\sigma$-linear endomorphism $\phi$ acting on $(1)$. On $H^*_{rig}(Y/{\mathfrak S}^{0},E)$ it satisfies $N\phi=q\phi N$; in particular, $N$ is nilpotent if $H^*_{rig}(Y/{\mathfrak S}^{0},E)$ is finite dimensional. If $k$ is perfect and $A_0=W(k)$, $\phi$ is bijective on each term of $(1)$. If $k$ is finite and if for any component intersection $M$ of $Y$ the eigenvalues of $\phi$ acting on the constant sheaf $E|_M$ are Weil numbers, then each term of $(1)$ is a mixed $F$-isocrystal: the eigenvalues of $\phi$ are Weil numbers.\\For the interaction between $N$ and the spectral sequence $(1)$ see \cite{mono}.

\addtocounter{satz}{1}{\bf \arabic{section}.\arabic{satz}}\newcounter{cryrig1}\newcounter{cryrig2}\setcounter{cryrig1}{\value{section}}\setcounter{cryrig2}{\value{satz}} Suppose $k$ is perfect, $A_0=W(k)$ and that even $Y$ is proper. Let $E$ be the constant sheaf $K$. Then ${\mathbb{R}}\Gamma_{conv}(Y/{\mathfrak S}^{0} )$ can be naturally identified (similarly to \cite{berfi} 1.9) with ${\mathbb{R}}\Gamma_{crys}(Y/{\mathfrak S}^{0} )\otimes_{A_0}K_0$, where ${\mathbb{R}}\Gamma_{crys}(Y/{\mathfrak S}^{0})$ is the logarithmic crystalline cohomology of $Y/{\mathfrak S}^{0}$. In \cite{hyoka}, the isomorphism $H^*_{dR}({\mathfrak X}'_{\mathbb{Q}})\cong H^*_{crys}(Y/{\mathfrak S}^{0} )\otimes_{A_0}K$ is constructed by crystalline methods. In case $K=K_0$ the method is to find a canonical section $s$ --- roughly by requiring that it be Frobenius compatible --- of $$H^*_{crys}(Y/\widetilde{S}^{DP})\otimes_{A_0}K_0\longrightarrow H^*_{crys}(Y/{\mathfrak S}^{0} )\otimes_{A_0}K_0,$$ and then to compose with$$H^*_{crys}(Y/\widetilde{S}^{DP})\otimes_{A_0}K_0\longrightarrow H^*_{crys}(Y/{\mathfrak S}^{\pi} )\otimes_{A_0}K_0=H^*_{dR}({\mathfrak X}'_{\mathbb{Q}}).$$ Here $\widetilde{S}^{DP}$ is the DP-envelope of $(t)$ in $\widetilde{S}$. In our approach, the section $s$ corresponds to the composite$$H^*_{conv}(M_{\bullet}/{\mathfrak S}^{0})\cong H^*_{rig}(M_{\bullet}/{\mathfrak S}^{0})\stackrel{(\beta^0)^{-1}}{\longrightarrow}H^*_{rig}((P_{\bullet},V_{\bullet})/{\mathfrak S})\longrightarrow H^*_{crys}(Y/\widetilde{S}^{DP})\otimes_{A_0}K_0.$$  

\section{The proof of Theorem \ref{hkrigis}}
\label{secreshea}
\addtocounter{satz}{1}{\bf \arabic{section}.\arabic{satz}}\newcounter{locset1}\newcounter{locset2}\setcounter{locset1}{\value{section}}\setcounter{locset2}{\value{satz}} We fix $\emptyset\ne J'\subset J\subset \Upsilon$ and assume the situation of \arabic{locemb1}.\arabic{locemb2}. We often drop subscripts $M$, i.e. $P^{J'}=P^{J'}_M$, $\widetilde{P}=\widetilde{P}_{M}$, $\widetilde{V}=\widetilde{V}_{M}$, $\widetilde{D}=\widetilde{D}_{M}$. Let ${\mathfrak M}$, resp. ${\mathfrak V}$, resp. ${\mathfrak P}$, be the weak completion of $\widetilde{M}$, resp. of $\widetilde{V}$, resp. of $\widetilde{P}$, and let ${\mathfrak V}^a={\mathfrak V}\times_{\mathfrak S}{\mathfrak S}^a$ for $a\in\{0,\pi\}$. Let ${\mathfrak M}_{{}\mathbb{Q}}$ (resp. ${\mathfrak P}_{\mathbb{Q}}$) be the generic fibre (as $K_0$-dagger space) of ${\mathfrak M}_{}$ (resp. of ${\mathfrak P}$). Write ${\nabla}^{0}_{{}}=]M_{}[_{{\mathfrak V}^0_{}}$. In$${\mathfrak M}_{{}\mathbb{Q}}\to{\nabla}^{0}_{{}}\to]M[_{\mathfrak P}\to{\mathfrak P}_{\mathbb{Q}}$$the first two arrows are closed immersions, the last one is an open immersion. Let $\Omega^{\bullet}_{\widetilde{P}/\widetilde{W}}$ denote the log de Rham complex of $\widetilde{P}\to\widetilde{W}$ and let $$\Omega^{\bullet}_{\widetilde{P}}=\frac{\Omega^{\bullet}_{\widetilde{P}/\widetilde{W}}}{\Omega^{\bullet-1}_{\widetilde{P}/\widetilde{W}}\wedge\dlog(t)}.$$Weak completion and tensoring with $\mathbb{Q}$ gives a sheaf complex $\Omega^{\bullet}_{\mathfrak{P},\mathbb{Q}}$ on ${\mathfrak P}_{\mathbb{Q}}$.

\begin{pro}\label{anarighomot} The canonical map $${\mathbb{R}}\Gamma_{rig}(M_{}/{\mathfrak S}^{0} )={\mathbb{R}}\Gamma({{\nabla}^{0}},(\Omega^{\bullet}_{\mathfrak{P},\mathbb{Q}}|_{]M[_{\mathfrak P}})\otimes{\mathcal O}_{{\nabla}^{0}})\longrightarrow {\mathbb{R}}\Gamma({\mathfrak M}_{{}\mathbb{Q}},(\Omega^{\bullet}_{\mathfrak{P},\mathbb{Q}}|_{]M[_{\mathfrak P}})\otimes{\cal O}_{{\mathfrak M}_{{}\mathbb{Q}}})$$is an isomorphism.
\end{pro}

{\sc Proof:} {\it Step 1:} Let $\widetilde{B}_{}[x_j]^{\dagger}_{j\in {J}}$ be the weak completion of $\widetilde{B}_{}[x_j]_{j\in {J}}$, i.e. ${\mathfrak V}=\spwf(\widetilde{B}_{}[x_j]^{\dagger}_{j\in {J}})$. For subsets $G\subset\Upsilon=J\sqcup\widehat{J}$ define the element$$t_G=(\prod_{j\in G\cap \widehat{J}}\widetilde{\nu}_{{}j})(\prod_{j\in G\cap J}x_j)\in \widetilde{B}_{}[x_j]_{j\in {J}}$$and the subset$$T_G=\{\widetilde{\nu}_{{}j};\,j\in G\cap \widehat{J}\}\cup\{x_j;\,j\in G\cap J\}\subset\widetilde{B}_{}[x_j]_{j\in {J}}.$$With our previous identification $x_j=\frac{z_{j,1}}{z_{j,0}}$ for $j\in {J}$, we have $t=t_{\Upsilon}$ and$${\nabla}^{0}_{{}}=\{x\in\spm(\frac{\widetilde{B}_{}[x_j]^{\dagger}_{j\in {J}}}{(t_{\Upsilon})}\otimes\mathbb{Q});\quad |x_j(x)|<1\mbox{ for all }j\in {J}\}.$$ For $\rho\in|K_0^{\times}|\otimes\mathbb{Q}, \rho<1$, define its admissible open $K_0$-dagger subspace $${\nabla}^{0}_{{}\rho}=\{x\in\spm(\frac{\widetilde{B}_{}[x_j]^{\dagger}_{j\in {J}}}{(t_{\Upsilon})}\otimes\mathbb{Q});\quad |x_j(x)|\le\rho\mbox{ for all }j\in {J}\}.$$
Then ${\nabla}^{0}_{{}}=\cup_{\rho<1}{\nabla}^{0}_{{}\rho}$ is an admissible open covering and it is enough to show that$${\mathbb{R}}\Gamma({{\nabla}^{0}_{\rho}},(\Omega^{\bullet}_{\mathfrak{P},\mathbb{Q}}|_{]M[_{\mathfrak P}})\otimes{\mathcal O}_{{\nabla}^{0}})\longrightarrow{\mathbb{R}}\Gamma({\mathfrak M}_{{}\mathbb{Q}},(\Omega^{\bullet}_{\mathfrak{P},\mathbb{Q}}|_{]M[_{\mathfrak P}})\otimes{\cal O}_{{\mathfrak M}_{{}\mathbb{Q}}})$$is an isomorphism for each such $\rho$ (the argument for this is like in step 1 of the proof of \ref{nochiso}). We fix such a $\rho$. Since ${\nabla}^{0}_{{}\rho}$ is affinoid, we have$$H^m({{\nabla}^{0}_{\rho}},(\Omega^{l}_{\mathfrak{P},\mathbb{Q}}|_{]M[_{\mathfrak P}})\otimes{\mathcal O}_{{\nabla}^{0}})=0=H^m({\mathfrak M}_{{}\mathbb{Q}},(\Omega^{l}_{\mathfrak{P},\mathbb{Q}}|_{]M[_{\mathfrak P}})\otimes{\cal O}_{{\mathfrak M}_{{}\mathbb{Q}}})$$for all $m>0,$ $l\ge0$, see \cite{crelle}. Each term of $\Omega^{\bullet}=\Gamma({{\nabla}^{0}_{\rho}},(\Omega^{\bullet}_{\mathfrak{P},\mathbb{Q}}|_{]M[_{\mathfrak P}})\otimes{\mathcal O}_{{\nabla}^{0}})$ is locally free over $R=\Omega^0=\Gamma({{\nabla}^{0}_{\rho}},{\mathcal O}_{{{\nabla}^{0}_{\rho}}})$. We have a natural map $g:\widetilde{B}_{}[x_j]_{j\in {J}}\to R$ and if for subsets $T\subset\widetilde{B}_{}[x_j]_{j\in {J}}$ we set$$R_{T}=\frac{R}{(\{g(r);\, r\in T\})}$$ then we need to show that $\Omega^{\bullet}\to\Omega^{\bullet}\otimes_RR_{T_J}$ is a quasi-isomorphism.\\{\it Step 2:} All the following tensor products are taken over $R$. We claim that for all $\emptyset\ne J_1\subset J_2\subset \Upsilon$ with $J_2-J_1\subset {J}$ the canonical map$$\mu:\Omega^{\bullet}\otimes R_{T_{J_1}}\stackrel{}{\longrightarrow}\Omega^{\bullet}\otimes R_{T_{J_2}}$$is a quasi-isomorphism. By induction we may suppose $J_2=J_1\cup\{j_0\}$ for some $j_0\in {J}$ with $j_0\notin J_1$. The ideal in $\widetilde{B}_{}[x_j]_{j\in {J}}$ generated by $T_{J_1}$ automatically contains $t_{\Upsilon}$ (because $\emptyset\ne J_1$); therefore, setting $L=J_1\cap \widehat{J}=J_2\cap \widehat{J}$ we find that $R_{T_{J_m}}$ for $m\in\{1,2\}$ is the ring of global sections on the $K_0$-dagger space$$\{x\in\spm(\frac{\widetilde{B}_{}}{(\{\widetilde{\nu}_{{}j};\,j\in L\})}[x_j]^{\dagger}_{j\in {J-(J_m\cap J)}});\quad |x_j(x)|\le\rho\mbox{ for all }j\in {J-(J_m\cap J)}\}.$$In particular, $R_{T_{J_1}}$ is the ring of overconvergent functions on the relative closed disk of radius $\rho$ over $R_{T_{J_2}}$ (with coordinate $x_{j_0}$), and the map $\mu$ has a natural section. We prove that this section induces surjective maps in cohomology. Denote again by ${\dlog(x_{j_0})}$ the class of $\dlog(x_{j_0})$ in $\Omega^1\otimes R_{T_{J_2}}$. An easy consideration with local coordinates (see for example \cite{colo}, proof of Theorem 3.14) shows that, at least after localization on $M$, we may choose a complement $\Omega_c^1$ of the $R_{T_{J_2}}$-submodule $<{\dlog(x_{j_0})}>$ of $\Omega^1\otimes R_{T_{J_2}}$ generated by ${\dlog(x_{j_0})}$ such that $\Omega_c^{\bullet}$, the sub-$R_{T_{J_2}}$-algebra of $\Omega^{\bullet}\otimes R_{T_{J_2}}$ generated by $\Omega_c^1$, is stable for the differential $d$. Denoting $<{\dlog(x_{j_0})}>^{\bullet}$ the sub-$R_{T_{J_2}}$-algebra of $\Omega^{\bullet}\otimes R_{T_{J_2}}$ generated by $<{\dlog(x_{j_0})}>$ we have$$\Omega^{\bullet}\otimes R_{T_{J_2}}=\Omega_c^{\bullet}\otimes<{\dlog(x_{j_0})}>^{\bullet}.$$Let $\omega\in\Omega^k\otimes R_{T_{J_1}}$. 
It can be written as$$\omega=\sum_{\lambda\ge 0}a_{\lambda}x^\lambda_{j_0}{\dlog(x_{j_0})}+\sum_{\lambda\ge0}b_{\lambda}x^{\lambda}_{j_0}$$with $a_{\lambda}\in\Omega_c^{k-1}$ and $b_{\lambda}\in\Omega_c^k$. Subtracting $d(\sum_{\lambda>0}\lambda^{-1}a_{\lambda}x^{\lambda}_{j_0})$ (the expression $\sum_{\lambda>0}\lambda^{-1}a_{\lambda}x^{\lambda}_{j_0}$ exists by overconvergence!) and renaming the coefficients we may write $\omega$ modulo exact forms as$$\omega=a_0{\dlog(x_{j_0})}+\sum_{\lambda\ge0}b_{\lambda}x^{\lambda}_{j_0}.$$If $d\omega=0$ we get $\omega=a_0{\dlog(x_{j_0})}+b_0$ which is an element of $\Omega^k\otimes R_{T_{J_2}}$ and the claim follows.\\{\it Step 3:} Now we can prove the statement to which we reduced the proposition in step 1 by entirely formal reasoning (which we also applied in a similar situation in \cite{colo} Theorem 3.14). We will show that in$$\Omega^{\bullet}=\Omega^{\bullet}\otimes R\stackrel{\alpha}{\longrightarrow}\Omega^{\bullet}\otimes R_{\{t_J\}}\stackrel{\beta}{\longrightarrow}\Omega^{\bullet}\otimes R_{T_J}$$both $\alpha$ and $\beta$ are quasi-isomorphisms. The exact sequences
$$0\longrightarrow R\longrightarrow R_{\{t_J\}}\bigoplus R_{\{t_{\widehat{J}}\}}\longrightarrow R_{\{t_J,t_{\widehat{J}}\}}\longrightarrow0$$$$0\longrightarrow R_{\{t_J\}}\longrightarrow R_{\{t_J\}}\bigoplus R_{\{t_J,t_{\widehat{J}}\}} \longrightarrow R_{\{t_J,t_{\widehat{J}}\}}\longrightarrow0$$ show that, to prove that $\alpha$ is a quasi-isomorphism, it is enough to show that $\Omega^{\bullet}\otimes R_{\{t_{\widehat{J}}\}}\to\Omega^{\bullet}\otimes R_{\{t_J,t_{\widehat{J}}\}}$ is a quasi-isomorphism. To see this, it is enough to show that both $\Omega^{\bullet}\otimes R_{\{t_{\widehat{J}}\}}\stackrel{\gamma}\to\Omega^{\bullet}\otimes R_{\{t_{\widehat{J}}\}\cup T_J}$ and $\Omega^{\bullet}\otimes R_{\{t_J,t_{\widehat{J}}\}}\stackrel{\delta}\to\Omega^{\bullet}\otimes R_{\{t_{\widehat{J}}\}\cup T_J}$ are quasi-isomorphisms. Consider the exact sequence\begin{gather}0\longrightarrow R_{\{t_{\widehat{J}}\}}\longrightarrow\bigoplus_{j\in \widehat{J}}R_{\{\widetilde{\nu}_{{}j}\}}\longrightarrow\bigoplus_{G\subset \widehat{J}\atop |G|=2}R_{T_{G}}\longrightarrow\ldots\longrightarrow R_{T_{\widehat{J}}}\longrightarrow 0\tag{$*$}\end{gather}Comparison of the exact sequences $(*)\otimes\Omega^{\bullet}$ and $(*)\otimes R_{\{t_{\widehat{J}}\}\cup T_J}\otimes\Omega^{\bullet}$ shows that to prove that $\gamma$ is a quasi-isomorphism, it is enough to show this for $\Omega^{\bullet}\otimes R_{T_{G}}\to\Omega^{\bullet}\otimes R_{T_{G\cup J}}$ for all $\emptyset\ne G\subset \widehat{J}$; but this has been done in step 2. Comparison of $(*)\otimes R_{\{t_J,t_{\widehat{J}}\}}\otimes\Omega^{\bullet}$ and $(*)\otimes R_{\{t_{\widehat{J}}\}\cup T_J}\otimes\Omega^{\bullet}$ shows that to prove that $\delta$ is a quasi-isomorphism, it is enough to show this for $\Omega^{\bullet}\otimes R_{\{t_J\}\cup T_G}\stackrel{\epsilon_G}{\to}\Omega^{\bullet}\otimes R_{T_J\cup T_G}$
 for all $\emptyset\ne G\subset \widehat{J}$. Consider the exact sequence \begin{gather}0\longrightarrow R_{\{t_J\}}\longrightarrow\bigoplus_{j\in {J}}R_{\{x_j\}}\longrightarrow \bigoplus_{F\subset J\atop |F|=2}R_{T_F}\longrightarrow\ldots\longrightarrow R_{T_{J}}\longrightarrow0\tag{$**$}\end{gather} The exact sequence $(**)\otimes R_{\{t_J\}\cup T_G}\otimes\Omega^{\bullet}$ shows that to prove that $\epsilon_G$ is a quasi-isomorphism, it is enough to show this for $\Omega^{\bullet}\otimes R_{T_{G\cup F}}\to\Omega^{\bullet}\otimes R_{T_{G\cup J}}$ for all $\emptyset\ne F\subset {J}$; but this has been done in step 2. The exact sequence $(**)\otimes\Omega^{\bullet}$ shows that to prove that $\beta$ is a quasi-isomorphism, it is enough to show this for $\Omega^{\bullet}\otimes R_{T_F}\to\Omega^{\bullet}\otimes R_{T_J}$ for all $\emptyset\ne F\subset {J}$; but this has been done in step 2.

\addtocounter{satz}{1}{\bf \arabic{section}.\arabic{satz}}\newcounter{proj1}\newcounter{proj2}\setcounter{proj1}{\value{section}}\setcounter{proj2}{\value{satz}} For a subset $I\subset J$ endow $({\mathbb P}_{A_0}^1)^I=\times_{\widetilde{W}}(\proj (A_0[z_{j,0},z_{j,1}]))_{j\in I}$ with the product log structure, where the individual factors carry the log structure defined by the divisor $\{0,\infty\}$, the zero and the pole of $x_j=\frac{z_{j,1}}{z_{j,0}}$. We always think of ${({\mathbb P}_{A_0}^1)}^I$ as a factor of $\widetilde{P}=\widetilde{M}\times(\mathbb{P}_{A_0}^1)^J$, i.e. as endowed with the projection $p:\widetilde{P}\to{({\mathbb P}_{A_0}^1)}^I$. (Throughout this section, we denote each projection map from a fibre product to any of its factors by $p$, and $p^*$ always denotes module theoretic pull back of structure sheaf modules). We denote the weak completion of $({\mathbb P}_{A_0}^1)^I$ again by $({\mathbb P}_{A_0}^1)^I$, and we view the generic fibre $({\mathbb P}_{K_0}^1)^I$ as a dagger space. Denote by $\widetilde{M}_{}^{\sharp}$ the log scheme with underlying scheme $\widetilde{M}_{}$ and with the log structure defined by the (relative to $\widetilde{W}$) normal crossings divisor ${\bf V}(\prod_{j\in \widehat{J}}\widetilde{\nu}_{{}j})=\widetilde{D}_{{}}\cap\widetilde{M}_{{}}$ on $\widetilde{M}_{}$. Let ${{\Omega}}^{\bullet}_{({\mathbb P}_{A_0}^1)^I}$ (resp. ${{\Omega}}^{\bullet}_{\widetilde{M}^{\sharp}_{}}$) be the logarithmic de Rham complex of $({\mathbb P}_{A_0}^1)^I\to\widetilde{W}$ (resp. of $\widetilde{M}_{}^{\sharp}\to{\widetilde{W}}$), and let ${{\Omega}}^{\bullet}_{({\mathbb P}_{K_0}^1)^I}$ (resp. $\Omega^{\bullet}_{{\mathfrak M}^{\sharp}_{{}\mathbb{Q}}}$) be the complex on $({\mathbb P}_{K_0}^1)^I$ (resp. on ${\mathfrak M}_{{}\mathbb{Q}}$) obtained by weak completion and tensoring with $\mathbb{Q}$. Thus ${{\Omega}}^{1}_{({\mathbb P}_{A_0}^1)^I}$ is free of rank $|I|$, a basis is $\{\dlog(x_{j});\,j\in I\}$. We have a canonical direct sum decomposition$${{{\Omega}}}^{1}_{\widetilde{P}_{}/\widetilde{W}}=p^*{{\Omega}}^1_{\widetilde{M}_{}^{\sharp}}\bigoplus p^*{{\Omega}}^{1}_{({\mathbb P}_{A_0}^1)^J}.$$Now fix $\emptyset\ne{J'}\subset {J}$ and $j'\in{J'}$, let $J_0={J}-\{j'\}$. Then the submodule $p^*{{\Omega}}^1_{\widetilde{M}_{}^{\sharp}}\oplus p^*{{\Omega}}^{1}_{({\mathbb P}_{A_0}^1)^{J_0}}$ of the right hand side maps isomorphically to the quotient ${{{\Omega}}}^{1}_{\widetilde{P}_{}}={{{\Omega}}}^{1}_{\widetilde{P}_{}/\widetilde{W}}/\dlog(t)$ of ${{{\Omega}}}^{1}_{\widetilde{P}_{}/\widetilde{W}}$ (note $\dlog(t)=\sum_{j\in {J}}\dlog(x_{j})+\sum_{j\in \widehat{J}}\dlog(\widetilde{\nu}_{{}j})$), i.e.$${{{\Omega}}}^{1}_{\widetilde{P}_{}}\cong p^*{{\Omega}}^1_{\widetilde{M}_{}^{\sharp}}\bigoplus p^*{{\Omega}}^1_{({\mathbb P}_{A_0}^1)^{J_0}}$$\begin{gather}{{\Omega}}^{\bullet}_{\widetilde{P}_{{}}}\cong
p^*{{\Omega}}^{\bullet}_{\widetilde{M}_{}^{\sharp}}\otimes_{{\cal O}_{\widetilde{P}_{{}}}} p^*{{\Omega}}^{\bullet}_{({\mathbb P}_{A_0}^1)^{J_0}}
\tag{$*$}\end{gather}

\begin{lem}\label{herzres} For $M^{\heartsuit}=M-(M\cap\cup_{j\in\widehat{J}}Y_j)$ and $a\in\{0,\pi\}$ the restriction maps \begin{gather} {\mathbb{R}}\Gamma_{rig}(M/{\mathfrak S}^{a} )\longrightarrow {\mathbb{R}}\Gamma_{rig}(M^{\heartsuit}/{\mathfrak S}^{a} )\tag{$1_a$}\\ {\mathbb{R}}\Gamma_{rig}((P^{{J'}}_M,V^{{J'}}_M)/{\mathfrak S})\longrightarrow {\mathbb{R}}\Gamma_{rig}((P^{{J'}}_{M^{\heartsuit}},V^{{J'}}_{M^{\heartsuit}})/{\mathfrak S})\tag{$2$}\end{gather}are isomorphisms.
\end{lem}

{\sc Proof:} In all three cases this is due to overconvergence: the corresponding statements with ${\mathbb{R}}\Gamma_{conv}$ instead of ${\mathbb{R}}\Gamma_{rig}$ are false. In general, for a strictly semistable weak formal scheme ${\mathfrak Y}$, if $M$ is the intersection of some of the irreducible components of its reduction and if $M^{\heartsuit}$ is the maximal open subscheme of $M$ which has empty intersection with all the other components, then the restriction map $H_{dR}^*(]M[_{\mathfrak Y})\to H_{dR}^*(]M^{\heartsuit}[_{\mathfrak Y})$ is an isomorphism, see \cite{findag} Theorem 2.3. Applied to ${\mathfrak Y}={\mathfrak V}^{\pi}$ we get $H^*_{rig}(M/{\mathfrak S}^{\pi})=H_{dR}^*(]M[_{{\mathfrak V}^{\pi}})=H_{dR}^*(]M^{\heartsuit}[_{{\mathfrak V}^{\pi}})=H^*_{rig}(M^{\heartsuit}/{\mathfrak S}^{\pi})$, i.e. that $(1_{\pi})$ is an isomorphism. Now consider $(2)$. Weak completion and tensoring with $\mathbb{Q}$ the decomposition $(*)$ of ${{\Omega}}^{\bullet}_{\widetilde{P}_{{}}}$ in \arabic{proj1}.\arabic{proj2} gives\begin{gather}\Omega^{\bullet}_{\mathfrak{P},\mathbb{Q}}=p^*\Omega^{\bullet}_{{\mathfrak M}^{\sharp}_{{}\mathbb{Q}}}\otimes p^*{{\Omega}}^{\bullet}_{({\mathbb P}_{K_0}^1)^{J_0}}.\tag{$**$}\end{gather}Let $N_{J_0,J''}$ be the closed subscheme of the reduction modulo $m_{A_0}$ of $({\mathbb P}_{A_0}^1)^{J_0}$ where all $x_j$ for $j\in J''=J-J'$ are defined and vanish. Then $$]P_{M}^{{J'}}[_{{\mathfrak P}_{{}}}={\mathfrak M}_{{}\mathbb{Q}}\times]N_{J_0,J''}[_{({\mathbb P}_{A_0}^1)^{J_0}}\times{\mathbb P}^1_{K_0}$$$$]P_{M^{\heartsuit}}^{{J'}}[_{{\mathfrak P}_{{}}}=]M^{\heartsuit}[_{{\mathfrak M}_{{}}}\times]N_{J_0,J''}[_{({\mathbb P}_{A_0}^1)^{J_0}}\times{\mathbb P}^1_{K_0}.$$Since $\widetilde{D}_{{}}\cap\widetilde{M}_{{}}\to\widetilde{M}_{}$ reduces to the embedding $(M-M^{\heartsuit})\to M$, it is immediate from the proofs of the comparison theorems in \cite{baci}, \cite{kidr} that the restriction map$${\mathbb{R}}\Gamma({\mathfrak M}_{{}\mathbb{Q}},\Omega^{\bullet}_{{\mathfrak M}^{\sharp}_{{}\mathbb{Q}}})\longrightarrow {\mathbb{R}}\Gamma(]M^{\heartsuit}[_{{\mathfrak M}_{{}}},\Omega^{\bullet}_{{\mathfrak M}^{\sharp}_{{}\mathbb{Q}}})$$is an isomorphism. In view of the above decompositions it follows from K\"unneth formulas (which in this case are easily proved, similarly to those in \cite{findag}; see also \ref{nochiso} below) that the restriction map$${\mathbb{R}}\Gamma(]P_M^{{J'}}[_{{\mathfrak P}_{{}}},\Omega^{\bullet}_{\mathfrak{P},\mathbb{Q}})\longrightarrow {\mathbb{R}}\Gamma(]P_{M^{\heartsuit}}^{{J'}}[_{{\mathfrak P}_{{}}},\Omega^{\bullet}_{\mathfrak{P},\mathbb{Q}})$$is an isomorphism. This proves $(2)$. Finally, by \ref{anarighomot} applied to both $M$ and $M^{\heartsuit}$ we see that to prove that $(1_0)$ is an isomorphism we only have to show that of$${\mathbb{R}}\Gamma({\mathfrak M}_{{}\mathbb{Q}},(\Omega^{\bullet}_{\mathfrak{P},\mathbb{Q}}|_{]M[_{\mathfrak P}})\otimes{\cal O}_{{\mathfrak M}_{{}\mathbb{Q}}})\longrightarrow{\mathbb{R}}\Gamma(]M^{\heartsuit}[_{{\mathfrak M}_{}},(\Omega^{\bullet}_{\mathfrak{P},\mathbb{Q}}|_{]M[_{\mathfrak P}})\otimes{\cal O}_{{\mathfrak M}_{{}\mathbb{Q}}}).$$ But $(\Omega^{\bullet}_{\mathfrak{P},\mathbb{Q}}|_{]M[_{\mathfrak P}})\otimes{\cal O}_{{\mathfrak M}_{{}\mathbb{Q}}}$ decomposes according to the decomposition $(**)$ of $\Omega^{\bullet}_{\mathfrak{P},\mathbb{Q}}|_{]M[_{\mathfrak P}}$. Thus the argument which proved that $(2)$ is an isomorphism works again.

\addtocounter{satz}{1}{\bf \arabic{section}.\arabic{satz}} We keep the setting and notations from \arabic{locset1}.\arabic{locset2} and assume in addition $M\cap \cup_{j\in\widehat{J}}Y_j=\emptyset$, hence ${J}=\Upsilon$. Fix an ordering of $J$. For a subset ${{I}}$ of ${J}$ let $H_{A_0,{{I}}}^0=A_0$ and for $s\ge 1$ let $H_{A_0,{I}}^s$ be the free $A_0$-module with basis the set of symbols$$\dlog(x_{j_1})\wedge\ldots\wedge\dlog(x_{j_s})$$with $j_r\in {{I}}$ and $j_1<\ldots<j_s$, and for $s\ge0$ let $H_{{I}}^s=H_{A_0,{{I}}}^s\otimes_{A_0}K_0$. Fix a non-empty subset ${J'}\subset {J}$, an element $j'\in{J'}$ and let ${J_0}={J}-\{j'\}$. Denote by $H_{dR}^{*}({\mathfrak M}_{{}\mathbb{Q}})$ the (non-logarithmic) de Rham cohomology of the dagger space ${\mathfrak M}_{{}\mathbb{Q}}$.

\begin{lem}\label{nochiso} For each $m\in\mathbb{Z}$ we have a canonical isomorphism$$H^m_{rig}((P^{{J'}}_M,V^{{J'}}_M)/{\mathfrak S})=H^m(]P_M^{J'}[_{{\mathfrak P}},\Omega^{\bullet}_{\mathfrak{P},\mathbb{Q}})=\underset{m=m_1+m_2}\bigoplus H_{dR}^{m_1}({\mathfrak M}_{{}\mathbb{Q}})\otimes_{K_0}H_{{J_0}}^{m_2}.$$
\end{lem}

{\sc Proof:} {\it Step 1:} For a subset $I\subset J$ we identify $\spec(A_0[x_j]_{j\in{{I}}})$ with the open affine space $({\mathbb A}_{A_0}^1)^{I}$ in $({\mathbb P}_{A_0}^1)^I$ (cf. \arabic{proj1}.\arabic{proj2}). For $\delta\in|K_0^{\times}|\otimes\mathbb{Q}$ let ${\mathbb D}^{I}_{\delta}$ be the affinoid open subspace of $({\mathbb A}_{A_0}^1)^{I}\subset({\mathbb P}_{A_0}^1)^{I}$ where $|x_j|\le\delta$ for all $j\in I$. Let ${J''}={J}-{J'}$ and define $${\nabla}^{{J'}}_{{}\epsilon}={\mathfrak M}_{{}\mathbb{Q}}\times {\mathbb D}^{J''}_{\epsilon}\times{\mathbb D}^{J'}_{1} .$$Let $\{\epsilon_n\}_{n\in{\mathbb N}}$ with $\epsilon_n\in|K_0^{\times}|\otimes\mathbb{Q}$, $\epsilon_n<1$, be an increasing sequence tending to $1$. The covering by $K_0$-dagger subspaces$$]P_M^{J'}[_{{\mathfrak P}}=\bigcup_{n}{\nabla}^{{J'}}_{{}\epsilon_n}$$is admissible open, hence a spectral sequence$$E_2^{m',m}={\mathbb R}^{m'}\lim_{\gets\atop n}H^{m}({\nabla}^{{J'}}_{{}\epsilon_n},\Omega^{\bullet}_{\mathfrak{P},\mathbb{Q}})\Longrightarrow H^{m+m'}(]P_M^{J'}[_{{\mathfrak P}},\Omega^{\bullet}_{\mathfrak{P},\mathbb{Q}}).$$We will show$$H^m({\nabla}^{{J'}}_{{}\epsilon},\Omega^{\bullet}_{\mathfrak{P},\mathbb{Q}})=\bigoplus_{m=m_1+m_2}H_{dR}^{m_1}({\mathfrak M}_{{}\mathbb{Q}})\otimes_{K_0}H_{{J_0}}^{m_2}$$for arbitrary $\epsilon\in|K_0^{\times}|\otimes\mathbb{Q}$, $\epsilon<1$. Granted this it follows in particular that all transition maps$$H^m({\nabla}^{{J'}}_{{}\epsilon_n},\Omega^{\bullet}_{\mathfrak{P},\mathbb{Q}})\longrightarrow H^m({\nabla}^{{J'}}_{{}\epsilon_{n'}},\Omega^{\bullet}_{\mathfrak{P},\mathbb{Q}})$$ for $n\ge{n'}$ are isomorphisms, hence ${\mathbb R}^{m'}\lim_{\gets\atop n}H^{m}({\nabla}^{{J'}}_{{}\epsilon_n},\Omega^{\bullet}_{\mathfrak{P},\mathbb{Q}})=0$ for all $m'\ne 0$, hence $$H^m(]P_M^{J'}[_{{\mathfrak P}},\Omega^{\bullet}_{\mathfrak{P},\mathbb{Q}})=\lim_{\gets\atop n}H^{m}({\nabla}^{{J'}}_{{}\epsilon_n},\Omega^{\bullet}_{\mathfrak{P},\mathbb{Q}})=H^{m}({\nabla}^{{J'}}_{{}\epsilon},\Omega^{\bullet}_{\mathfrak{P},\mathbb{Q}})$$for an arbitrary $\epsilon\in|K_0^{\times}|\otimes\mathbb{Q}$, $\epsilon<1$, and the Lemma follows. We now fix such an $\epsilon$.\\{\it Step 2:} Here we compute the coherent sheaf cohomology $H^n({\nabla}^{{J'}}_{{}\epsilon},\Omega^{m}_{\mathfrak{P},\mathbb{Q}})$. For this we choose finite type $A_0$-scheme models for ${\nabla}^{{J'}}_{{}\epsilon}$ and ${\mathbb D}^{J''}_{\epsilon}$ and use the scheme theoretic K\"unneth formula. Choose $r\in\mathbb{N}$ and $\alpha\in A_0$ such that $|\alpha|=\epsilon^r$. The scheme $$\widetilde{\mathbb D}^{J''}_{\epsilon}=\spec(\frac{A_0[x_j,y_j]_{j\in {J''}}}{(\alpha.y_j-x_j^r)_{j\in {J''}}})$$comes with an obvious map $\iota_{\epsilon}$ to $({\mathbb P}_{A_0}^1)^{J''}$. Let ${{\Omega}}^{\bullet}_{\widetilde{\mathbb D}^{J''}_{\epsilon}}=\iota_{\epsilon}^*{{\Omega}}^{\bullet}_{({\mathbb P}_{A_0}^1)^{J''}}$, a graded sheaf of coherent ${\cal O}_{\widetilde{\mathbb D}^{J''}_{\epsilon}}$-modules, and then form the graded sheaf of coherent ${\cal O}_{(\widetilde{\mathbb D}^{J''}_{\epsilon})_{\widetilde{M}}}$-modules $\Omega^{\bullet}_{(\widetilde{\mathbb D}^{J''}_{\epsilon})_{\widetilde{M}}}=p^*{{\Omega}}^{\bullet}_{\widetilde{\mathbb D}^{J''}_{\epsilon}}\otimes p^*{{\Omega}}^{\bullet}_{\widetilde{M}_{}^{\sharp}}$ on $(\widetilde{\mathbb D}^{J''}_{\epsilon})_{\widetilde{M}}=\widetilde{M}_{{}}\times\widetilde{\mathbb D}^{J''}_{\epsilon}$. Let ${J'_0}={J'}-\{j'\}$ and define$$\widetilde{P}_{\epsilon}=(\widetilde{\mathbb D}^{J''}_{\epsilon})_{\widetilde{M}}\times ({\mathbb P}_{A_0}^1)^{J'_0}\times ({\mathbb P}_{A_0}^1)^{\{j'\}}.$$From \arabic{proj1}.\arabic{proj2} it follows that the graded sheaf ${{\Omega}}^{\bullet}_{\widetilde{P}_{\epsilon}}=p^{*}{{\Omega}}^{\bullet}_{\widetilde{P}_{{}}}$ on $\widetilde{P}_{\epsilon}$ decomposes as$${{\Omega}}^{\bullet}_{\widetilde{P}_{\epsilon}}=p^*\Omega^{\bullet}_{(\widetilde{\mathbb D}^{J''}_{\epsilon})_{\widetilde{M}}}\otimes p^*{{\Omega}}^{\bullet}_{({\mathbb P}_{A_0}^1)^{J'_0}}.$$Now for the "missing" factor $j'$ we have $H^0(({\mathbb P}_{A_0}^1)^{\{j'\}},{\cal O}_{({\mathbb P}_{A_0}^1)^{\{j'\}}})=A_0$ and $H^n(({\mathbb P}_{A_0}^1)^{\{j'\}},{\cal O}_{({\mathbb P}_{A_0}^1)^{\{j'\}}})=0$ if $n\ne0$. Hence, by the scheme theoretic K\"unneth formula,$$H^n(\widetilde{P}_{\epsilon},{{\Omega}}^{m}_{\widetilde{P}_{\epsilon}})=\underset{m=m_1+m_2\atop n=n_1+n_2}\bigoplus H^{n_1}((\widetilde{\mathbb D}^{J''}_{\epsilon})_{\widetilde{M}},\Omega^{m_1}_{(\widetilde{\mathbb D}^{J''}_{\epsilon})_{\widetilde{M}}})\otimes_{A_0}H^{n_2}(({\mathbb P}_{A_0}^1)^{J'_0},{{\Omega}}^{m_2}_{({\mathbb P}_{A_0}^1)^{J'_0}}).$$For the tensor factors we find: 
$H^{n}((\widetilde{\mathbb D}^{J''}_{\epsilon})_{\widetilde{M}},\Omega^{m}_{(\widetilde{\mathbb D}^{J''}_{\epsilon})_{\widetilde{M}}})=0$ whenever $n>0$ since $(\widetilde{\mathbb D}^{J''}_{\epsilon})_{\widetilde{M}}$ is affine. Furthermore $$H^{n}(({\mathbb P}_{A_0}^1)^{J'_0},{{\Omega}}^{m}_{({\mathbb P}_{A_0}^1)^{J'_0}})=\bigoplus_{n=\sum_{j\in{J'_0}}n_j}\bigoplus_{m=\sum_{j\in{J'_0}}m_j}(\otimes_{A_0}(H^{n_j}(({\mathbb P}_{A_0}^1)^{\{j\}},{{\Omega}}^{m_j}_{({\mathbb P}_{A_0}^1)^{\{j\}}}))_{j\in{J'_0}}).$$
Now ${{\Omega}}^{\bullet}_{({\mathbb P}_{A_0}^1)^{\{j\}}}$ is the logarithmic de Rham complex on ${\mathbb P}^1_{A_0}$ with logarithmic poles at $0$ and $\infty$, thus ${{\Omega}}^{m}_{({\mathbb P}_{A_0}^1)^{\{j\}}}\cong{\cal O}_{{\mathbb P}_{A_0}^1}$ if $m\in\{0,1\}$, and ${{\Omega}}^{m}_{({\mathbb P}_{A_0}^1)^{\{j\}}}=0$ if $m\notin\{0,1\}$. Hence $H^0(({\mathbb P}_{A_0}^1)^{\{j\}},{{\Omega}}^{0}_{({\mathbb P}_{A_0}^1)^{\{j\}}})=A_0$ and  $H^0(({\mathbb P}_{A_0}^1)^{\{j\}},{{\Omega}}^{1}_{({\mathbb P}_{A_0}^1)^{\{j\}}})$ is freely generated by $\dlog(x_j)$, and $H^n(({\mathbb P}_{A_0}^1)^{\{j\}},{{\Omega}}^{m}_{({\mathbb P}_{A_0}^1)^{\{j\}}})=0$ if $n\ne0$ or $m\notin\{0,1\}$. 
Therefore $H^{n}(({\mathbb P}_{A_0}^1)^{J'_0},{{\Omega}}^{m}_{({\mathbb P}_{A_0}^1)^{J'_0}})=0$ if $n>0$, and $H^{0}(({\mathbb P}_{A_0}^1)^{J'_0},{{\Omega}}^{m}_{({\mathbb P}_{A_0}^1)^{J'_0}})=H_{A_0,{J'_0}}^m$. We thus obtain $H^n(\widetilde{P}_{\epsilon},{{\Omega}}^{m}_{\widetilde{P}_{\epsilon}})=0$ if $n>0$, and$$H^0(\widetilde{P}_{\epsilon},{{\Omega}}^{m}_{\widetilde{P}_{\epsilon}})=\underset{m=m_1+m_2}\bigoplus H^{0}((\widetilde{\mathbb D}^{J''}_{\epsilon})_{\widetilde{M}},\Omega^{m_1}_{(\widetilde{\mathbb D}^{J''}_{\epsilon})_{\widetilde{M}}})\otimes_{A_0}H_{A_0,{J'_0}}^{m_2}.$$Let $C=\Gamma((\widetilde{\mathbb D}^{J''}_{\epsilon})_{\widetilde{M}},{\cal O}_{(\widetilde{\mathbb D}^{J''}_{\epsilon})_{\widetilde{M}}})$ and let $C^{\dagger}$ be its weak completion. Let ${\mathfrak P}_{\epsilon}$ (resp. ${\Omega}^{\bullet}_{{\mathfrak P}_{\epsilon}}$) be the weak completion of $\widetilde{P}_{\epsilon}$ (resp. of ${{\Omega}}^{\bullet}_{\widetilde{P}_{\epsilon}}$). By Meredith's GAGA-theorem for weak completion (\cite{mer} p.25) applied to the projective morphism $\widetilde{P}_{\epsilon}\to(\widetilde{\mathbb D}^{J''}_{\epsilon})_{\widetilde{M}}$ we have $$H^n({\mathfrak P}_{\epsilon},\Omega^m_{{\mathfrak P}_{\epsilon}})=H^n(\widetilde{P}_{\epsilon},{{\Omega}}^{m}_{\widetilde{P}_{\epsilon}})\otimes_{C}C^{\dagger}.$$ Since ${\mathfrak P}_{\epsilon}$ (resp. the sheaf ${\Omega}^{m}_{{\mathfrak P}_{\epsilon}}$ on it) is an integral model for ${\nabla}^{{J'}}_{{}\epsilon}$ (resp. for the sheaf $\Omega^{m}_{\mathfrak{P},\mathbb{Q}}|_{{\nabla}^{{J'}}_{{}\epsilon}}$ on it), we obtain$$H^n({\nabla}^{{J'}}_{{}\epsilon},\Omega^{m}_{\mathfrak{P},\mathbb{Q}})=H^n({\mathfrak P}_{\epsilon},\Omega^m_{{\mathfrak P}_{\epsilon}})\otimes\mathbb{Q}=0\quad\quad\mbox{if}\quad n>0,$$ $$H^0({\nabla}^{{J'}}_{{}\epsilon},\Omega^{m}_{\mathfrak{P},\mathbb{Q}})=\bigoplus_{m=m_1+m_2}H^{0}((\widetilde{\mathbb D}^{J''}_{\epsilon})_{\widetilde{M}},\Omega^{m_1}_{(\widetilde{\mathbb D}^{J''}_{\epsilon})_{\widetilde{M}}})\otimes_CC^{\dagger}\otimes_{A_0}H_{{J'_0}}^{m_2}.$$Together this means that $H^m({\nabla}^{{J'}}_{{}\epsilon},\Omega^{\bullet}_{\mathfrak{P},\mathbb{Q}})$ is the $m$-th cohomology group of the complex of $K_0$-vector spaces $$[\bigoplus_{m=m_1+m_2}H^{0}((\widetilde{\mathbb D}^{J''}_{\epsilon})_{\widetilde{M}},\Omega^{m_1}_{(\widetilde{\mathbb D}^{J''}_{\epsilon})_{\widetilde{M}}})\otimes_CC^{\dagger}\otimes_{A_0}H_{{J'_0}}^{m_2}]_{m\in\mathbb{Z}}.$${\it Step 3:} Let $\Omega^{\bullet}_{{\mathbb D}^{J''}_{\epsilon}}$ (resp. $\Omega^{\bullet}_{{\mathfrak M}_{{}\mathbb{Q}}\times {\mathbb D}^{J''}_{\epsilon}}$) be the logarithmic de Rham complex on the affinoid dagger space ${\mathbb D}^{J''}_{\epsilon}$ (resp. ${\mathfrak M}_{{}\mathbb{Q}}\times {\mathbb D}^{J''}_{\epsilon}$) with logarithmic poles along the (respective) divisor ${\bf V}(\prod_{j\in {J''}}x_j)$. As in \cite{findag} section 2 Lemma 3 one easily proves $$H^n({\mathbb D}^{J''}_{\epsilon},\Omega^{\bullet}_{{\mathbb D}^{J''}_{\epsilon}})=H_{{J''}}^n$$and the K\"unneth formula$$H^n({\mathfrak M}_{{}\mathbb{Q}}\times {\mathbb D}^{J''}_{\epsilon},\Omega^{\bullet}_{{\mathfrak M}_{{}\mathbb{Q}}\times {\mathbb D}^{J''}_{\epsilon}})=\bigoplus_{n_1+n_2=n}H_{dR}^{n_1}({\mathfrak M}_{{}\mathbb{Q}})\otimes_{K_0}H_{{J''}}^{n_2}.$$Since ${\mathfrak M}_{{}\mathbb{Q}}\times {\mathbb D}^{J''}_{\epsilon}$ is affinoid this is also the $n$-th cohomology group of the complex $$[H^0({\mathfrak M}_{{}\mathbb{Q}}\times {\mathbb D}^{J''}_{\epsilon},\Omega^{m}_{{\mathfrak M}_{{}\mathbb{Q}}\times {\mathbb D}^{J''}_{\epsilon}})]_{m\in\mathbb{Z}}.$$ By construction we have$$H^0({\mathfrak M}_{{}\mathbb{Q}}\times {\mathbb D}^{J''}_{\epsilon},\Omega^{m}_{{\mathfrak M}_{{}\mathbb{Q}}\times {\mathbb D}^{J''}_{\epsilon}})=H^{0}((\widetilde{\mathbb D}^{J''}_{\epsilon})_{\widetilde{M}},\Omega^{m}_{(\widetilde{\mathbb D}^{J''}_{\epsilon})_{\widetilde{M}}})\otimes_CC^{\dagger}\otimes\mathbb{Q}$$and the proof of the lemma is finished in view of what we saw in step 2.  

\addtocounter{satz}{1}{\bf \arabic{section}.\arabic{satz}}\newcounter{frodec1}\newcounter{frodec2}\setcounter{frodec1}{\value{section}}\setcounter{frodec2}{\value{satz}} Suppose the endomorphism $\phi$ of ${\mathfrak M}$ lifts the $q$-th power Frobenius endomorphism of $M$ for some $q\in p^{\mathbb{N}}$. It induces an endomorphism $\phi_{m,triv}$ of $H_{dR}^{m}({\mathfrak M}_{{}\mathbb{Q}})$. Extend $\phi$ to an endomorphism of ${\mathfrak V}$ by sending $x_j\mapsto x_j^q$. This induces endomorphisms $\phi$ of ${\mathfrak V}^0$ and ${\nabla}^{0}_{{}}$ and thus an endomorphism $\phi_{m,0}$ of 
$H_{rig}^m(M/{\mathfrak S}^0)=H^m({{\nabla}^{0}},(\Omega^{\bullet}_{\mathfrak{P},\mathbb{Q}}|_{]M[_{\mathfrak P}})\otimes{\mathcal O}_{{\nabla}^{0}})$.

\begin{lem}\label{restiso} $H_{dR}^{m}({\mathfrak M}_{{}\mathbb{Q}})=H_{rig}^m(M/K_0)$ (non logarithmic rigid cohomology) and\begin{gather} H_{rig}^m(M/{\mathfrak S}^a)=\underset{m=m_1+m_2}\bigoplus H_{dR}^{m_1}({\mathfrak M}_{{}\mathbb{Q}})\otimes_{K_0}(H_{{J_0}}^{m_2}\otimes_{K_0}K_0(a)).\tag{$*_a$}\end{gather}for $a\in\{0,\pi\}$. The decomposition $(*_0)$ is a decomposition into $\phi_{m,0}$-stable subspaces, and on the $(m_1,m_2)$-summand $\phi_{m,0}$ acts as $\phi_{m_1,triv}\otimes q^{m_2}$. 
\end{lem}

{\sc Proof:} Let ${\bf D}^0$ be the open unit disk, viewed as a $K_0$-dagger space. Let $\{x_j\}_{j\in J}$ be standard coordinates on the $|J|$-dimensional polydisk $({\bf D}^0)^{|J|}$. Let $t$ be a standard coordinate on ${\bf D}^0$ and define $\lambda:({\bf D}^0)^{|J|}\to {\bf D}^0$ by sending $t\mapsto t=\prod_{j\in {J}}x_{j}$. Endow $E=({\bf D}^0)^{|J|}$ and ${\bf D}^0$ with the log structures defined by the respective normal crossings divisor ${\bf V}(t)$. Denote by $E^{a}={\bf V}(t-a)\to\spm(K_0(a))$ the morphism obtained by the base change of $\lambda$ with $t\mapsto a$. Let $\Omega^{\bullet}_{E^{a}}$ be the corresponding relative logarithmic de Rham complex (of course, $\Omega^{\bullet}_{E^{\pi}}$ is the {\it usual} de Rham complex, without any additional log poles). Note that $\Omega^{1}_{E^{a}}$ is a free ${\cal O}_{E^{a}}$-module, one basis is $\{\dlog(x_{j}); j\in {J_0}\}$. Since ${J}=\Upsilon$, the $t$ defined here can be identified with the $t$ defined earlier and we get canonical identifications$${\nabla}^{a}_{{}}=]M[_{{\mathfrak V}^a_{}}={\mathfrak M}_{{}\mathbb{Q}}\times E^{a}$$(fibre products of dagger spaces over $\spm(K_0)$) and isomorphisms of complexes
$$(\Omega^{\bullet}_{\mathfrak{P},\mathbb{Q}}|_{]M[_{\mathfrak P}})\otimes{\mathcal O}_{{\nabla}^{a}}\cong p^*\Omega^{\bullet}_{{\mathfrak M}^{\sharp}_{{}\mathbb{Q}}}\otimes p^*\Omega^{\bullet}_{E^a}.$$ Explicit computations show$$H^s(E^a,\Omega^{\bullet}_{E^a})=H_{{J_0}}^s\otimes_{K_0}K_0(a).$$Thus $(*_{\pi})$ is the K\"unneth formula for $H^m_{dR}({\nabla}^{\pi}_{{}})=H^m({\nabla}^{\pi}_{{}},(\Omega^{\bullet}_{\mathfrak{P},\mathbb{Q}}|_{]M[_{\mathfrak P}})\otimes{\mathcal O}_{{\nabla}^{\pi}})$ proven in \cite{findag} section 2 Lemma 3; its proof can be literally repeated to prove $(*_0)$. The statement on $\phi_{m,0}$ holds since $\phi$ acts on each $\dlog(x_{j})$ by multiplication with $q$.\\

{\sc Proof of \ref{hkrigis}:} The isomorphy claim is a local statement, so we may assume the setting of \arabic{locset1}.\arabic{locset2}. Then by \ref{herzres}, restriction from $M^{\heartsuit}$ to $M$ does not change the cohomology objects in question. Therefore we may also assume $M\cap \cup_{j\in\widehat{J}}Y_j=\emptyset$; but then we may even assume ${J}=\Upsilon$. Choosing an auxiliary $j'\in{J'}$ the theorem follows from the explicit computations in \ref{nochiso} and \ref{restiso}.

\section{Weight filtration and monodromy on $H^*_{rig}(./{\mathfrak S}^{0})$}

\label{fromo} Here we explain some statements from \arabic{coueco1}.\arabic{coueco2} and \arabic{phien1}.\arabic{phien2}.

\addtocounter{satz}{1}{\bf \arabic{section}.\arabic{satz}}\newcounter{adli1}\newcounter{adli2}\setcounter{adli1}{\value{section}}\setcounter{adli2}{\value{satz}} An {\it admissible weak formal lift} of the semistable $k$-log scheme $(Y,{\cal N}_Y)$ is a weak formal ${\mathfrak S}$-log scheme $({\mathfrak Z},{\cal N}_{\mathfrak Z})$ together with an isomorphism of ${S}^{0} $-log schemes$$(Y,{\cal N}_Y)\cong({\mathfrak Z},{\cal N}_{\mathfrak Z})\times_{\mathfrak S}{S}^{0} $$satisfying the following conditions: On underlying weak formal schemes ${\mathfrak Z}$ is smooth over ${\mathfrak W}$, flat over ${\mathfrak S}$ and its reduction modulo $m_{A_0}$ is generically smooth over $S$; the fibre ${\mathfrak Y}={\bf V}(t)$ above $t=0$ is a divisor with normal crossings on ${\mathfrak Z}$, and ${\cal N}_{\mathfrak Z}$ is the log structure defined by this divisor. We denote an admissible weak formal lift by $({\mathfrak Z},{\mathfrak Y})$. Locally on $Y$, admissible lifts exist. Indeed, by \cite{fkato} 11.3 we locally find embeddings of $Y$ as a normal crossings divisor into smooth $k$-schemes $Z$. Assuming $Z$ is affine we can lift $Z$ to a smooth affine $A_0$-scheme $\widetilde{Z}$, see \cite{elk}, and we let ${\mathfrak Z}$ be its weak completion. Then we lift equations of $Y$ in ${\cal O}_{Z}$ (which form part of a local system of coordinates on $Z$) to equations in ${\cal O}_{\mathfrak Z}$: these define ${\mathfrak Y}$.

\addtocounter{satz}{1}{\bf \arabic{section}.\arabic{satz}} Choose an open covering $Y=\cup_{h\in H}U_h$ of $Y$, together with admissible liftings $({\mathfrak Z}_h,{\mathfrak Y}_h)$ of the $U_h$ (so $U_h$ is the reduction of ${\mathfrak Y}_h$). For a subset $G\subset H$ let $U_G=\cap_{h\in G}U_h$. Let $\{U_{G,\beta}\}_{\beta\in \Upsilon_G}$ be the set of irreducible components of $U_G$. For $h\in G$ and $\beta\in \Upsilon_G$ let ${\mathfrak Y}_{h,\beta}$ be the unique $A_0$-flat irreducible component of ${\mathfrak Y}_h$ with $U_{G,\beta}={\mathfrak Y}_{h,\beta}\cap U_G$. Let ${\mathfrak K}'_G$ be the blowing up of $\times_{\mathfrak W}({\mathfrak Z}_h)_{h\in G}$ along $\sum_{\beta\in \Upsilon_G}(\times_{\mathfrak W}({\mathfrak Y}_{h,\beta})_{h\in G})$, let ${\mathfrak K}_G$ be the complement of the strict transforms in ${\mathfrak K}'_G$ of all ${\mathfrak Y}_{h_0,\beta}\times (\times({\mathfrak Z}_h)_{h\in G-\{h_0\}})$ (i.e. all $h_0\in G$, all $\beta\in \Upsilon_G$), and let ${\mathfrak Y}_G$ be the exceptional divisor in ${\mathfrak K}_G$. It is a normal crossings divisor, and its   $A_0$-flat irreducible components are indexed by $\Upsilon_G$: they are the inverse images of the $A_0$-flat irreducible components of ${\mathfrak Y}_h$, for any $h\in G$. By construction, the diagonal embedding $U_G\to\times_{{\mathfrak W}}({\mathfrak Y}_h)_{h\in G}$ lifts canonically to an embedding $U_G\to{\mathfrak Y}_G\to {\mathfrak K}_G$. Viewing ${\mathfrak K}_G$ as a weak formal ${\mathfrak S}$-log scheme (with log structure defined by ${\mathfrak Y}_G$), this is an exact closed embedding of (weak formal) ${\mathfrak S}$-log schemes. 
Denote by $\tilde{\omega}_{{\mathfrak K}_G}^{\bullet}$ the logarithmic de Rham complex of ${\mathfrak K}_G\to {\mathfrak W}$. Write $\theta=\dlog(t)$ and let $$\tilde{\omega}_{{\mathfrak Y}_G}^{\bullet}=\tilde{\omega}_{{\mathfrak K}_G}^{\bullet}\otimes{\cal O}_{{\mathfrak Y}_G}\quad\quad\quad\quad\quad\omega_{{\mathfrak Y}_G}^{\bullet}=\frac{\tilde{\omega}_{{\mathfrak Y}_G}^{\bullet}}{\tilde{\omega}_{{\mathfrak Y}_G}^{\bullet-1}\wedge\theta}.$$So $\omega_{{\mathfrak Y}_G}^{\bullet}$ is the logarithmic de Rham complex of the morphism of weak formal log schemes ${\mathfrak Y}_G\to {\mathfrak S}^{0} $. Let ${\mathfrak Y}_{G,\mathbb{Q}}$ be the generic fibre of ${\mathfrak Y}_G$, a $K_0$-dagger space, and let $\tilde{\omega}_{{\mathfrak Y}_{G,\mathbb{Q}}}^{\bullet}$ resp. ${\omega}_{{\mathfrak Y}_{G,\mathbb{Q}}}^{\bullet}$ denote the sheaf complexes on ${\mathfrak Y}_{G,\mathbb{Q}}$ obtained from $\tilde{\omega}_{{\mathfrak Y}_G}^{\bullet}$ resp. ${\omega}_{{\mathfrak Y}_G}^{\bullet}$ by tensoring with $\mathbb{Q}$. Let $E\in{\cal LS}(Y,K_0)$. On the admissible open subspace $]U_G[_{{\mathfrak Y}_G}$ of ${\mathfrak Y}_{G,\mathbb{Q}}$ we define the sheaf complexes $$E\otimes_{K_0}\tilde{\omega}_{{\mathfrak Y}_{G,\mathbb{Q}}}^{\bullet}=sp^{-1}E\otimes_{K_0}\tilde{\omega}_{{\mathfrak Y}_{G,\mathbb{Q}}}^{\bullet}|_{]U_G[_{{\mathfrak Y}_G}}$$$$E\otimes_{K_0}{\omega}_{{\mathfrak Y}_{G,\mathbb{Q}}}^{\bullet}=sp^{-1}E\otimes_{K_0}{\omega}_{{\mathfrak Y}_{G,\mathbb{Q}}}^{\bullet}|_{]U_G[_{{\mathfrak Y}_G}}$$where $sp:]U_G[_{{\mathfrak Y}_G}\to U_G\subset Y$ is the specialization map. For $G_1\subset G_2$ we have natural transition maps $]U_{G_2}[_{{\mathfrak Y}_{G_2}}\to]U_{G_1}[_{{\mathfrak Y}_{G_1}}$. Hence a site $(]U_{G}[_{{\mathfrak Y}_{G}})_{G\subset H}=]U_{\bullet}[_{{\mathfrak Y}_{\bullet}}$ with sheaf complexes $E\otimes_{K_0}\tilde{\omega}_{{\mathfrak Y}_{\bullet}}^{\bullet}$ and $E\otimes_{K_0}{\omega}_{{\mathfrak Y}_{\bullet}}^{\bullet}$ on it. Clearly\begin{gather}{\mathbb{R}}\Gamma_{rig}(Y/{\mathfrak S}^{0},E)={\mathbb{R}}\Gamma(]U_{\bullet}[_{{\mathfrak Y}_{\bullet}},E\otimes_{K_0}{\omega}_{{\mathfrak Y}_{\bullet}}^{\bullet}).\tag{$1$}\end{gather}
Now we draw on a construction of Steenbrink \cite{steen}. For $j\ge 0$ let $$P_j\tilde{\omega}_{{\mathfrak K}_{\bullet}}^{k}=\bi(\tilde{\omega}_{{\mathfrak K}_{\bullet}}^{j}\otimes\Omega^{k-j}_{{\mathfrak K}_{\bullet}}\longrightarrow \tilde{\omega}_{{\mathfrak K}_{\bullet}}^{k})$$where $\Omega^{\bullet}_{{\mathfrak K}_{\bullet}}$ denotes the non-logarithmic de Rham complex on the simplicial weak formal scheme ${\mathfrak K}_{\bullet}$. Then let $$P_j\tilde{\omega}_{{\mathfrak Y}_{\bullet}}^{\bullet}=\frac{P_j\tilde{\omega}_{{\mathfrak K}_{\bullet}}^{\bullet}}{\tilde{\omega}_{{\mathfrak K}_{\bullet}}^{\bullet}\otimes{\mathfrak I}_{{\mathfrak Y}_{\bullet}}}$$where ${\mathfrak I}_{{\mathfrak Y}_{\bullet}}$ is the ideal of ${{\mathfrak Y}_{\bullet}}$ in ${\mathfrak K}_{\bullet}$. On ${\mathfrak Y}_{\bullet}$ these complexes give rise to a filtration $P_{\bullet}\tilde{\omega}_{{\mathfrak Y}_{\bullet}}^{\bullet}$ of $\tilde{\omega}_{{\mathfrak Y}_{\bullet}}^{\bullet}$. The graded pieces are computed by means of residue maps:\begin{gather}\Gr_j(\tilde{\omega}_{{\mathfrak Y}_{G,\mathbb{Q}}}^{\bullet})=\underset{{{\mathfrak N}\in\Theta_{j,G}}}\bigoplus\Omega^{\bullet}_{{\mathfrak N}_{\mathbb{Q}}}[-j].\tag{$2$}\end{gather}Here the sum runs through the set $\Theta_{j,G}$ of all intersections ${\mathfrak N}$ of $j$ different $A_0$-flat irreducible components of ${\mathfrak Y}_{G}$, and $\Omega^{\bullet}_{{\mathfrak N}_{\mathbb{Q}}}$ denotes the non-logarithmic de Rham complex on the smooth $K_0$-dagger space ${\mathfrak N}_{\mathbb{Q}}$, the generic fibre of ${\mathfrak N}$. The set of $A_0$-flat irreducible components of ${\mathfrak Y}_{G}$ corresponds bijectively to the set of irreducible components of the reduction of ${\mathfrak Y}_{G}$, and then further (intersect with the diagonally embedded $U_G$) to the set of irreducible components of $U_G$. Using the Poincar\'{e} lemma we find\begin{gather}{\mathbb{R}}\Gamma(]U_{\bullet}[_{{\mathfrak Y}_{\bullet}},E\otimes_{K_0}\Gr_j(\tilde{\omega}_{{\mathfrak Y}_{\bullet}}^{\bullet}))\cong\underset{N\in \Theta_j}\bigoplus E(N)\otimes_{K_0}{\mathbb{R}}\Gamma_{rig}(N/K_0)[-j]\tag{$3$}\end{gather}Here the sum runs through the set $\Theta_j$ of all intersections $N$ of $j$ different irreducible components of $Y$, on the right hand side we mean rigid cohomology with respect to trivial log structures as in \arabic{rigvebe1}.\arabic{rigvebe2}, and $E(N_j)$ means the value of the restriction of $E$ to $N_j$ (where it is constant). On ${\mathfrak Y}_{G,\mathbb{Q}}$ define the double complex $A_G^{\bullet\bullet}$ as follows: let$$A_G^{ij}=\frac{\tilde{\omega}_{{\mathfrak Y}_{G,\mathbb{Q}}}^{i+j+1}}{P_j(\tilde{\omega}_{{\mathfrak Y}_{G,\mathbb{Q}}}^{i+j+1})},$$as differentials $A_G^{ij}\to A_G^{(i+1)j}$ take those induced by $(-1)^jd$, and as differentials $A_G^{ij}\to A_G^{i(j+1)}$ take those induced by $\omega\mapsto\omega\wedge\theta$. Let $A_G^{\bullet}$ be the associated total complex. We claim that the augmentation $\tilde{\omega}_{{\mathfrak Y}_{G,\mathbb{Q}}}^{\bullet}\to A_G^{\bullet 0}$ defined by $\omega\mapsto\omega\wedge\theta$ induces a quasi-isomorphism ${\omega}_{{\mathfrak Y}_{G,\mathbb{Q}}}^{\bullet}\to A_G^{\bullet}$. As in \cite{steen} it suffices to prove that$$0\longrightarrow\Gr_0(\tilde{\omega}_{{\mathfrak Y}_{G,\mathbb{Q}}}^{\bullet})\stackrel{\wedge\theta}{\longrightarrow}\Gr_1(\tilde{\omega}_{{\mathfrak Y}_{G,\mathbb{Q}}}^{\bullet})[1]\stackrel{\wedge\theta}{\longrightarrow}\Gr_2(\tilde{\omega}_{{\mathfrak Y}_{G,\mathbb{Q}}}^{\bullet})[2]\stackrel{\wedge\theta}{\longrightarrow}\ldots$$is exact. In view of $(2)$ this means proving that$$0\longrightarrow \Omega^{\bullet}_{{\mathfrak Y}_{G,\mathbb{Q}}}\longrightarrow\underset{{\mathfrak N}\in\Theta_{1,G}}\bigoplus\Omega^{\bullet}_{{\mathfrak N}_{\mathbb{Q}}}\longrightarrow \underset{{\mathfrak N}\in\Theta_{1,G}}\bigoplus\Omega^{\bullet}_{ {\mathfrak N}_{\mathbb{Q}}}\longrightarrow\ldots$$is exact; but this is a general fact on differential modules on normal crossings intersections of smooth spaces. Combining for varying $G$ we get a quasi-isomorphism ${\omega}_{{\mathfrak Y}_{\bullet}}^{\bullet}\to A_{\bullet}^{\bullet}$ of sheaf complexes on ${\mathfrak Y}_{\bullet}$. Define the weight filtration on $A_{\bullet}^{\bullet}$ by setting$$P_kA^{ij}_G=\frac{ P_{2j+k+1}(\tilde{\omega}_{{\mathfrak Y}_{G,\mathbb{Q}}}^{i+j+1})}{P_j(\tilde{\omega}_{{\mathfrak Y}_{G,\mathbb{Q}}}^{i+j+1})}.$$The associated spectral sequence for the cohomology of $A_{\bullet}^{\bullet}$ then reads, using the (quasi)isomorphisms $(1)$ and ${\omega}_{{\mathfrak Y}_{\bullet}}^{\bullet}\cong A_{\bullet}^{\bullet}$, and the computation $(3)$:\begin{gather}E_1^{-k,i+k}=\underset{{j\ge0\atop j\ge -k}}\bigoplus\underset{N\in\Theta_{2j+k+1}}\bigoplus E(N)\otimes_{K_0}H_{rig}^{i-2j-k}(N/K_0)\Longrightarrow H_{rig}^i(Y/{\mathfrak S}^{0},E)\tag{$4$}\end{gather}

\begin{satz}\label{maint} (i) If $Y$ is quasi-compact, $H^*_{rig}(Y/{\mathfrak S}^{0},E)$ is finite dimensional.\\(ii) If all irreducible components of $Y$ are proper, the canonical morphism ${\mathbb{R}}\Gamma_{rig}(Y/{\mathfrak S}^{0},E)\to {\mathbb{R}}\Gamma_{conv}(Y/{\mathfrak S}^{0},E)$ is an isomorphism.\\(iii) If $k$ is perfect, $A_0=W(k)$ and $\sigma$ are as in \arabic{fro1}.\arabic{fro2}, and $(E,\phi)\in F\mbox{-}{\cal LS}(Y,K_0)$, then the endomorphism $\phi$ on $H^*_{rig}(Y/{\mathfrak S}^{0},E)$ is bijective.\\(iv) In (iii), if $k$ is finite and if for any intersection $N$ of irreducible components of $Y$ the eigenvalues of $\phi$ acting on $E(N)$ are Weil numbers, then $H^*_{rig}(Y/{\mathfrak S}^{0},E)$ is a mixed $F$-isocrystal: the eigenvalues of $\phi$ on $H^*_{rig}(Y/{\mathfrak S}^{0},E)$ are Weil numbers.  
\end{satz} 

{\sc Proof:} Assertions (i), (iii) and (iv) follow easily from the spectral sequence $(4)$ (which in cases (iii) and (iv) is Frobenius equivariant) and the corresponding results for the rigid cohomology with constant coefficients of (classically smooth) $k$-schemes, see \cite{berfi} \cite{chiwe}. For (ii) observe that we can repeat all constructions using rigid spaces instead of dagger spaces, obtaining the spectral sequence\begin{gather}E_1^{-k,i+k}=\underset{j\ge0\atop j\ge -k}\bigoplus\underset{N\in\Theta_{2j+k+1}}\bigoplus E(N)\otimes_{K_0}H_{conv}^{i-2j-k}(N/K_0)\Longrightarrow H_{conv}^i(Y/{\mathfrak S}^{0},E).\tag{$4'$}\end{gather}This reduces the problem to proving that the maps $H_{rig}^{*}(N/K_0)\to H_{conv}^{*}(N/K_0)$ are isomorphisms. Since the $N$ are proper, this is done (in view of \arabic{rigvebe1}.\arabic{rigvebe2}) in \cite{berfi}.

\addtocounter{satz}{1}{\bf \arabic{section}.\arabic{satz}}\newcounter{mode1}\newcounter{mode2}\setcounter{mode1}{\value{section}}\setcounter{mode2}{\value{satz}} ($k$ arbitrary) By construction, we have a short exact sequence\begin{gather}0\longrightarrow E\otimes_{K_0}{\omega}_{{\mathfrak Y}_{\bullet}}^{\bullet}[-1]\stackrel{\wedge\theta}{\longrightarrow} E\otimes_{K_0}\tilde{\omega}_{{\mathfrak Y}_{\bullet}}^{\bullet}\longrightarrow E\otimes_{K_0}{\omega}_{{\mathfrak Y}_{\bullet}}^{\bullet}\longrightarrow 0.\tag{$5$}\end{gather}By definition, the {\it monodromy operator}$$N:H^*_{rig}(Y/{\mathfrak S}^{0},E)\longrightarrow H^*_{rig}(Y/{\mathfrak S}^{0},E)$$is the connecting homomorphism in cohomology associated with $(5)$.

\begin{pro}\label{nfpfn} Suppose $q, \sigma$ are as in \arabic{fro1}.\arabic{fro2}, and $(E,\phi)\in F\mbox{-}{\cal LS}(Y,K_0)$. Then we have $N\phi=q\phi N$ on $H^*_{rig}(Y/{\mathfrak S}^{0},E )$.  
\end{pro} 

{\sc Proof:} Let $\sigma$ be the unique endomorphism of the weak formal log scheme ${\mathfrak S}$ which equals $\sigma$ on scalars, which sends $t$ to $t^{q}$, and for which multiplication with $q$ on the standard chart $\mathbb{N}$ of ${\mathfrak S}$ is a chart. Specializing $t\mapsto 0$ we get the endomorphism $\sigma:{\mathfrak S}^{0} \to {\mathfrak S}^{0} $ defined in \arabic{fro1}.\arabic{fro2}. We may assume that for each $h\in H$ we are given an endomorphism of the weak formal ${\mathfrak S}$-log scheme ${\mathfrak Z}_h$ which lifts the $q$-th power Frobenius endomorphism of the reduction of ${\mathfrak Z}_h$, which sends equations for the divisors ${\mathfrak Y}_{h,\beta}$ on ${\mathfrak Z}_h$ to their $q$-th powers, and which is compatible with $\sigma$ on ${\mathfrak S}$. Then we also get an endomorphism of the site $({\mathfrak K}_G)_{G\subset H}$. It gives rise to endomorphisms $\phi$ on $\tilde{\omega}_{{\mathfrak Y}_{\bullet}}^{\bullet}$ and on ${\omega}_{{\mathfrak Y}_{\bullet}}^{\bullet}$. Tensoring with $\phi:E\to E$ we get endomorphisms $\phi$ on $E\otimes_{K_0}\tilde{\omega}_{{\mathfrak Y}_{\bullet}}^{\bullet}$ and on $E\otimes_{K_0}{\omega}_{{\mathfrak Y}_{\bullet}}^{\bullet}$ (inducing $\phi$ on $H^*_{rig}(Y/{\mathfrak S}^{0},E)$ as defined in \arabic{fro1}.\arabic{fro2}). The claim now follows from the following fact: If in the exact sequence $(5)$ we act on the first (non-zero) term by $q\phi$, on the second and the third term by $\phi$, then these actions are compatible with the maps in $(5)$ (note $\phi(\theta)=q\theta$).\\ 
\begin{pro}\label{corit} Let $M$ be the intersection of some irreducible components of $Y$. If $M$ is proper, the natural morphism$${\mathbb{R}}\Gamma_{rig}(M/{\mathfrak S}^{0} )\longrightarrow {\mathbb{R}}\Gamma_{conv}(M/{\mathfrak S}^{0} )$$is an isomorphism. Moreover $\phi$ on $H^*_{rig}(M/{\mathfrak S}^{0} )$ is bijective; if $k$ is finite, then $H^*_{rig}(M/{\mathfrak S}^{0} )$ is mixed.
\end{pro}

{\sc Proof:} This is a slight modification of the proof of \ref{maint}. We keep notations from above. For $G\subset H$ let $M_G=M\cap U_G$, let ${\mathfrak M}_G$ be the intersection of the $A_0$-flat irreducible components of ${\mathfrak Y}_G$ containing $M_G$, and let ${\mathfrak M}_{G,\mathbb{Q}}$ be its generic fibre (as a dagger space). Varying $G$ we get a site as before. The natural map\begin{gather}{\mathbb{R}}\Gamma_{rig}(M/{\mathfrak S}^{0} )={\mathbb{R}}\Gamma(]M_{\bullet}[_{{\mathfrak Y}_{\bullet}},{\omega}_{{\mathfrak Y}_{\bullet}}^{\bullet})\longrightarrow {\mathbb{R}}\Gamma(]M_{\bullet}[_{{\mathfrak Y}_{\bullet}},{\omega}_{{\mathfrak Y}_{\bullet}}^{\bullet}\otimes{\cal O}_{{\mathfrak M}_{\bullet}})\tag{$1_M$}\end{gather}is an isomorphism. Indeed, one can check this for every $G$ separately; applying the Poincar\'{e} lemma to both sides one reduces to the case where $G$ consists of a single element; but then the claim is proven in \ref{anarighomot} (which is formulated only for the embedding of $M$ into a particular admissible weak formal lift of $V_M$; but the same proof applies here). For $G\subset H$ consider the filtration of $\tilde{\omega}_{{\mathfrak Y}_{G,\mathbb{Q}}}^{\bullet}\otimes{\cal O}_{{\mathfrak M}_{G,\mathbb{Q}}}$ defined by $$P_j(\tilde{\omega}_{{\mathfrak Y}_{G,\mathbb{Q}}}^{\bullet}\otimes{\cal O}_{{\mathfrak M}_{G,\mathbb{Q}}})=\frac{P_j\tilde{\omega}_{{\mathfrak Y}_{G,\mathbb{Q}}}^{\bullet}}{{\mathfrak I}_{{\mathfrak M}_{G,\mathbb{Q}}}\otimes \tilde{\omega}_{{\mathfrak Y}_{G,\mathbb{Q}}}^{\bullet}\cap P_j\tilde{\omega}_{{\mathfrak Y}_{G,\mathbb{Q}}}^{\bullet}}$$where ${\mathfrak I}_{{\mathfrak M}_{G,\mathbb{Q}}}\subset {\cal O}_{{\mathfrak Y}_{G,\mathbb{Q}}}$ is the ideal of ${\mathfrak M}_{G,\mathbb{Q}}$ in ${\mathfrak Y}_{G,\mathbb{Q}}$. From $(2)$ one derives\begin{gather}\Gr_j(\tilde{\omega}_{{\mathfrak Y}_{G,\mathbb{Q}}}^{\bullet}\otimes{\cal O}_{{\mathfrak M}_{G,\mathbb{Q}}})=\underset{{\mathfrak N}\in\Theta_{j,G}}\bigoplus\Omega^{\bullet}_{{\mathfrak M}_{G,\mathbb{Q}}\cap {\mathfrak N}_{\mathbb{Q}}}[-j],\tag{$2_M$}\end{gather} with $\Omega^{\bullet}_{{\mathfrak M}_{G,\mathbb{Q}}\cap {\mathfrak N}_{\mathbb{Q}}}$ the non-logarithmic de Rham complex on the smooth $K_0$-dagger space ${\mathfrak M}_{G,\mathbb{Q}}\cap {\mathfrak N}_{\mathbb{Q}}$. Using the Poincar\'{e} lemma we this time find\begin{gather}{\mathbb{R}}\Gamma(]M_{\bullet}[_{{\mathfrak Y}_{\bullet}},\Gr_j(\tilde{\omega}_{{\mathfrak Y}_{\bullet}}^{\bullet}\otimes{\cal O}_{{\mathfrak M}_{\bullet}}))\cong\underset{N\in\Theta_j}\bigoplus {\mathbb{R}}\Gamma_{rig}(M\cap N/K_0)[-j].\tag{$3_M$}\end{gather} For $G\subset H$ define on ${\mathfrak Y}_{G,\mathbb{Q}}$ the double complex $A_{M,G}^{\bullet\bullet}$ as follows: let$$A_{M,G}^{ij}=\frac{\tilde{\omega}_{{\mathfrak Y}_{G,\mathbb{Q}}}^{i+j+1}\otimes{\cal O}_{{\mathfrak M}_{G,\mathbb{Q}}}}{P_j(\tilde{\omega}_{{\mathfrak Y}_{G,\mathbb{Q}}}^{i+j+1}\otimes{\cal O}_{{\mathfrak M}_{G,\mathbb{Q}}})},$$as differentials $A_{M,G}^{ij}\to A_{M,G}^{(i+1)j}$ take those induced by $(-1)^jd$, and as differentials $A_{M,G}^{ij}\to A_{M,G}^{i(j+1)}$ take those induced by $\omega\mapsto\omega\wedge\theta$. If $A_{M,G}^{\bullet}$ is the associated total complex, the augmentation $\tilde{\omega}_{{\mathfrak Y}_{G,\mathbb{Q}}}^{\bullet}\otimes{\cal O}_{{\mathfrak M}_{G,\mathbb{Q}}}\to A_{M,G}^{\bullet 0}$ defined by $\omega\mapsto\omega\wedge\theta$ induces a quasi-isomorphism ${\omega}_{{\mathfrak Y}_{G,\mathbb{Q}}}^{\bullet}\otimes{\cal O}_{{\mathfrak M}_{G,\mathbb{Q}}}\cong A_{M,G}^{\bullet}$. This time this assertion is reduced, using $(2_M)$, to proving that $$0\longrightarrow \Omega^{\bullet}_{{\mathfrak M}_{G,\mathbb{Q}}}\longrightarrow\underset{{\mathfrak N}\in \Theta_{1,G}}\bigoplus\Omega^{\bullet}_{{\mathfrak M}_{G,\mathbb{Q}}\cap {\mathfrak N}_{\mathbb{Q}}}\longrightarrow \underset{{\mathfrak N}\in \Theta_{1,G}}\bigoplus\Omega^{\bullet}_{{\mathfrak M}_{G,\mathbb{Q}}\cap {\mathfrak N}_{\mathbb{Q}}}\longrightarrow\ldots$$is exact; see \cite{mono} for an elementary proof in a similar situation. Hence, filtering the complex $A_{M,\bullet}^{\bullet}$ on ${\mathfrak Y}_{\bullet}$ analogously as before and applying ${\mathbb{R}}\Gamma(]M_{\bullet}[_{{\mathfrak Y}_{\bullet}},.)$ we get\begin{gather}E_1^{-k,i+k}=\underset{j\ge0\atop j\ge -k}\bigoplus\underset{N\in \Theta_{2j+k+1}}\bigoplus H_{rig}^{i-2j-k}(M\cap N/K_0)\Longrightarrow H_{rig}^i(M/{\mathfrak S}^{0} ).\tag{$4_M$}\end{gather}Now all this can be literally repeated using rigid spaces instead of dagger spaces. The only point where the argument must be varied is the proof of the rigid space version of \ref{anarighomot} (which we need to get the rigid space version of $(1_M)$). In fact, the rigid space version of \ref{anarighomot} is even easier: the reason why in \ref{anarighomot} we worked with the affinoid covering ${\nabla}^{0}_{q}=\cup_{\rho<1}{\nabla}^{0}_{q,\rho}$ is that we do not know the acyclicity of ${\nabla}^{0}_{q}$ for coherent ${\cal O}_{{\nabla}^{0}_{q}}$-modules. But for coherent modules over the associated rigid structure sheaf we know it (by Kiehl's theorem), i.e. in the rigid space context we simply do without the passage to the affinoid covering (which, of course, even would not work in the rigid space context). So we get a spectral sequence\begin{gather}E_1^{-k,i+k}=\underset{j\ge0\atop j\ge -k}\bigoplus\underset{N\in \Theta_{2j+k+1}}\bigoplus H_{conv}^{i-2j-k}(M\cap N/K_0)\Longrightarrow H_{conv}^i(M/{\mathfrak S}^{0} )\tag{$4'_M$}\end{gather}and we conclude as in \ref{maint}. We are done.

\section{The symmetric space and its weak formal model} We specialize to the case where $K$ is a finite extension of $\mathbb{Q}_p$ and $A_0=W(k)$. Let $q=|k|$ and for this $q$ let $\sigma$ be as in \arabic{fro1}.\arabic{fro2}. We still fix a uniformizer $\pi\in A={\cal O}_K$. 

\addtocounter{satz}{1}{\bf \arabic{section}.\arabic{satz}} Let $d\ge 2$. We describe a weak formal $A$-scheme ${\mathfrak Q}$ whose generic fibre is Drinfel'd's $p$-adic symmetric space $\Omega^{(d+1)}_K$ of dimension $d$ over $K$. (-- The description given in \cite{gen} which we work out here is not sufficiently detailed and (strictly speaking) not entirely correct, it seems to us. --) Let ${\mathfrak T}$ be the Bruhat-Tits building of ${\rm PGL}_{d+1}(K)$. We fix a copy ${\bf P}^d_K$ of projective $d$-space over $K$. The set of vertices of ${\mathfrak T}$ is in natural bijection with the set of equivalence classes of pairs $(\widetilde{P},\phi)$, where $\widetilde{P}$ is an $A$-scheme isomorphic to ${\bf P}^d_{A}$, and $\phi$ is an isomorphism of $K$-schemes $\widetilde{P}\otimes_{A}K\cong{\bf P}^d_K$. Two pairs $(\widetilde{P}_1,\phi_1)$ and $(\widetilde{P}_2,\phi_2)$ are equivalent if there exists an isomorphism $\widetilde{P}_1\cong \widetilde{P}_2$ respecting the $\phi_i$. For a vertex $v$ in ${\mathfrak T}$ let $(\widetilde{P}_v,\phi_v)$ be the corresponding pair. Each $k$-rational linear subscheme of the reduction $P_{v}$ of $\widetilde{P}_v$ defines another pair $(\widetilde{P}',\phi')$; namely, $\widetilde{P}'$ is the projective space over $A$ spanned by the kernel of the map of projective coordinate rings to which the inclusion gives rise. The vertices in ${\mathfrak T}$ adjacent to $v$ are precisely those corresponding to the pairs $(\widetilde{P}',\phi')$ which can be obtained in this way. This remark is used below.\\
We fix a vertex $v_0$ of ${\mathfrak T}$. For $n\ge0$ let $V_n$ be the set of vertices $v$ with $d(v_0,v)\le n$. Here we write $d(v_0,v)$ for the minimal number $d\ge0$ for which there is a sequence $v_{0},v_1,\ldots,v_{d}=v$ with $v_i$ adjacent to $v_{i-1}$ for all $i\ge1$. For $n\ge 1$ let $W_n=V_n-V_{n-1}$. For $m\ge1$ let $V_n^m$ be the union of $V_{n}$ with the set of vertices $v$ in $W_{n+1}$ such that there exists an element $w\in V_{n}$ adjacent to $v$ such that the subscheme of $\widetilde{P}_w$ corresponding to $v$ has dimension at most $m-1$. In particular $V_n^d=V_{n+1}$. For a vertex $w$ let $H_w$ be the set of all vertices adjacent to (but different from) $w$.\\
Now we can begin. Blow up $\widetilde{Q}^{0,0}=\widetilde{P}_{v_0}$ in all $k$-rational points of the reduction $P_{v_{0}}=P'_{v_{0}}$ of $\widetilde{P}_{v_0}=\widetilde{Q}^{0,0}$, call the result $\widetilde{Q}^{0,1}$. Then blow up $\widetilde{Q}^{0,1}$ in the strict transforms of all $k$-rational lines of $P'_{v_{0}}$, call the result $\widetilde{Q}^{0,2}$; ...; finally blow up $\widetilde{Q}^{0,d-1}$ in the strict transforms of all $(d-1)$-dimensional $k$-rational linear subschemes of $P'_{v_{0}}$, call the result $\widetilde{Q}^{0,d}=\widetilde{Q}^{1,0}$. For $w\in W_1$ let $P'_{w}$ be the succesive blowing up of the reduction $P_{w}$ of $\widetilde{P}_w$ along its subschemes which correspond to the vertices in $H_w\cap V_1$ (first blow up the points, then the strict transforms of the lines, etc.). Taking strict transforms under this blowing up sequence, the vertices of $H_w-(H_w\cap V_1)$ correspond to subschemes of $P'_{w}$. This $k$-scheme $P'_{w}$ can be identified with a subscheme of $\widetilde{Q}^{1,0}$. Now blow up $\widetilde{Q}^{1,0}$ in all subschemes of $P'_{w}$ corresponding to vertices in $H_w\cap (V_1^1-V_1)$, for all $w\in W_1$. Call the result $\widetilde{Q}^{1,1}$. Then blow up $\widetilde{Q}^{1,1}$ in the strict transforms of all subschemes of $P'_{w}$ corresponding to vertices in $H_w\cap (V_1^2-V_1^1)$, for all $w\in W_1$. Call the result $\widetilde{Q}^{1,2}$; ...; finally blow up $\widetilde{Q}^{1,{d-1}}$ in the strict transforms of all subschemes of $P'_{w}$ corresponding to vertices in $H_w\cap (V_1^d-V_1^{d-1})$, for all $w\in W_1$. Call the result $\widetilde{Q}^{1,d}=\widetilde{Q}^{2,0}$. Keep going. We get a sequence$$\ldots\longrightarrow \widetilde{Q}^{n,0}\longrightarrow \widetilde{Q}^{n-1,0}\longrightarrow\ldots\longrightarrow \widetilde{Q}^{1,0}\longrightarrow \widetilde{Q}^{0,0}.$$Remove from $\widetilde{Q}^{n,0}$ the centers of the sequence of blowing ups $\widetilde{Q}^{n+1,0}\to \widetilde{Q}^{n,0}$, call the result $\widetilde{Q}^{(n)}$. These $\widetilde{Q}^{(n)}$ form a sequence of open immersions$$\widetilde{Q}^{(0)}\longrightarrow \widetilde{Q}^{(1)}\longrightarrow\ldots\longrightarrow \widetilde{Q}^{(n)}\longrightarrow\ldots.$$Let $\widetilde{Q}^{(\infty)}$ be its inductive limit and let ${\mathfrak Q}$ be the weak completion of $\widetilde{Q}^{(\infty)}$. 
   
\addtocounter{satz}{1}{\bf \arabic{section}.\arabic{satz}}\newcounter{compde1}\newcounter{compde2}\setcounter{compde1}{\value{section}}\setcounter{compde2}{\value{satz}} From our description of ${\mathfrak Q}$ we deduce the following facts (cf. \cite{ito}, sect. 6). It is a strictly semistable weak formal $A$-scheme. Let $Q$ be the special fibre of ${\mathfrak Q}$. Any non-empty intersection of distinct irreducible components of $Q$ is isomorphic to the product of $k$-schemes each of which results from the following procedure (for some $r\ge0$): First blow up all $k$-rational points in projective $r$-space ${\bf P}^{r}_k$; then blow up the strict transforms (in this blow up) of all $k$-rational lines in ${\bf P}^{r}_k$; ...; finally blow up the strict transforms of all $k$-rational linear subschemes of ${\bf P}^{r}_k$ of dimension $r-2$.

\begin{satz}\label{compur} Let $M$ be the intersection of $i$ different irreducible components of $Q$ and let $s\ge0$. Then we have $\phi=q^{s}$ on $H_{rig}^s(M/{\mathfrak S}^0)$.
\end{satz}

{\sc Proof:} Let $M^{\heartsuit}$ be the open complement in $M$ of the intersection of $M$ with the union of all irreducible components of $Q$ not containing $M$. Then $H_{rig}^s(M/{\mathfrak S}^0)\cong H_{rig}^s(M^{\heartsuit}/{\mathfrak S}^0)$ by \ref{herzres}, so we need to show $\phi=q^{s}$ on $H_{rig}^s(M^{\heartsuit}/{\mathfrak S}^0)$. In view of \ref{restiso} we only need to show $\phi=q^{m}$ on $H_{rig}^m(M^{\heartsuit}/K_0)$ for all $m\ge 0$. Choosing a product decomposition of $M$ according to \arabic{compde1}.\arabic{compde2} we get a product decomposition of $M^{\heartsuit}$ into $k$-schemes each of which is the complement in ${\bf P}^{r}_k$ (for some $r$) of all $k$-rational linear hyperplanes. By the K\"unneth formula we may treat every factor separately. We claim: For a finite non-empty set $G$ of $k$-rational linear hyperplanes in ${\bf P}^{r}$, if we let $\cup_G=\cup_{H\in G}H$, then $\phi=q^m$ on $H^m_{rig}({\bf P}_k^{r}-\cup_G/K_0)$ for all $m\ge 0$, all $r\ge1$. We induce on the cardinality of $G$. For $|G|=1$ we have ${\bf P}_k^{r}-\cup_G=\mathbb{A}^r_k$, but $H_{rig}^m(\mathbb{A}^r_k/K_0)=0$ for $m>0$, and $\phi=\id$ on $H_{rig}^0({\mathbb{A}}^r_k/K_0)$. If $|G|>1$ pick $H'\in G$, let $G'=G-\{H'\}$ and consider the $\phi$-equivariant exact Gysin sequence$$H^m_{rig}({\bf P}^{r}-\cup_{G'}/K_0)\longrightarrow H^m_{rig}({\bf P}^{r}-\cup_G/K_0)\to H^{m-1}_{rig}(H'-(H'\cap\cup_{G'})/K_0)(-1)$$(see \cite{chiwe}). By induction hypothesis we have $\phi=q^m$ on both outer terms.

\addtocounter{satz}{1}{\bf \arabic{section}.\arabic{satz}} Let $\Omega^{(d+1)}_K=X$ be the generic fibre of ${\mathfrak Q}$, a $K$-dagger analytic space. Since it is a {\it Stein} space its coherent and de Rham cohomology is (by \cite{crelle}) the same as that of the associated rigid space, Drinfel'd's $p$-adic symmetric space of dimension $d$ over $K$. With notation $]Q^{r}[_{\mathfrak Q}$ as in \arabic{lifsit1}.\arabic{lifsit2} we have:
 
\begin{pro}\label{desh} (de Shalit) For any $s\in\mathbb{Z}$, the complex $$H^s_{dR}(]Q^{1}[_{\mathfrak Q})\longrightarrow H^s_{dR}(]Q^{2}[_{\mathfrak Q})\longrightarrow H^s_{dR}(]Q^{3}[_{\mathfrak Q})\longrightarrow\ldots$$is a resolution of $H^s_{dR}(X)$.
\end{pro}

{\sc Proof:} This follows from Theorems 5.7 and 7.7 in the paper \cite{ds} of Ehud de Shalit. Strictly speaking, de Shalit does not work with the (weak) formal scheme ${\mathfrak Q}$. Instead, he works with the Bruhat-Tits building ${\mathfrak T}$ of ${\rm PGL}_{d+1}(K)$, so we give some remarks on how his setting translates into ours. ${\mathfrak T}$ can be regarded as dual to $Q$: A vertex $v$ of ${\mathfrak T}$ corresponds to an irreducible component $C_{v}$ of $Q$. Let $C_{v}^{\heartsuit}$ be the maximal open subscheme of $C_{v}$ which has empty intersection with all other irreducible components of $Q$. There is a "reduction map" $r:X(\overline{K})\to |{\mathfrak T}|.$ The preimage of a vertex $v$ under $r$ is the set of $\overline{K}$-valued points of the preimage of $C_{v}^{\heartsuit}$ under the specialization map $sp$. The preimage of the star of a vertex $v$ under $r$ is the set of $\overline{K}$-valued points of the preimage of $C_{v}$ under $sp$. In a similar way, preimages of $k$-cells (for $k\ge0$) correspond to preimages of: intersections of $k+1$ distinct irreducible components minus their intersection with other components. Also note that since the spaces $]Q^{i}[_{\mathfrak Q}$ are partially proper their de Rham cohomology is the same for the dagger space version as for the rigid space version (by \cite{crelle}).

\begin{kor}\label{sympur} For each $s\ge0$ we have $\phi=q^{s}$ on $H_{rig}^s(Q/{\mathfrak S}^0)$ and there is a canonical isomorphism $H_{dR}^{s}(X)\cong H_{rig}^s(Q/{\mathfrak S}^0)\otimes_{K_0}K$. 
\end{kor}

{\sc Proof:} For the spectral sequences \begin{gather}E_{1}^{rs}=H^s_{rig}(Q^{r+1}/{\mathfrak S}^0)\Longrightarrow H^{r+s}_{rig}(Q/{\mathfrak S}^0)\tag*{$(U)_{0}$}\\E_{1}^{rs}=H^s_{dR}(]Q^{r+1}[_{\mathfrak Q})\Longrightarrow H^{r+s}_{dR}(X)\tag*{$(U)_{\pi}$}\end{gather}we have $(U)_{0}\otimes_{K_0}K\cong(U)_{\pi}$ by \ref{rigs}. In $(U)_{0}$ we know $\phi=q^{s}$ on $E_1^{0s}$ by \ref{compur}, so it is enough to show $E_2^{rs}=0$ whenever $r>0$. But this we can equivalently check in $(U)_{\pi}$ where it follows from \ref{desh}.

\section{Hodge decomposition on quotients of $\Omega^{(d+1)}_K=X$}

\addtocounter{satz}{1}{\bf \arabic{section}.\arabic{satz}} There is a natural action of ${\rm PGL}_{d+1}(K)$ on ${\mathfrak Q}$. We fix a cocompact discrete subgroup $\Gamma\subset {\rm PGL}_{d+1}(K)$. Passing to a subgroup with finite index in $\Gamma$ we may and will assume that $\Gamma$ is torsionfree and that the quotient ${\mathfrak Q}_{\Gamma}=\Gamma\backslash {\mathfrak Q}$ has strictly semistable reduction. By \cite{mus} it even algebraizes to a projective $A$-scheme, and the covering map ${\mathfrak Q}\to {\mathfrak Q}_{\Gamma}$ is \'{e}tale. Let $X_{\Gamma}=\Gamma\backslash X$ be the generic fibre, let $Q_{\Gamma}=\Gamma\backslash Q$ be the special fibre of ${\mathfrak Q}_{\Gamma}$. Furthermore we fix a finite dimensional $K[\Gamma]$-module ${\bf F}$. It gives rise to an element $F\in{\cal LS}(X_{\Gamma},K)$; namely, for admissible open $U\subset X_{\Gamma}$ let $F(U)={F_X}(X\times_{X_{\Gamma}}U)^{\Gamma}$ where $F_X$ is the constant sheaf on $X$ associated with ${\bf F}$. The main tool in \cite{schn} for studying $H_{dR}^{*}(X_{\Gamma},F)$ is the covering spectral sequence\begin{gather} E_2^{rs}=H^r(\Gamma,{\bf F}\otimes_K H_{dR}^s(X))\Longrightarrow H_{dR}^{r+s}(X_{\Gamma},F).\tag*{$(G)_{\pi}$}\end{gather} Only the cohomology $H_{dR}^{d}(X_{\Gamma},F)$ in middle degree $d$ is interesting (\cite{schn}). Let $$H_{dR}^d(X_{\Gamma},F)=F^0_{\Gamma}\supset F^1_{\Gamma}\supset\ldots\supset F^{d+1}_{\Gamma}=0$$ be the descending filtration which $(G)_{\pi}$ defines on it. The other analytically defined filtration on $H_{dR}^{d}(X_{\Gamma},F)$ is the canonical Cech filtration: the descending filtration $(F_{C}^r)_{r\ge0}$ on $H^{d}_{dR}(X_{\Gamma},F)$ induced by the spectral sequence\begin{gather}E_{1}^{rs}=F(]Q_{\Gamma}^{r+1}[_{{\mathfrak Q}_{\Gamma}})\otimes_{K}H^s_{dR}(]Q_{\Gamma}^{r+1}[_{{\mathfrak Q}_{\Gamma}})\Longrightarrow H^{s+r}_{dR}(X_{\Gamma},F).\tag*{$(C)_{\pi}$}\end{gather}
 
\addtocounter{satz}{1}{\bf \arabic{section}.\arabic{satz}} Letting $E=sp_* F\in{\cal LS}(Q_{\Gamma},K)$ we may view $(E,\id_E)$ as an object in $F\mbox{-}{\cal LS}(Q_{\Gamma},K_0)$. In particular we get $K_0$-linear endomorphisms $\phi$ and $N$ on $H^d_{rig}(Q_{\Gamma}/{\mathfrak S}^{0},E)$ satisfying $N\phi=q\phi N$. Via the isomorphism $$H^d_{dR}(X_{\Gamma},F)\cong H^d_{rig}(Q_{\Gamma}/{\mathfrak S}^{0},E)$$ from \ref{rigs} we view $\phi$ and $N$ as endomorphisms $\phi$ and $N$ on $H^d_{dR}(X_{\Gamma},F)$.

\begin{satz}\label{charfil} The spectral sequences $(G)_{\pi}$ and $(C)_{\pi}$ degenerate in $E_2$. The filtrations $(F_{\Gamma}^r)_{r\ge0}$ and $(F_{C}^r)_{r\ge0}$ on $H^{d}_{dR}(X_{\Gamma},F)$ coincide. They are stable for the action of $\phi$, and we have $\phi=q^{d-r}$ on $F_{\Gamma}^r/F^{r+1}_{\Gamma}=F_{C}^r/F^{r+1}_{C}$. 
\end{satz}

{\sc Proof:} From \ref{rigs} we get an isomorphism of spectral sequences between $(C)_{\pi}$ and \begin{gather}E_{1}^{rs}={\bf F}\otimes_{K}H^s_{rig}(Q_{\Gamma}^{r+1}/{\mathfrak S}^0)\Longrightarrow H^{r+s}_{rig}(Q_{\Gamma}/{\mathfrak S}^0,E)\tag*{$(C)_{0}$}.\end{gather}$\phi$ acts on $(C)_0$, and on its $E_1^{rs}$-term we have $\phi=q^{s}$, by \ref{compur}. On the other hand, the locally constant sheaf $E$ on $Q_{\Gamma}$ is obtained by taking $\Gamma$-invariants of the constant sheaf on $Q$ associated with ${\bf F}$. Therefore we get a spectral sequence \begin{gather} E_2^{rs}=H^r(\Gamma,{\bf F}\otimes_{K_0}H_{rig}^s(Q/{\mathfrak S}^0))\Longrightarrow H_{rig}^{r+s}(Q_{\Gamma}/{\mathfrak S}^0,E).\tag*{$(G)_{0}$}\end{gather}The isomorphism ${\mathbb{R}}\Gamma_{dR}({X})\cong K\otimes_{K_0}{\mathbb{R}}\Gamma_{rig}(Q/{\mathfrak S}^0)$ from \ref{rigs} and hence the obtained isomorphism ${\bf F}\otimes_K {\mathbb{R}}\Gamma_{dR}({X})\cong {\bf F}\otimes_{K_0} {\mathbb{R}}\Gamma_{rig}(Q/{\mathfrak S}^0)$ is $\Gamma$-equivariant and induces the isomorphism ${\mathbb{R}}\Gamma_{dR}({X_{\Gamma}},F)\cong {\mathbb{R}}\Gamma_{rig}(Q_{\Gamma}/{\mathfrak S}^0,E)$ from \ref{rigs} (because the construction \ref{rigs} is local and $Q\to Q_{\Gamma}$ is \'{e}tale). Thus $(G)_{\pi}$ and $(G)_{0}$ are isomorphic. Also here: $\phi$ acts on $(G)_0$, and on its $E_1^{rs}$-term we have $\phi=q^{s}$, by \ref{sympur}. By transport of structure, $\phi$ acts on the spectral sequences $(G)_{\pi}$ and $(C)_{\pi}$, and on their $E^{rs}$-terms it is multiplication with $q^{s}$. Thus $(G)_{\pi}$ and $(C)_{\pi}$ degenerate in $E_2$ for weight reasons and the induced filtrations on the abutment satisfy the described property with respect to $\phi$, hence coincide because this property is characterizing.\\
The degeneration of $(G)_{\pi}$ was proven by another argument in \cite{schn}.   

\addtocounter{satz}{1}{\bf \arabic{section}.\arabic{satz}} The {\it Hodge filtration} $(F_{Hdg}^j)_{j\ge 0}$ on $H_{dR}^{d}(X_{\Gamma},F)$ is the one induced by the stupid filtration $(F\otimes_K\Omega^{\bullet\ge j}_{X_{\Gamma}})_{j\ge 0}$ of $F\otimes_K\Omega^{\bullet}_{X_{\Gamma}}$. Peter Schneider conjectures \cite{schn} that $(F_{Hdg}^j)_{j\ge 0}$ is opposite to $(F_{\Gamma}^r)_{r\ge0}$, i.e. that we have the Hodge-type decomposition $$H_{dR}^{d}(X_{\Gamma},F)=F_{Hdg}^r\bigoplus F_{\Gamma}^{d+1-r}\quad\mbox{ for any }r\in\mathbb{Z}.$$We attempt to understand this conjecture in terms of $p$-adic Hodge theory.

\addtocounter{satz}{1}{\bf \arabic{section}.\arabic{satz}} Let $D$ be a finite dimensional $K$-vector space, endowed with $K$-linear automorphisms $\phi$ and $N$ satisfying $N\phi=q\phi N$, and with a descending exhaustive and separated filtration $(Fil^iD)_{i\in\mathbb{Z}}$. For $i\in\mathbb{Z}$ define the {\it Hodge number} $h_H(D,i)=\dim_K(Fil^iD/Fil^{i+1}D)$. Let $P_0=\quot(W(\overline{k}))$ and denote by $\sigma$ the functorial lifting to $P_0$ of the $q$-power map of the algebraic closure $\overline{k}$ of $k$. For $\alpha=r/s\in\mathbb{Q}$ define the {\it Newton number} $h_N(D,\alpha)=[K:K_0]^{-1}\dim_{K_0}D_{[\alpha]}$ where $D_{[\alpha]}$ is the sub $K_0$-vector space of $D\otimes_{K_0}P_0$ generated by the elements $x$ satisfying $(\phi\otimes\sigma)^sx=q^{r}x$. Then set$$t_N(D)=\sum_{\alpha\in{\mathbb Q}}\alpha h_N(D,\alpha)\quad\quad t_H(D)=\sum_{i\in\mathbb{Z}}i h_H(D,i).$$We say $D$ is {\it weakly admissible} if $t_N(D)=t_H(D)$ and if for all sub $K$-vector spaces $D'\subset D$ stable for $N$ and $\phi$, if we endow them with the induced filtration, we have $t_N(D')\ge t_H(D')$. We say $D$ is {\it ordinary} if it is weakly admissible and if $h_N(D,i)=h_H(D,i)$ for all $i$ in $\mathbb{Z}$, and $h_N(D,\alpha)=0$ for all $\alpha\in\mathbb{Q}-\mathbb{N}$.\\In \cite{peri} 1.2 it is defined another notion of ordinary filtered $(\phi,N)$-module; in particular, $\phi$ in \cite{peri} means a non-iterated Frobenius endomorphism on a $K_0$-lattice. It is easy to see that, given an ordinary filtered $(\phi,N)$-module in the sense of \cite{peri}, by taking the $K$-linear extension of its $\log_pq$-fold iterated Frobenius operator $\phi$ we obtain an ordinary object $D$ as defined above.\\Clearly, if $D$ is ordinary (in our sense) than the slope filtration (appropriately numbered) associated with $\phi$ is opposite to $(Fil^iD)_{i\in\mathbb{Z}}$. Now consider $D=H_{dR}^{d}(X_{\Gamma},F)$ with $Fil^iD=F_{Hdg}^i$ and with $N$, $\phi$ obtained from the Hyodo-Kato isomorphism, as before. By \ref{charfil} the slope filtration (w.r.t. $\phi$) is just $(F_{\Gamma}^r)_{r\ge0}$. So we obtain: 

\begin{pro} If $H_{dR}^{d}(X_{\Gamma},F)$ is ordinary, then the filtrations $(F_{Hdg}^j)_{j\ge 0}$ and $(F_{\Gamma}^r)_{r\ge0}$ on $H_{dR}^{d}(X_{\Gamma},F)$ are opposite.
\end{pro}

\begin{satz}\label{conjhod} $H_{dR}^{d}(X_{\Gamma})$ is ordinary. In particular, the filtrations $(F_{Hdg}^j)_{j\ge 0}$ and $(F_{\Gamma}^r)_{r\ge0}$ on $H_{dR}^{d}(X_{\Gamma})$ are opposite.
\end{satz}

{\sc Proof:} Hyodo defines proper semistable {\it ordinary $A$-schemes}, and explicitly mentions the $A$-schemes ${\mathfrak Q}_{\Gamma}$ as examples (\cite{hyoin} p.544). (-- The ordinarity of ${\mathfrak Q}_{\Gamma}$ also follows from \cite{illord} 1.10 (or \cite{mokr} 3.23) and the fact that all irreducible strata of ${Q}_{\Gamma}$ and their intersections are ordinary $k$-schemes: this last fact follows from our description \arabic{compde1}.\arabic{compde2} and \cite{illord} 1.6. --) In view of our remarks on ordinary filtered $(\phi,N)$-modules our claim now follows from the general fact that for a proper semistable ordinary $A$-scheme $Z$ in the sense of Hyodo, the de Rham cohomology, endowed with its Hodge filtration and with $N$ and $\phi$ coming from the comparison isomorphism with Hyodo-Kato cohomology, is an ordinary filtered $(\phi,N)$-module in the sense of \cite{peri}. (-- We are not aware of a "direct" proof of this last fact, one may however argue like this: From \cite{illast} Therorem 2.7 (which is due to Hyodo) it follows that $\underline{D}_{st}(H^*_{et}(Z\otimes_A{\overline{K}},\mathbb{Q}_p))$ is an ordinary filtered $(\phi,N)$-module in the sense of \cite{peri}. By the now proven Fontaine-Jannsen conjecture, $\underline{D}_{st}(H^*_{et}(Z\otimes_A{\overline{K}},\mathbb{Q}_p))$ is the filtered $(\phi,N)$-module obtained by glueing $H^*_{crys}(Z\otimes_A k/{\mathfrak S}^0)\otimes_{A_0}K_0$ with $H_{dR}^{*}(Z\otimes_A K)$ via the Hyodo-Kato isomorphism.--)

\addtocounter{satz}{1}{\bf \arabic{section}.\arabic{satz}} (i) That $(F_{Hdg}^j)_{j\ge 0}$ and $(F_{\Gamma}^r)_{r\ge0}$ are opposite to each other whenever ${\bf F}$ admits an integral lattice has been shown earlier by Iovita and Spiess \cite{iovspi} by completely different methods, and yet another proof is due to Alon and de Shalit.\\(ii) \label{pdiv} Also for other {\it pure slope} Frobenius structures on the sheaf $E=sp_*F$ (i.e. not necessarily the identity), our arguments show that the filtration $(F_{\Gamma}^r)_{r\ge0}$ on $H_{dR}^{d}(X_{\Gamma},F)$ is the (scalar extended) slope filtration for the $\phi$-action. For example, the Dieudonn\'{e} module of the universal $p$-divisible group over ${\mathfrak Q}$ (base extended to the completion of a maximal unramified extension of $K$; see \cite{rz}) defines such a sheaf $E$ with pure slope Frobenius structure. The above Hodge type decomposition conjecture thus translates into a conjecture on the resulting filtered $(\phi,N)$-module. Via Falting's proof of the $C_{st}$-conjecture with coefficients it then translates into a conjecture on the $p$-adic \'{e}tale cohomology of the relative Tate module of the universal $p$-divisible group. We hope to address this problem in our future work.

\begin{kor} The canonical map $H^d(X_{\Gamma},K)\to H^d(X_{\Gamma},{\cal O}_{X_{\Gamma}})$ is an isomorphism.
\end{kor}

{\sc Proof:} Here $H^d(X_{\Gamma},K)$ is to be understood with respect to the rigid Grothendieck topology on $X_{\Gamma}$. It can be identified with the term $E_2^{d0}$ in the spectral sequence $(C)_{\pi}$ (with ${\bf F}=K$ there). This follows for example from the fact that the tubes $]Q_{\Gamma}^{s}[$ for $s\ge 1$ are disjoint unions of contractible spaces, when viewed as analytic spaces in the sense of Berkovich. Now $(C)_{\pi}$ degenerates in $E_2$ for weight reasons, thus $E_2^{d0}=F_C^d$ in $H_{dR}^d(X_{\Gamma})$. From \ref{conjhod} it follows that $F_C^d=F_{\Gamma}^d\to H_{dR}^d(X_{\Gamma})/F_{Hdg}^1$ is an isomorphism. But $H_{dR}^d(X_{\Gamma})/F_{Hdg}^1$ is $H^d(X_{\Gamma},{\cal O}_{X_{\Gamma}})$, by the degeneration of the Hodge-de Rham spectral sequence.

\section{The monodromy operator}
\label{mosec} 

\addtocounter{satz}{1}{\bf \arabic{section}.\arabic{satz}} Specializing to ${\bf F}=K$, the trivial representation of $\Gamma$, we wish to relate the monodromy operator $N$ on $H_{dR}^{d}(X_{\Gamma})$ to the filtration $(F_{\Gamma}^r)_{r\ge0}$ on $H_{dR}^{d}(X_{\Gamma})$. As remarked in \arabic{cryrig1}.\arabic{cryrig2} we have a canonical isomorphism $H^d_{crys}(Q_{\Gamma}/{\mathfrak S}^0)\otimes_{A_0}K_0\cong H^d_{rig}(Q_{\Gamma}/{\mathfrak S}^0)$. Under this isomorphism, $N$ and $\phi^{\log_pq}$ on $H^d_{crys}(Q_{\Gamma}/{\mathfrak S}^0)\otimes_{A_0}K_0$ (with $N$ and $\phi$ as defined in \cite{hyoka} or \cite{mokr}: this $\phi$ is a {\it non iterated} Frobenius) correspond to $N$ and $\phi$ on $H^d_{rig}(Q_{\Gamma}/{\mathfrak S}^0)$ as we defined it. We keep {\it our} $\phi$. Associated with $N$ is the {\it monodromy filtration} $(M_r)_{r\in\mathbb{Z}}$ on $H^d_{crys}(Q_{\Gamma}/{\mathfrak S}^0)\otimes_{A_0}K_0$, the convolution of the image filtration and the kernel filtration for $N$. More precisely, $$M_r=\sum_i\ke(N^{i+1})\cap \bi(N^{i-r}).$$ On the other hand, $H^d_{crys}(Q_{\Gamma}/{\mathfrak S}^0)\otimes_{A_0}K_0$ is a mixed $F$-isocrystal for the action of $\phi$. Let $(P_r)_{r\in\mathbb{Z}}$ be the corresponding weight filtration on $H^d_{crys}(Q_{\Gamma}/{\mathfrak S}^0)\otimes_{A_0}K_0$, shifted by the number $d$: the uniquely determined $\phi$-stable filtration such that $P_r/P_{r-1}$ is a pure $F$-isocrystal of weight $d+r$. From \ref{charfil} it follows that for all $j\in\mathbb{Z}$ we have $$P_{2j-d}(H_{crys}^d(Q_{\Gamma}/{\mathfrak S}^0)\otimes_{A_0}K_0)=P_{2j-d+1}(H_{crys}^d(Q_{\Gamma}/{\mathfrak S}^0)\otimes_{A_0}K_0),$$ and that after tensoring $\otimes_{K_0}K$ this subspace corresponds to $F^{d-j}_{\Gamma}$ in $H_{dR}^{d}(X_{\Gamma})$. Equivalently, $P_{d-2j}=P_{d-2j+1}$ and this corresponds to $F^j_{\Gamma}$. The {\it monodromy-weight conjecture} for $Q_{\Gamma}$ says $M_r=P_r$ for all $r\in\mathbb{Z}$. It has recently been proven by T. Ito \cite{ito}, and independently, relying on our results in the present paper, by de Shalit \cite{desh}. We obtain:

\begin{satz}\label{moncov} The filtration $(F_{\Gamma}^j)_{j\ge0}$ on $H_{dR}^{d}(X_{\Gamma})$ coincides with the monodromy filtration for $N$ on $H_{dR}^{d}(X_{\Gamma})$:$$F^j_{\Gamma}=\sum_i\ke(N^{i+1})\cap \bi(N^{i-d+2j})=\sum_i\ke(N^{i+1})\cap \bi(N^{i-d+2j-1}).$$ 
\end{satz}

\begin{lem}\label{kebi} Let $d\in\mathbb{N}$ and let $N$ be an endomorphism of an abelian group such that $N^{d+1}=0$. Suppose that for all $j\ge 0$ we have$$\sum_i\ke(N^{i+1})\cap \bi(N^{i-d+2j})=\sum_i\ke(N^{i+1})\cap \bi(N^{i-d+2j-1}).$$Denote these subgroups by $F^j$. If moreover $\ke(N)=F^d$, then also for all $j\ge 0$ $$\ke(N^{d+1-j})=\bi(N^j)=F^j.$$
\end{lem}

{\sc Proof:} First observe that $N^{d+1}=0$ implies $F^j\subset\ke(N^{d+1-j})$ for all $j$. By descending induction on $j$ we prove that this inclusion is an equality. For $j=d$ this is our assumption. Let now $j<d$ and $x\in\ke(N^{d+1-j})$. Then $x'=Nx\in F^{j+1}$ by induction hypothesis. That is, $x'=\sum_ix'_i$ with $x'_i\in \ke(N^{i+1})\cap \bi(N^{i-d+2j+2})$. Write $x'_i=N^{i-d+2j+2}y_i$ and let $z_i=N^{i-d+2j+1}y_i$. Then $z_i\in \ke(N^{i+2})\cap \bi(N^{i-d+2j+1})$, thus $z=\sum_iz_i\in(\sum_i\ke(N^{i+1})\cap \bi(N^{i-d+2j}))=F^j$. By construction also $x-z\in\ke(N)=F^d\subset F^j$. Thus $x=(x-z)+z\in F^j$. This finishes the proof of $\ke(N^{d+1-j})=F^j$. The identity with $\bi(N^j)$ follows easily from this. 

\addtocounter{satz}{1}{\bf \arabic{section}.\arabic{satz}} In $H_{dR}^{d}(X_{\Gamma})$ we define a subspace $C^d$ as follows. If $d$ is odd, we set $C^d=0$. Now let $d$ be even. From $N\phi=q\phi N$ it follows that $\ke(N)$ is stable for $\phi$, and together with \ref{charfil} it follows that $F_{\Gamma}^d$ is the weight-zero-subspace of $\ke(N)$. Let $C^d\subset \ke(N)$ be its $\phi$-stable complement. In the proof of \ref{netcoh} below we will see --- assuming the monodromy-weight conjecture for $Q_{\Gamma}$ --- that $$\dim_{K}(C^d)=1,\quad\quad\quad\quad C^d\cap F_{\Gamma}^{d/2+1}=\emptyset$$ (hence $\phi=q^{d/2}$ on $C^d$). For arbitrary $d$ we define$$\overline{H}^{d}_{dR}(X_{\Gamma})=H^{d}_{dR}(X_{\Gamma})/C^d.$$In particular $\overline{H}^{d}_{dR}(X_{\Gamma})=H^{d}_{dR}(X_{\Gamma})$ if $d$ is odd. The operators $N$ and $\phi$ on $H^{d}_{dR}(X_{\Gamma})$ induce operators $\overline{N}$ and $\overline{\phi}$ on $\overline{H}^{d}_{dR}(X_{\Gamma})$. We {\it define}: the filtration $({\overline F}_{\Gamma}^r)_{r\ge 0}$ on $\overline{H}^{d}_{dR}(X_{\Gamma})$ is the image of the filtration $(F_{\Gamma}^r)_{r\ge0}$ on $H^{d}_{dR}(X_{\Gamma})$; the filtration $({\overline F}_{Hdg}^j)_{j\ge 0}$ on $\overline{H}^{d}_{dR}(X_{\Gamma})$ is the image of the filtration $(F_{Hdg}^j)_{j\ge 0}$ on $H^{d}_{dR}(X_{\Gamma})$.

\begin{satz}\label{netcoh} (a) The filtrations $({\overline F}_{\Gamma}^r)_{r\ge0}$ and $({\overline F}_{Hdg}^j)_{j\ge 0}$ are opposite.\\(b) The filtration $({\overline F}_{\Gamma}^r)_{r\ge0}$ is stable for $\overline{\phi}$; we have $\overline{\phi}=q^{d-r}$ on ${\overline F}_{\Gamma}^r/\overline {F}^{r+1}_{\Gamma}$.\\(c) The filtration $({\overline F}_{\Gamma}^r)_{r\ge0}$ coincides with both the kernel and the image filtration for $\overline{N}$: for all $r$ we have $$\overline {F}^r_{\Gamma}=\ke(\overline{N}^{d+1-r})=\bi(\overline{N}^{r}).$$  
\end{satz}

{\sc Proof:} Assertions (a) and (b) follow from \ref{conjhod} and \ref{charfil}. For (c) we need the computation \cite{ss} p.93 of the dimensions of the graded pieces for the filtration $({F}_{\Gamma}^j)_{j\ge0}$. Namely, for odd $d$ and $d\ge j\ge0$ we have$$\dim_K(F_{\Gamma}^j/{F}_{\Gamma}^{j+1})=\dim_K({H}_{dR}^{d}(X_{\Gamma})/{F}_{\Gamma}^1).$$If $d$ is even, $d=2t$, the same holds for all $d\ge j\ge 0$, $j\ne t$, and moreover $$\dim_K(F_{\Gamma}^{t}/{F}_{\Gamma}^{t+1})=\dim_K({H}_{dR}^{d}(X_{\Gamma})/{F}_{\Gamma}^1)+1.$$From \ref{moncov} we easily deduce that the iterates of $N$ induce surjective maps\begin{gather}N^k:F_{\Gamma}^j/{F}_{\Gamma}^{j+1}\longrightarrow F_{\Gamma}^{j+k}/{F}_{\Gamma}^{j+k+1}\tag{$1$}\end{gather}for all $j,k$ with $d-k\ge j\ge0$. Now consider first the case where $d$ is odd. Then all maps $(1)$ must be bijective, for dimension reasons. By descending induction on $e$ we prove $F_{\Gamma}^e=\bi(N^e)$ for all $d\ge e\ge 1$. By \ref{moncov} we have $$F^e_{\Gamma}=\sum_i\ke(N^{i+1})\cap \bi(N^{i-d+2e}).$$For the summation index $i=d-e$ we have $\ke(N^{i+1})\cap \bi(N^{i-d+2e})=\bi(N^e)$ since $N^{d+1}=0$. For summation indices $i>d-e$ we have $\ke(N^{i+1})\cap \bi(N^{i-d+2e})\subset \bi(N^{i-d+2e})\subset \bi(N^e)$. For summation indices $0\le i<d-e$ we consider the isomorphism $(1)$ with $k=i+1$ and $j=e$: it tells us $\ke(N^{i+1})\cap \bi(N^{i-d+2e})\subset F^{e+1}$, but by induction hypothesis $F^{e+1}=\bi(N^{e+1})\subset \bi(N^e)$; the induction is finished. That $F_{\Gamma}^e=\ke({N}^{d+1-e})$ for all $e$ follows easily from this and the identities$$F^e_{\Gamma}=\sum_i\ke(N^{i+1})\cap \bi(N^{i-d+2e-1})$$from \ref{moncov}. Now let $d$ be even, $d=2t$. Here we first observe \begin{gather}\ke(N)\ne F_{\Gamma}^d\tag{$2$}.\end{gather}Indeed, if we had $\ke(N)=F_{\Gamma}^d$ we could apply \ref{kebi} (we know $N^{d+1}=0$ from \cite{mokr}) which in particular would give us ${F}_{\Gamma}^1=\ke(N^d)=\ke(N^{2t})$, $F_{\Gamma}^{t}=\bi(N^t)$ and ${F}_{\Gamma}^{t+1}=\ke(N^t)$, implying that $N^t$ induces an isomorphism ${H}_{dR}^{d}(X_{\Gamma})/{F}_{\Gamma}^1\cong F_{\Gamma}^{t}/{F}_{\Gamma}^{t+1}$, contradicting our dimension estimates. However, for $d$ odd one shows \begin{gather}F_{\Gamma}^e=\bi(N^e)\quad\quad \mbox{ for all } d\ge e\ge t+1. \tag{$3$}\end{gather}Together with the identities$$F^e_{\Gamma}=\sum_i\ke(N^{i+1})\cap \bi(N^{i-d+2e-1})$$from \ref{moncov} one derives $\ke(N)\cap\bi(N)\subset F_{\Gamma}^e$ for all $d\ge e\ge t+1$ by an easy induction on $e$ (beginning with $e=t+1$); in particular $F_{\Gamma}^d=\ke(N)\cap\bi(N)$. By definition, $C^d\subset\ke(N)$ is the uniquely determined $\phi$-stable complement of the weight-zero-part $F_{\Gamma}^d$ of $\ke(N)$. From $F_{\Gamma}^d=\ke(N)\cap\bi(N)$ we get $C^d\cap\bi(N)=\emptyset$ which implies $\ke(\overline{N})=\ke(N)\mod C^d$ for the operator $\overline{N}=(N\mod C^d)$ on $\overline{H}^{d}_{dR}(X_{\Gamma})$. Therefore the identities from \ref{moncov} pass to the identities $$\overline {F}^j_{\Gamma}=\sum_i\ke(\overline{N}^{i+1})\cap \bi(\overline{N}^{i-d+2j})=\sum_i\ke(\overline{N}^{i+1})\cap \bi(\overline{N}^{i-d+2j-1}).$$ Since $F^{t+1}_{\Gamma}\subset \bi(N)$ and $\ke(N)\subset F^t_{\Gamma}$ (again from \ref{moncov}) we also get $C^d\cap F^{t+1}_{\Gamma}=\emptyset$ and $C^d\subset F_{\Gamma}^t$. On the other hand we know $C^d\ne0$ from (2). Therefore we obtain the estimates$$\dim_K({\overline F}_{\Gamma}^j/{\overline F}_{\Gamma}^{j+1})=\dim_K(F_{\Gamma}^j/{F}_{\Gamma}^{j+1})=\dim_K({H}_{dR}^{d}(X_{\Gamma})/{F}_{\Gamma}^1) \mbox{ for all }d\ge j\ge 0, j\ne t $$$$\dim_K({\overline F}_{\Gamma}^t/{\overline F}_{\Gamma}^{t+1})\le\dim_K(F_{\Gamma}^t/{F}_{\Gamma}^{t+1})-1=\dim_K({H}_{dR}^{d}(X_{\Gamma})/{F}_{\Gamma}^1).$$Now the same proof as in the case where $d$ is odd gives assertion (c) also in the case where $d$ is even. We are done.     

  %@@

\end{document}